\documentclass{article}
\usepackage[dvips]{graphicx}
\usepackage{subfigure}
\usepackage{pstricks, amssymb, xr, multind, multicol}
\def\hfl#1{\smash{\mathop{\hbox to 10mm{\rightarrowfill}}\limits^{\textstyle
#1}}}

\makeindex{NN}
\makeindex{TT}
\newcommand{\indexT}[1]{#1\index{TT}{#1}}

\newtheorem{proposition}[equation]{Proposition}
\newtheorem{corollary}[equation]{Corollary}
\newtheorem{theorem}[equation]{Theorem}
\newtheorem{exa}[equation]{Example}
\newtheorem{ex}[equation]{Exercise}
\newtheorem{s-ex}[equation]{Side-exercise}
\newtheorem{exas}[equation]{Examples}
\newtheorem{lemma}[equation]{Lemma}
\newtheorem{sublemma}[equation]{Sublemma}
\newtheorem{remar}[equation]{Remark}
\newtheorem{remars}[equation]{Remarks}
\newtheorem{nota}[equation]{Notation}
\newtheorem{sremar}[equation]{Side-remark}
\newtheorem{definitio}[equation]{Definition}

\newenvironment{notation}{\begin{nota} \rm }{\end{nota}}
\newenvironment{remark}{\begin{remar} \rm }{\end{remar}}

\newenvironment{example}{\begin{exa} \rm }{\end{exa}}
\newenvironment{definition}{\begin{definitio} \rm }{\end{definitio}}

\newcommand{\CA}{{\cal A}}

\newcommand{\CI}{{\cal I}}

\newcommand{\CL}{{\cal L}}

\newcommand{\CT}{{\cal T}}

\newcommand{\ZZ}{\mathbb{Z}}
\newcommand{\RR}{\mathbb{R}}
\newcommand{\QQ}{\mathbb{Q}}
\newcommand{\CC}{\mathbb{C}}
\newcommand{\NN}{\mathbb{N}}
\newcommand{\bp}{\noindent {\sc Proof: }}
\newcommand{\eop}{\nopagebreak
			\hspace*{\fill}{$\diamond$}
			\medskip}
\newcommand{\const}{\cite{lesconst}}
\newcommand{\refconst}[1]{\cite[#1]{lesconst}}
\newcommand{\inconstdefformfund}{Definition~\ref{defformfund} in \const}
\newcommand{\inconstlemdepboun}{Lemma~\ref{lemdepboun} in \const}
\newcommand{\inconstlemexisfunap}{Lemma~\ref{lemexisfunap} in \const}
\newcommand{\inconstlemhctwo}{Lemma~\ref{lemhctwo} in \const}
\newcommand{\inconstlemnolosseta}{Lemma~\ref{lemnolosseta} in \const}
\newcommand{\inconstlemtranspform}{Lemma~\ref{lemtranspform} in \const}

\newcommand{\inconstproppont}{Proposition~\ref{proppont} in \const}
\newcommand{\inconstproporrev}{Proposition~\ref{proporrev} in \const}
\newcommand{\inconstpropxiun}{Proposition~\ref{propxiun} in \const}
\newcommand{\inconstpropzad}{Proposition~\ref{propzad} in \const}
\newcommand{\inconstpropkeytriv}{Proposition~\ref{propkeytriv} in \const}
\newcommand{\inconstproppontproptrois}{Propositions~\ref{proppont} and \ref{proptrois} in \const}
\newcommand{\inconstproppontproptroispropkeytriv}{Propositions~\ref{proppont}, \ref{proptrois},  and \ref{propkeytriv} in \const}
\newcommand{\inconstsubapseckkt}{Subsection~\ref{subapseckkt} in \const}
\newcommand{\inconstsubpont}{Subsection~\ref{subpont} in \const}
\newcommand{\inconstsubsketchpkkt}{Subsection~\ref{subsketchpkkt} in \const}

\newcommand{\rconstsecsta}{\refconst{Section~\ref{secsta}}}
\newcommand{\rconstsubfundform}{\cite[Subsection~\ref{subfundform}]{lesconst}}
\newcommand{\rconstthkktfra}{\refconst{Theorem~\ref{thkktfra}}}
\newcommand{\tata}{\begin{pspicture}[0.2](0,0)(.5,.4)
\pscircle(0.25,0.2){.25}
\psline{*-*}(0.05,.2)(.45,.2)
\end{pspicture}}

\begin{document}
\thispagestyle{empty}
\begin{center}
\vskip4cm

\huge{Splitting formulae for the Kontsevich-Kuperberg-Thurston invariant of rational homology $3$-spheres}

\vskip1cm

\large{Christine Lescop}

\vskip1cm

\begin{large}
Pr\'epublication de l'Institut Fourier n$^{o}$ 656 (2004) \\
http://www-fourier.ujf-grenoble.fr/prepublications.html \\
\end{large}
\normalsize(Second version: Minor modifications in the abstract and in the introduction) 
\end{center}

\vskip 4cm

\begin{abstract}
M.~Kontsevich proposed a topological construction for an invariant $Z$ of rational homology $3$-spheres using configuration space integrals. 
G.~Kuperberg and D.~Thurston proved that $Z$ is a universal real finite type invariant for integral homology spheres in the sense of Ohtsuki, Habiro and Goussarov.

We discuss the behaviour of $Z$ under rational homology handlebodies replacements.
The explicit formulae that we present generalize a sum formula
obtained by the author for the Casson-Walker invariant in 1994. They allow us
to identify the degree one term of $Z$ with the Walker invariant for rational homology spheres.

\vskip1cm

\noindent {\bf Keywords:} 3-manifolds, configuration space integrals, homology spheres, finite type invariants, 
Jacobi diagrams, clovers, claspers, Casson-Walker invariant\\ 
{\bf A.M.S. subject classification:} 57M27 57N10 55R80 57R20
\end{abstract}

\vskip4cm

\noindent Institut Fourier \\
Laboratoire de math\'ematiques de l'Universit\'e Joseph Fourier, UMR 5582 du CNRS \\
B.P.74 \\
38402 Saint-Martin-d'H\`eres cedex (France)\\
mail: lescop@ujf-grenoble.fr

\baselineskip.5cm
\newpage
\tableofcontents

\newpage
\section{Introduction}

In 1995 in \cite{oht}, Tomotada Ohtsuki introduced a notion of finite type invariants for homology $3$-spheres (that are compact oriented 3-manifolds with the same homology with integral coefficients as the standard $3$-sphere $S^3$), following the model of the notion of Vassiliev invariants
for links in the ambient space $\RR^3$. He defined a filtration of the real vector space freely generated by homology $3$-spheres and began the study of the
associated graded space. In \cite{le}, Thang Le finished identifying this graded space
to an algebra of Jacobi diagrams called $\CA(\emptyset)$ whose definition is recalled in Subsection~\ref{subjac}. To do this, Le proved that the
LMO invariant of $3$-manifolds $Z_{LMO}$
\index{NN}{ZLMO@$Z_{LMO}$} that he constructed with the help of J.~Murakami and Ohtsuki
in \cite{lmo} induces an isomorphism from the Ohtsuki graded space to $\CA(\emptyset)$. In \cite{kt}, following Witten, Axelrod, Singer, Kontsevich, Bott and Cattaneo, Greg Kuperberg and Dylan Thurston constructed another (possibly equal) invariant $Z_{KKT}$ \index{NN}{ZKKT@$Z_{KKT}$} of rational homology $3$-spheres valued in $\CA(\emptyset)$, and they proved that $Z_{KKT}$ also induces the already mentioned Le isomorphism. 
All real-valued finite type invariants in the Ohtsuki sense factor through $Z=Z_{KKT}$ \index{NN}{Z@$Z$}(or $Z_{LMO}$).
Therefore $Z_{KKT}$ and $Z_{LMO}$ are called universal finite type invariants
of homology $3$-spheres. They play the same role as the Kontsevich integral does play in the theory of Vassiliev link invariants.

In this article,
%%% BLOC monograph 
 we prove explicit formulae on the behaviour of $Z_{KKT}$ under Lagrangian-preserving rational homology handlebodies replacements. The precise statement is given in Theorem~\ref{mainth} after the needed definitions. This behaviour had been
observed in the case of Torelli replacements by Kuperberg and Thurston in \cite{kt}, it is the key point in their proof of universality for the invariant $Z_{KKT}$.

The obtained formulae described below generalize the formulae obtained in \cite{les} for the Casson-Walker invariant. They enlight the relationships between finite type invariants, Jacobi diagrams, intersection forms and linking forms. 
%and the information encoded in the weight systems of finite type invariants.
They also allow us to identify the degree one part of $Z_{KKT}$
with the Walker invariant for any rational homology sphere in Theorem~\ref{thcompwal}.

In the case of integral homology spheres, it is proved in \cite{al} that the splitting formulae obtained in this article
%%% BLOC monograph 
 follow from
the Kuperberg-Thurston formulae for Torelli replacements.

The detailed proofs of the formulae in the general case are given in Sections~\ref{secskp} to \ref{secnormap}. Their sketch, that is given in Section~\ref{secskp}, is the 
now-standard sketch in this kind of proofs. But filling in the details in the general case was surprisingly complicated to me.
The detailed proofs are given here with full generality, they substantially simplify
in the case of Torelli replacements. Since the proofs heavily rely on the 
Kontsevich-Kuperberg-Thurston construction, this construction 
has been recalled in 
%%% BLOC 
%is recalled in Section~\ref{secsta}, 
%\cite[Section~\ref{secsta}]{lesconst} 
\rconstsecsta ~and all the precise statements that are needed in our proof
are given and proved in 
%Sections~\ref{secsta} to \ref{seccomp}. 
\const.
%%% BLOC
%No originality is claimed
%%in Sections~\ref{secsta} to \ref{seccomp} that are 
%in \const that is
%just a painful work of writing up the details
%used in 
%%the following more original sections. 
%%% BLOC
%this article.
%However 
%%these sections 
%%% BLOC
% \const
%may be useful as a detailed partial review of \cite{kt}.
I hope that the technical work contained here will help finding other properties
for the invariant $Z_{KKT}$.

I thank Dylan Thurston for very useful and pleasant conversations. 
 
%Now, let us be more specific and give the exact statement of the main theorem.

%%% BLOC
%\section{Statement of the main original result}
\section{Statement of the main result}
\setcounter{equation}{0}
\label{secstamain}

\subsection{Jacobi diagrams}
\label{subjac}
Here, a {\em  \indexT{Jacobi diagram}\/} $\Gamma$ is a trivalent graph $\Gamma$ without simple loop like $\begin{pspicture}[.2](0,0)(.6,.4)
\psline{-*}(0.05,.2)(.25,.2)
\pscurve{-}(.25,.2)(.4,.05)(.55,.2)(.4,.35)(.25,.2)
\end{pspicture}$. The set of vertices of such a $\Gamma$ will be denoted by $V(\Gamma)$,
\index{NN}{VGamma@$V(\Gamma)$}
its set of edges will be denoted by $E(\Gamma)$.
\index{NN}{EGamma@$E(\Gamma)$}
A {\em \indexT{half-edge}\/} $c$ of $\Gamma$ is an element of
$$H(\Gamma)=\{c=(v(c);e(c)) | v(c) \in V(\Gamma); e(c) \in E(\Gamma);v(c) \in e(c)\}.$$
\index{NN}{HGamma@$H(\Gamma)$}
An {\em automorphism\/} of $\Gamma$ 
\index{TT}{Jacobi diagram!automorphism of} 
is a permutation $b$ of $H(\Gamma)$
such that for any $c,c^{\prime} \in H(\Gamma)$,
$$v(c)=v(c^{\prime}) \Longrightarrow v(b(c))=v(b(c^{\prime}))\;\;\mbox{and}\;\;e(c)=e(c^{\prime}) \Longrightarrow e(b(c))=e(b(c^{\prime})).$$
The number of automorphisms of $\Gamma$ will be
denoted by $\sharp \mbox{Aut}(\Gamma)$. 
\index{NN}{AutGamma@$\sharp \mbox{Aut}(\Gamma)$}
For example, $ \sharp \mbox{Aut}(\tata)=12$.
{\em An orientation\/} of a vertex of such a diagram $\Gamma$ 
 is a cyclic order of the three
half-edges that meet at that vertex.
\index{TT}{Jacobi diagram!orientation of}
A Jacobi diagram $\Gamma$ is {\em oriented\/} if all its vertices are oriented (equipped with an orientation).
The {\em degree} of such a diagram is 
half the number of its vertices. 

Let $\CA_n(\emptyset)$ 
\index{NN}{An@$\CA_n(\emptyset)$}
denote the real vector space generated by the degree $n$ oriented Jacobi diagrams, quotiented out by the following relations AS and IHX:

$$ {\rm AS :}  \begin{pspicture}[.2](0,-.2)(.8,1)
\psset{xunit=.7cm,yunit=.7cm}
\psarc[linewidth=.5pt](.5,.5){.2}{-70}{15}
\psarc[linewidth=.5pt](.5,.5){.2}{70}{110}
\psarc[linewidth=.5pt]{->}(.5,.5){.2}{165}{250}
\psline{*-}(.5,.5)(.5,0)
\psline{-}(.1,.9)(.5,.5)
\psline{-}(.9,.9)(.5,.5)
\end{pspicture}
+
\begin{pspicture}[.2](0,-.2)(.8,1)
\psset{xunit=.7cm,yunit=.7cm}
\pscurve{-}(.9,.9)(.3,.7)(.5,.5)
\pscurve[border=2pt]{-}(.1,.9)(.7,.7)(.5,.5)
\psline{*-}(.5,.5)(.5,0)
\end{pspicture}=0,\;\;\mbox{and IHX :} 
\begin{pspicture}[.2](0,-.2)(.8,1)
\psset{xunit=.7cm,yunit=.7cm}
\psline{-*}(.1,1)(.35,.2)
\psline{*-}(.5,.5)(.5,1)
\psline{-}(.75,0)(.5,.5)
\psline{-}(.25,0)(.5,.5)
\end{pspicture}
+
\begin{pspicture}[.2](0,-.2)(.8,1)
\psset{xunit=.7cm,yunit=.7cm}
\psline{*-}(.5,.6)(.5,1)
\psline{-}(.8,0)(.5,.6)
\psline{-}(.2,0)(.5,.6)
\pscurve[border=2pt]{-*}(.1,1)(.3,.3)(.7,.2)
\end{pspicture}
+
\begin{pspicture}[.2](0,-.2)(.8,1)
\psset{xunit=.7cm,yunit=.7cm}
\psline{*-}(.5,.35)(.5,1)
\psline{-}(.75,0)(.5,.35)
\psline{-}(.25,0)(.5,.35)
\pscurve[border=2pt]{-*}(.1,1)(.2,.75)(.7,.75)(.5,.85)
\end{pspicture}
=0. 
$$
\index{NN}{AS}
\index{NN}{IHX}

Each of these relations relate diagrams which can be represented by planar immersions that are identical outside the part of them represented in the pictures. Here, the orientation of vertices is induced by the counterclockwise order of the half-edges. For example, 
AS identifies the sum of two diagrams which only differ by the orientation
at one vertex to zero.
$\CA_0(\emptyset)$ is equal to $\RR$ generated by the empty diagram.

\subsection{The universal finite type invariant $Z$}

Let $\Lambda$ be $\ZZ$, $\ZZ/2\ZZ$ or $\QQ$.
A {\em $\Lambda$-sphere\/} is a  compact oriented 3-manifold $M$ such that $H_{\ast}(M;\Lambda)=H_{\ast}(S^3;\Lambda)$. A \index{TT}{Zsphere@$\ZZ$-sphere} $\ZZ$-sphere is also called a \indexT{homology sphere} while a \indexT{rational homology sphere} is a \index{TT}{Qsphere@$\QQ$-sphere} $\QQ$-sphere.
Following Witten, Axelrod, Singer, Kontsevich, Bott and Cattaneo, Greg Kuperberg and Dylan Thurston constructed invariants $Z_n=(Z_{KKT})_n$ of oriented $\QQ$-spheres valued in $\CA_n(\emptyset)$ and they proved that these invariants have the following property:

\begin{theorem}[Kuperberg-Thurston \cite{kt}]
\label{thktone}
An  invariant $\nu$ of\/ $\ZZ$-spheres valued in a real vector space $X$ is
of degree $\leq n$ if and only if there exist
linear maps $$\phi_k(\nu): \CA_k(\emptyset) \longrightarrow X,$$ for any $k \leq n$, such that $$\nu=\sum_{k=0}^n
\phi_k(\nu) \circ Z_k.$$
\end{theorem}

A {\em real \indexT{finite type invariant}\/} of $\ZZ$-spheres is a topological invariant of 
$\ZZ$-spheres valued in a real vector space $X$ which is of degree less than
some natural integer $n$.
The  Kontsevich-Kuperberg-Thurston construction is recalled in \rconstsecsta. In this article,
%%% BLOC monograph, 
Theorem~\ref{thktone} is used as a definition of degree $\leq n$
real-valued invariants of $\ZZ$-spheres. 

A degree $\leq n$ invariant $\nu$ is of {\em degree\/} $n$ if $\phi_n(\nu) \neq 0$. In this case,
$\phi_n(\nu)$ is the {\em \indexT{weight system}\/} of $\nu$ and is denoted by $W_{\nu}$.

\begin{remark} The above definition coincides with the Ohtsuki definition of real finite type invariants \cite{oht}. The Ohtsuki degree (that is always a multiple of $3$) is three times the above degree. We shall not discuss the concept of finite-type invariants any further here. See \cite{oht,ggp,hab,al} and references therein.
\end{remark}

\subsection{Generalized clovers}

Unless otherwise mentioned, manifolds are compact and oriented.
Boundaries are oriented with the outward normal first convention.
A {\em genus $g$ \indexT{rational homology handlebody} or \index{TT}{Qhandlebody@$\QQ$-handlebody} $\QQ$-handlebody}\/ (resp. a {\em genus $g$ \indexT{integral homology handlebody} or \index{TT}{Zhandlebody@$\ZZ$-handlebody} $\ZZ$-handlebody\/}) is an (oriented, compact) 3-manifold $A$ with
the same homology with rational (resp. integral) coefficients as the standard (solid) handlebody $H_g$ below. \index{NN}{Hg@$H_g$}
$$H_g = \begin{pspicture}[.4](0,0)(4.5,.95)
\psset{xunit=.5cm,yunit=.5cm}  
\psecurve{-}(5.7,1.3)(5.2,1.3)(3.9,1.8)(2.6,1.3)(1.3,1.8)(.1,1)(1.3,.1)(2.6,.7) (3.9,.1)(5.2,.7)(5.7,.7)
\pscurve{-}(.8,1.2)(1,.9)(1.3,.8)(1.6,.9)(1.8,1.2)
\pscurve{-}(1,.9)(1.3,1.2)(1.6,.9)
\rput[r](.9,-.2){$a_1$}
\psecurve{->}(1.6,.4)(1.3,.8)(1.05,.4)(1.3,.1)
\psecurve{-}(1.3,.8)(1.05,.4)(1.3,.1)(1.6,.4)
\psecurve[linestyle=dashed,dash=3pt 2pt](1,.4)(1.3,.8)(1.55,.4)(1.3,.1)(1,.4)
\pscurve{-}(3.4,1.2)(3.6,.9)(3.9,.8)(4.2,.9)(4.4,1.2)
\pscurve{-}(3.6,.9)(3.9,1.2)(4.2,.9)
\rput[r](3.5,-.2){$a_2$}
\psecurve{->}(4.2,.4)(3.9,.8)(3.65,.4)(3.9,.1)
\psecurve{-}(3.9,.8)(3.65,.4)(3.9,.1)(4.2,.4)
\psecurve[linestyle=dashed,dash=3pt 2pt](3.6,.4)(3.9,.8)(4.15,.4)(3.9,.1)(3.6,.4)
\rput(5.8,1.3){\dots}
\rput(5.8,.7){\dots}
\psecurve{-}(5.8,1.3)(6.4,1.3)(7.7,1.8)(8.9,1)(7.7,.1)(6.4,.7)(5.8,.7) 
\pscurve{-}(7.2,1.2)(7.4,.9)(7.7,.8)(8,.9)(8.2,1.2)
\pscurve{-}(7.4,.9)(7.7,1.2)(8,.9)
\rput[l](6.8,-.2){$a_g$}
\psecurve{->}(8,.4)(7.7,.8)(7.45,.4)(7.7,.1)
\psecurve{-}(7.7,.8)(7.45,.4)(7.7,.1)(8,.4)
\psecurve[linestyle=dashed,dash=3pt 2pt](7.4,.4)(7.7,.8)(7.95,.4)(7.7,.1)(7.4,.4)
\end{pspicture}$$
Note that the boundary of such a $\QQ$-handlebody $A$ is homeomorphic to the
boundary $(\partial H_g =\Sigma_g)$ of $H_g$. The intersection form on a surface ${\Sigma}$
is denoted by $\langle,\rangle_{\Sigma}$. \index{NN}{int@$\langle,\rangle_{\Sigma}$}
For a (compact, oriented) 3-manifold $A$ with boundary $\partial A$, 
${\cal L}_A$ \index{NN}{LA@${\cal L}_A$}denotes the kernel of the map induced by the inclusion:
$$ H_1(\partial A;\QQ) \longrightarrow H_1( A;\QQ).$$
It is a Lagrangian of $(H_1(\partial A;\QQ),\langle,\rangle_{\partial A})$, we call it the {\em \indexT{Lagrangian}}\/ of $A$. 

A {\em \indexT{rational generalised clover}\/} is a 4-tuple  %$$D=(M;k;(A^i)_{i=1,\dots,k};(B^i)_{i=1,\dots,k};(\phi_i)_{i=1,\dots,k})$$
$$D=(M;k;(A^i)_{i=1,\dots,k};(B^i)_{i=1,\dots,k})$$
where
\begin{enumerate}
\item $M$ is a rational homology sphere,\\
\item for any $i=1,2, \dots k$, $A^i$ and $B^i$ are $\QQ$-handlebodies
whose boundaries are identified by implicit diffeomorphisms (we shall write $\partial B^i=\partial A^i$) so that $\CL_{B^i}=\CL_{A^i}$,
\item the disjoint union of the $A^i$ is embedded in $M$. We shall write
$$\sqcup_{i=1}^k A^i \subset M.$$
\end{enumerate}
The integral number $k$ is called the {\em degree\/} of $D$. \index{TT}{rational
generalised clover!degree of} 
Such a rational generalised clover $D$ is an {\em \indexT{integral generalised clover}\/}
if furthermore
$M$ is an integral homology sphere, and if $B_1$, $B_2$, \dots, $B_k$  are integral homology handlebodies.

For such a rational generalised clover $D$, if $J$ is a subset of $\{1, \dots, k\}$,
$M_J(D)$ denotes the rational homology sphere obtained by replacing $A^i$ by $B^i$ for every element $i$ of $J$.
\index{NN}{MJD@$M_J(D)$}
$$M_J(D)=\left(M \setminus \cup_{i \in J} \mbox{Int}(A^i) \right) \cup_{ (\cup_{i \in J} \partial A^i)} \left( \cup_{i \in J} B^i \right)$$

If $I$ is a topological invariant of integral (resp. rational) homology spheres
valued in an abelian group, and if $D$ is an integral (resp. rational) generalised
clover, then we define $I(D)$ as $I(D)= \sum_{J \subseteq \{1, \dots, k\}}(-1)^{\sharp J}I(M_J(D))$.

\begin{remark} 
The terminology {\em generalised clover\/} may not be a very happy one. I use it for the following reasons. 
The {\em generalised clovers\/} generalise the \cite{ggp} clovers. 
In \cite{hab}, Habiro independently developed a {\em clasper calculus\/} that encloses the {\em clover calculus\/} and also allows for more general modifications. 
In the Habiro terminology, {\em clovers\/} are called {\em allowable graph claspers.\/} 
I feel that the terminology {\em generalised clasper\/} cannot be used for something that does not generalise all the Habiro claspers, and I do not feel
like saying {\em generalised allowable graph claspers.\/}
\end{remark}

\subsection{Generalised clovers and Jacobi diagrams}

Let $\Gamma$ be an oriented degree $n$ Jacobi diagram.
Let $V(\Gamma)$ \index{NN}{VGamma@$V(\Gamma)$} and $E(\Gamma)$ \index{NN}{EGamma@$E(\Gamma)$} denote the set of vertices of $\Gamma$ and the set of edges of $\Gamma$, respectively.
%A {\em half-edge\/} $c$ of $\Gamma$ is a pair $c=(v(c);e(c))$ 
%where $v(c) \in V(\Gamma)$, $e(c) \in E(\Gamma)$ and $v(c)$ belongs to $e(c)$.
The set of half-edges of $\Gamma$ is denoted by $H(\Gamma)$ \index{NN}{HGamma@$H(\Gamma)$} and its two natural  projections onto $V(\Gamma)$ and $E(\Gamma)$ are denoted by $v$ and $e$, respectively.

Let $D=(M;2n;(A^i)_{i=1,\dots,2n};(B^i)_{i=1,\dots,2n})$
be a rational generalised clover.
Let $\sigma: V(\Gamma) \longrightarrow \{1,2,\dots,2n\}$ be a bijection.
Let us define the {\em linking number $\ell(D;\Gamma;\sigma)$ 
\index{NN}{lDGammasigma@$\ell(D;\Gamma;\sigma)$} of $D$ with respect to
$\Gamma$ and $\sigma$.\/}

%The boundary of an oriented manifold is always oriented 
%with the outward normal first convention.
The Mayer-Vietoris boundary map 
$$\partial_{i,MV}:H_2(A^i \cup_{\partial A^i} -B^i) \longrightarrow \CL_{A^i}$$
that maps the homology class of an oriented surface to the oriented boundary of its intersection with $A^i$ is an isomorphism.
This isomorphism carries the intersection form of the closed $3$-manifold $(A^i \cup_{\partial A^i} -B^i)$ on $\otimes^3 H_2(A^i \cup_{\partial A^i} -B^i)$ to a linear form $\CI(A^i,B^i)$ \index{NN}{IAB@$\CI(A^i,B^i)$} on
$\otimes_{j=1}^3 \CL^{(j)}_{A^i}$ which is antisymmetric with respect to the permutation of two factors, where $\CL^{(j)}_{A^i}$ denote the $j^{\mbox{\small th}}$ copy of $\CL_{A^i}$.
The linear form $\CI(A^i,B^i)$ may be seen canonically as an element of
$\otimes_{j=1}^3 \left(\CL^{(j)}_{A^i}\right)^{\ast}$ where $\left(\CL^{(j)}_{A^i}\right)^{\ast}$ denotes the dual $\mbox{Hom}(\CL^{(j)}_{A^i};\QQ)$ of $\CL^{(j)}_{A^i}$.

For each vertex $w$, number the three half-edges that contain $w$
with a bijection 
$$b(w): v^{-1}(w) \longrightarrow \{1,2,3\}$$ 
that induces the given cyclic order of these half-edges.

Let $c$ be a half-edge.
Assign it the space $$X(c)=\left(\CL^{(b(v(c))(c)}_{A^{\sigma(v(c))}}\right)^{\ast}.$$
The linear form $\CI(A^i,B^i)$ belongs to 
$\otimes_{c \in H(\Gamma); \sigma(v(c))=i}X(c)$. The tensor product of all
the $\CI(A^i,B^i)$, for $i=1,2,\dots,2n$, belongs to
$$\otimes_{c \in H(\Gamma)}X(c).$$

For $\{i,j\} \subseteq  \{1,2,\dots,2n\}$, the linking number in $M$ induces a bilinear form on $H_1(A^i;\QQ) \times H_1(A^j;\QQ)$,
where $H_1(A^i)$ is canonically isomorphic to $\frac{H_1(\partial A^i;\QQ)}{\CL_{A^i}}$.
Furthermore, the intersection form $\langle , \rangle_{\partial A^i}$ induces the map
$$\begin{array}{llll} \langle,.\rangle:& H_1(\partial A^i;\QQ) &\longrightarrow & \CL_{A^i}^{\ast}\\
& x & \mapsto &\langle .,x \rangle\end{array}$$
that in turn induces an isomorphism from 
$\frac{H_1(\partial A^i;\QQ)}{\CL_{A^i}}$ to $\CL_{A^i}^{\ast}$.

Thus, for each edge $f \in E(\Gamma)$ made of two half-edges $c$ and $d$, (so that $e^{-1}(f)=\{c,d\}$) the linking number
yields a contraction
$$ \ell_f: X(c) \otimes X(d) \longrightarrow \QQ.$$
Applying all these contractions to our big tensor maps this tensor to the {\em linking number $\ell(D;\Gamma;\sigma)$ of $D$ with respect to
$\Gamma$ and $\sigma$.\/} 

Finally, we define the {\em linking number $\ell(D;\Gamma)$ \index{NN}{lDGamma@$\ell(D;\Gamma)$} of $D$ with respect to
$\Gamma$\/} as the sum running over all the bijections $\sigma$ from $V(\Gamma)$ to $\{1,2,\dots,2n\}$ of the $\ell(D;\Gamma;\sigma)$.
Note that the product $\ell(D;\Gamma)[\Gamma]$ does not depend on the vertex-orientation of $\Gamma$.

\subsection{Statement of the theorem}

The main theorem of this article
%%% BLOC monograph 
is the following one.
\begin{theorem}
\label{mainth}
Let $n$ and $k$ be two integers such that $k \geq 2n \geq 0$.
Let $D$ be a degree $k$ rational generalised clover.\index{NN}{lDGamma@$\ell(D;\Gamma)$} \index{NN}{ZnD@$Z_{n}(D)$}
$$\begin{array}{llll}
Z_{n}(D)&=& 0 &\mbox{if}\; k > 2n,\\
Z_{n}(D)&=&\sum_{\Gamma}\frac{\ell(D;\Gamma)}{\sharp \mbox{Aut}(\Gamma)}[\Gamma] &\mbox{if}\; k=2n,\end{array}$$
where the sum runs over all degree $n$ Jacobi diagrams $\Gamma$.
\end{theorem}

The above theorem has the following immediate corollary.

\begin{corollary}
Let $n$ and $k$ be two integral numbers such that $k \geq 2n$.
Let $\nu$ be a degree $n$ invariant of homology spheres valued in a real vector space. Let $D$ be a degree $k$ integral generalised clover.
$$\begin{array}{llll}
\nu(D)&=& 0 &\mbox{if}\; k > 2n,\\
\nu(D)&=&\sum_{\Gamma}\frac{\ell(D;\Gamma)}{\sharp \mbox{Aut}(\Gamma)}W_{\nu}(\Gamma) &\mbox{if}\; k=2n,\end{array}$$
where the sum runs over all degree $n$ Jacobi diagrams $\Gamma$. The product $\ell(D;\Gamma)W_{\nu}(\Gamma)$ does not depend on the vertex-orientation.
\end{corollary}
Note that this corollary applies to $(Z_{LMO})_n$.

The following corollary of Theorem~\ref{mainth} is proved in Section~\ref{seccompwal}. 
\begin{theorem}
\label{thcompwal}
For any rational homology sphere $M$, if $\lambda_W$ \index{NN}{lambdaW@$\lambda_W$}
denotes the Walker invariant normalized as in~\cite{wal}, then 
$$Z_{1}(M)=\frac{\lambda_W(M)}{4}[\tata]$$
\end{theorem}

As another corollary, we could describe the generalised clovers in the setting of the Habiro-Goussarov filtration of integral homology spheres,
and give an algebraic version of the clover calculus for integral homology $3$-spheres over $\QQ$. This is done
in \cite{al}, where an algebraic version of the clover calculus for integral homology $3$-spheres over $\ZZ$ is also given.

Theorem~\ref{mainth} was observed by D.~Thurston and G.~Kuperberg 
when the rational homology handlebodies $B^i$ are obtained from the $A^i$ by 
composition of the identification of the boundaries by a Torelli diffeomorphism 
that is a diffeomorphism that induces the identity in homology in \cite{kt}.
This particular case is enough to conclude that $Z$ is a universal finite type invariant of integral homology spheres.
Together with the fact that $Z(S^3)=1$, it implies Theorem~\ref{thcompwal}
for integral homology spheres that is also due to D.~Thurston and G.~Kuperberg.

The proof of Theorem~\ref{mainth} strongly relies on the Kontsevich-Kuperberg-Thurston construction of 
%$Z$. We shall therefore recall their results in the next sections.
$Z$ that is given in \const ~and not repeated here.
\newpage
\section{Sketch of the proof of Theorem~\ref{mainth}}
\label{secskp}
\setcounter{equation}{0}

We refer to the construction of $Z$ given in \rconstsecsta.
However we choose the homogeneous volume form $\omega_{S^2}$ \index{NN}{omegaS2@$\omega_{S^2}$} on $S^2$ with total volume $1$ once for all, and we only consider forms on $C_2(M)$  that are antisymmetric (with respect to the exchange of two points).
In particular, 
in this article, 
%%% BLOC from now on,
we say that a $2$-form $\omega_M$ on $C_2(M)$ is {\em fundamental\/} 
\index{TT}{form!fundamental} with respect to a trivialisation $\tau_M$ of $T(M \setminus \infty)$ that is standard near $\infty$ if it is antisymmetric and fundamental with respect to 
$\tau_M$ and $\omega_{S^2}$ in the sense
of \inconstdefformfund.
Similarly, 
here,
%%% BLOC from now on,
 a two-form $\omega_M$ on $C_2(M)$ or on $\partial C_2(M)$ is {\em admissible\/} if \index{TT}{form!admissible} \\
$\bullet$ its restriction to $\partial C_2(M) \setminus ST(B_M)$ is $p_M(\tau_M)^{\ast}(\omega_{S^2})$ for some trivialisation $\tau_M$ of $T(M \setminus \infty)$ standard near $\infty$, and,\\ 
$\bullet$ it is closed,
and antisymmetric.\\
In particular, all fundamental or admissible forms coincide on $\partial C_2(M) \setminus ST(B_M)$.

\medskip
 
Fix a rational generalised clover 
$D=(M;\sharp N;(A^i)_{i=1,\dots,\sharp N};(B^i)_{i=1,\dots,\sharp N})$.\\
For $I \subseteq N=\{1,2,\dots,\sharp N\}$, set $M_I=M_I(D)$. \index{NN}{MI@$M_I$}

For any $i \in N=\{1,\dots,\sharp N\}$, fix disjoint simple closed curves $(a^i_j)_{j=1,\dots,g_i}$ and simple closed curves  $(z^i_j)_{j=1,\dots,g_i}$ on $\partial A^i$, such that
$$\CL_{A^i}=\oplus_{j=1}^{g_i} [a^i_j],$$ 
%$a^i_j$ and $z^i_k$ are disjoint if $j \neq k$, $a^i_j$ and $z^i_j$ 
%intersect once
so that 
$$\langle a^i_j,z^i_k \rangle_{\partial A^i}=\delta_{jk}=\left\{\begin{array}{ll} 0 &  \mbox{if} \; j \neq k\\
1 & \mbox{if} \; j = k.\end{array}\right.$$

Let $\tau_M$ \index{NN}{tauM@$\tau_M$} be a trivialisation of $M \setminus \infty$ that is standard near $\infty$. 
%Let $\omega_{S^2}$ denote the homogeneous two-form on $S^2$ 
%with total area $1$. 
Define $\omega(\tau_M)=p_M(\tau_M)^{\ast}(\omega_{S^2})$ \index{NN}{omegatauM@$\omega(\tau_M)$}
on $\partial C_2(M)$.

For $j \in N$,
define an antisymmetric closed 2-form $\omega_j$ \index{NN}{omegaj@$\omega_j$} on $ST(B^j)$ that coincides with $\omega(\tau_M)$ on $ST(B^j)_{|\partial B^j}$, and define a trivialisation
$\tau_j^{\CC}$ \index{NN}{taujC@$\tau_j^{\CC}$}
of $TB^j \otimes \CC$ that is the complexification of $\tau_M$ on $\partial B^j$ as follows.

When the restriction of $\tau_M$ to $\partial B^j$ extends to $B^j$ as a trivialisation $\tau_j$, \index{NN}{tauj@$\tau_j$}
simply set $\omega_j=p_{M_j}(\tau_j)^{\ast}(\omega_{S^2})$, and $\tau_j^{\CC}
= \tau_j \otimes \CC$.

When $\tau_{M|\partial B^j}$ does not extend to $B^j$, there exists a curve $c^0_j$ in $\partial B^j$ such that the twist 
$\CT_{c^0_j} \circ \tau_{M|\partial B^j}$ 
across $c^0_j$ of $\tau_{M|\partial B^j}$ (see Definition~\ref{deftwicur}) extends to $B^j$. The curve $c^0_j$ inherits a framing from  $\partial B^j$. Let $N(c_j)=[a,b] \times c_j \times [-1,1]$ denote a neighborhood of a curve $c_j$ parallel to $c^0_j$ inside $B^j$ such that $\{b\} \times c_j \times [-1,1] \subset \partial B^j$ and $c^0_j=\{b\} \times c_j \times \{0\}$. 

Then $\tau_{M|\partial B^j}$ extends to the closure of 
$(B^j \setminus N(c_j))$ as a trivialisation $\tau_j$, \index{NN}{tauj@$\tau_j$} and $\tau_M$ extends to $N(c_j)$ so that 
$$\tau_{M|\partial N(c_j)}=\left\{\begin{array}{ll} \tau_j & \mbox{over}\; \partial N(c_j)\setminus (\{a\} \times c_j \times [-1,1]) \\
\CT_{c_j}^{-1} \circ \tau_j  & \mbox{over}\; \{a\} \times c_j \times [-1,1].
\end{array} \right.$$
Then with the notation of Subsection~\ref{subspecad}, set
$\omega_j=\omega(c_j;\tau_j,\tau_{M| N(c_j)})$ \index{NN}{omegaj@$\omega_j$} and \index{NN}{taujC@$\tau_j^{\CC}$}
$\tau_j^{\CC}=\tau_{\CC}(c_j;\tau_j,\tau_{M| N(c_j)})$.

\begin{remark}
\label{rkfirstap}
When $A$ is a compact oriented $3$-manifold bounded by a compact connected surface, set \index{NN}{LAZ2@$\CL_A^{\ZZ/2\ZZ}$}
$$\CL_A^{\ZZ/2\ZZ}=\mbox{Ker}(H_1(\partial A; \ZZ/2\ZZ) \longrightarrow H_1( A; \ZZ/2\ZZ))$$
and set
$$\CL_A^{\ZZ}=\mbox{Ker}(H_1(\partial A; \ZZ) \longrightarrow H_1( A; \ZZ))$$
so that $\CL_A=\CL_A^{\ZZ} \otimes \QQ$.
When $A$ is a $\ZZ$-handlebody, 
$\CL_A^{\ZZ/2\ZZ}=\CL_A^{\ZZ} \otimes {\ZZ/2\ZZ}$.
If $\CL_{A^j}^{\ZZ/2\ZZ}=\CL_{B^j}^{\ZZ/2\ZZ}$, then the restriction of $\tau_M$ to $\partial B^j$ extends to $B^j$.
This would always be the case, if only $\ZZ$-handlebodies were involved. This would also be the case if $B^j$ were obtained from $A^j$ by twisting the identification of $\partial A^j$ by a {\em \indexT{Torelli homeomorphism}\/} of $\partial A^j$, that is a homeomorphism that induces the identity in homology.
However, in the general case, $\tau_{M|\partial B^j}$ may fail to extend to $B^j$. See Example~\ref{exabadlag}.
 
Nevertheless, in a first approach to the proof, the reader can assume that $\tau_{M|\partial B^j}$ extends to $B^j$ and forget about Subsections~\ref{subspecad} and~\ref{subproofpunspec}
and the previous paragraph that are useless in the case when $\tau_{M|\partial B^j}$ extends to $B^j$.
\end{remark}
 
For $I \subseteq N$, equip $\partial C_2(M_I)$ with the admissible 
 $2$-form 
$\omega(M_I)$ \index{NN}{omegaMI@$\omega(M_I)$} that coincides with $\omega(\tau_M)$  on $ST(B_{M_I}) \setminus \cup_{j \in I} ST(B^j)$ and that 
is equal to $\omega_j$ on $ST(B^j)$. 
Similarly, equip $T(M_I \setminus \infty) \otimes \CC$ with the trivialisation $\tau_I^{\CC}$
that is the complexification of $\tau_M$ over $M \setminus (\infty \cup (\cup_{j \in I} B^j))$, and that equals $\tau_j^{\CC}$ over $B^j$.
Define $Z(\omega(M_I))=(Z_n(\omega(M_I))_{n \in \NN}$ \index{NN}{ZomegaMI@$Z(\omega(M_I))$} as in \inconstpropzad.
Then according to Proposition~\ref{propintpunspec}, 
%(or to Theorem~\ref{thkktf} 
(or to \rconstthkktfra ~when $\tau_{M|\partial B^j}$ extends to $B^j$), \index{NN}{ZMI@$Z(M_I)$}
$$Z(M_I)= Z(\omega(M_I)) \exp(\frac{p_1(\tau_I^{\CC})}{4}\xi).$$
%%%SIGNE

Furthermore, the $p_1(\tau_I^{\CC})$ are related by the following lemma.

\begin{lemma}
\label{lemtrivrel} 
Let $p(i)=p_1(\tau(M_{\{i\}}))-p_1(\tau(M))$. \index{NN}{p(i)@$p(i)$} Then,
for any subset $I$ of $N$,\index{NN}{ponetauIC@$p_1(\tau_I^{\CC})$}
$$p_1(\tau_I^{\CC})=p_1(\tau(M))+\sum_{i \in I}p(i).$$
\end{lemma}
\bp
Fix $j \in N$. 
Let $Y$ be a rational homology handlebody with boundary $-\partial A^j$
and with the same lagrangian $\CL_Y$ as $(M_I \setminus A^j)$ for any
$I \subseteq (N\setminus \{j\})$. Embed the standard neighborhood  $(S^3 \setminus B^3)$  of $\infty$ into $Y$. Let $\tau_Y$ be a trivialisation of 
$T(Y \setminus \infty) \otimes \CC$ that coincides with $\tau_j^{\CC}$ on $\partial A^j$ and that is standard near $\infty$.
Let $\tau(\tau_Y,\tau_j^{\CC})$ denote the trivialisation of 
$T( B^j \cup Y \setminus \infty) \otimes \CC$ that coincides with $\tau_Y$ over $Y\setminus \infty$ and with
$\tau_j^{\CC}$ over $B^j$. Similarly define the  trivialisation $\tau(\tau_Y,\tau_{M|A^j})$ of 
$T(A^j \cup Y \setminus \infty) \otimes \CC$.
Then it is enough to prove that $$p_1(\tau(\tau_Y,\tau_j^{\CC}))-p_1(\tau(\tau_Y,\tau_{M|A^j}))$$ is independent of the rational homology handlebody $Y$ with boundary $-\partial A^j$
and with lagrangian $\CL_Y$, and that it is independent of $\tau_Y$.

As any $3$-manifold, the manifold $\left(-A^j \cup ([0,1] \times  \partial A^j) \cup B^j\right)$ bounds a $4$-manifold $W_j$. Then 
$$W_j \cup_{[0,1] \times  \partial A^j} [0,1] \times \overline{Y \setminus (S^3 \setminus B^3)}$$ is a cobordism between $B_{Y \cup A^j}$ and $B_{Y \cup B^j}$
whose signature is independent of the rational homology handlebody $Y$ with boundary $-\partial A^j$ and with lagrangian $\CL_Y$ and can be adjusted to zero for any such after
a connected sum with copies of $\CC P^2$ or $\overline{\CC P^2}$ in the interior of $W_j$.
After this adjustment, $p_1(\tau(\tau_Y,\tau_j^{\CC}))-p_1(\tau(\tau_Y,\tau_{M|A^j}))$
is just the obstruction to extend the trivialisation of $TW_j$ on 
$\partial W_j$ induced by the given trivialisations $\tau_j^{\CC}$
and $\tau_{M|A^j}$. Therefore it is independent of $Y$ and $\tau_Y$.

%Let us first show that $B^j$ is obtained from $A^j$ by surgery 
%on a framed link in the interior of $A^j$, so that the identifications of the 
%boundaries are induced by the surgery. Of course, it is enough to prove this %when $B=B^j$ is the standard handlebody with meridians the $a^i_j$.
%Consider the $g^j$-component framed link ${\bf L}$ of $A^j$ made of curves
%$K^i$ parallel to the $a^i_j$ (framed by the surfaces). 
%Then surgery on ${\bf L}$ produces the connected sum of $B$ and a closed %manifold $M$ where the connected sum is performed inside $B$.
%Now, we conclude by performing in $(M \setminus B^3)$ a surgery that transforms
%$M$ into $S^3$.
\eop

Let $[-4,4] \times (\coprod_{i \in N} \partial A^i)$ be a tubular neighborhood of $(\coprod_{i \in N} \partial A^i)$ 
in $M$. This neighborhood intersects $A^i$ as $[-4,0] \times \partial A^i$.
Let $[-4,0] \times \partial A^i$ be a neighborhood of $\partial B^i=\partial A^i$ in $B^i$.
The manifold $M_{\{i\}}=M_i$ is obtained from $M$ by removing $\left(A^i \setminus
(]-4,0] \times \partial A^i)\right)$ and by gluing back $B^i$ along $(]-4,0] \times \partial A^i)$.

Let $C_I^i \subset M_I$, \index{NN}{CIi@$C_I^i$} $C_I^i=A^i$ if $i \notin I$,  $C_I^i=B^i$ if $i \in I$.
Let $\eta_{[-1,1]}$ \index{NN}{etamin@$\eta_{[-1,1]}$} be a one-form with compact support in $]-1,1[$ such that  $\int_{[-1,1]}\eta_{[-1,1]}=1$. 
Let $(a^i_j  \times [-1,1])$ be a tubular neighborhood of  $a^i_j$ in $\partial A^i$. 
Let $\eta(a^i_j)$ \index{NN}{etaaij@$\eta(a^i_j)$} be a closed one-form on $C_I^i$ such that
the support of $\eta(a^i_j)$ intersects $[-4,0] \times \partial A^i$ inside  $[-4,0] \times (a^i_j  \times [-1,1])$, where $\eta(a^i_j)$ can be written
$$\eta(a^i_j)=\pi_{[-1,1]}^{\ast}(\eta_{[-1,1]}).$$ 
Note that the forms $\eta(a^i_j)$ on $A^i$ and $B^i$ induce closed one-forms still denoted by $\eta(a^i_j)$ \index{NN}{etaaij@$\eta(a^i_j)$} on $(A^i \cup_{\partial A^i} -B^i)$ that restrict to the previous ones.

%Let $S(a^i_j)$ be a 2-chain in $C_I^i$ with boundary $a^i_j$ that intersects 
%$[-4,4], \times \partial A^i$ along $[-4,0] \times a^i_j$.

%surface in $(A^i \cup_{\partial A^i} -B^i)$.
%Let $S^A(a^i_j)=S(a^i_j)\cap A^i$.
% and $S^A(a^i_j)=S(a^i_j)\cap (-B^i)$.
%Assume that $\partial S^A(a^i_j)=a^i_j$.
%Let $S(a^i_j) \times [-1,1]  \subset(A^i \cup_{\partial A^i} -B^i)$ be a tubular %neighborhood of  $S(a^i_j)$ that intersects $\partial A^i$ along $(a^i_j \subset %S(a^i_j)) \times [-1,1]$ (so that the two ways of writing the latter product %coincide).
%%Let $S^I(a^i_j)$ be a surface embedded in $C_I^i$ whose boundary is $\partial %S^I(a^i_j) = a^i_j$, and let  $S^I(a^i_j) \times [-1,1]  \subset C_I^i$ be a %tubular neighborhood of  $S^I(a^i_j)$.
%Let $\eta(a^i_j)$ be a one-form on $(A^i \cup_{\partial A^i} -B^i)$ supported in %this tubular neighborhood 
%that may be written 
%as $$\eta(a^i_j)=\pi_{[-1,1]}^{\ast}(\eta).$$

The following innocent-looking proposition~\ref{propnormasim} is the key proposition.
As of yet, I do not have a simple proof for it. Fix a degree $\sharp N$
generalised clover $D$, and, for $I \subseteq N=\{1,2,\dots,\sharp N\}$, set $M_I=M_I(D)$.
Let $$p_{12}: C_{2}(M_I) \longrightarrow (M_I)^2$$ be the natural projection. For $X \subset M_I$, $C_2(X)$ \index{NN}{CtwoX@$C_2(X)$} denotes $p_{12}^{-1}(X^2) \subset C_2(M_I)$.

\begin{proposition}
\label{propnormasim}
There exist admissible $2$-forms $\omega(M_I)$ \index{NN}{omegaMI@$\omega(M_I)$} on $C_2(M_I)$, that extend the $2$-forms $\omega(M_I)$ defined on $\partial C_2(M_I)$ such that:\\
$\bullet$ For any $I,J \subseteq N$, $\omega(M_I)=\omega(M_J)$ wherever it makes sense in $C_{2}(M)$,\\
that is
on $C_2\left( \left(M \setminus \bigcup_{i \in I \cup J} \mbox{Int}(A^i)\right) \cup  \bigcup_{j \in I \cap J} B^j\right)$.\\
$\bullet$ On $C_I^i \times C_I^k$,\\
$$\omega(M_I)=\sum_{\begin{array}{cc}j=1,\dots,g_i\\ \ell=1,\dots,g_k \end{array}}\ell(z^i_j,z^k_{\ell})\pi_{C_I^i}^{\ast}(\eta(a^i_j)) \wedge \pi_{C_I^k}^{\ast}(\eta(a^k_{\ell}))$$
where $\ell(.,.)$ \index{NN}{l@$\ell$} stands for the linking number.
\end{proposition}

Using these forms, we can easily prove the following lemma.

\begin{lemma}
\label{lemthindiv}
Let $n \in \NN$. Let $\Gamma$ be a degree $n$ oriented Jacobi diagram:
\index{NN}{lDGamma@$\ell(D;\Gamma)$}
$$\begin{array}{llll}
\sum_{I \subseteq N}(-1)^{\sharp I}I_{\Gamma}(\omega(M_I))&=& 0 &\mbox{if}\; \sharp N > 2n,\\
&=&\ell(D;\Gamma) &\mbox{if}\; \sharp N=2n.\end{array}$$
\end{lemma}
\bp
Let $E$ be the set of edges of $\Gamma$.
We want to compute $$\Delta=\sum_{I \subseteq N}(-1)^{\sharp I}I_{\Gamma}(\omega(M_I)).$$
Number the vertices of $\Gamma$ so that $V(\Gamma)=\{1,2,\dots,2n\}$
and $$\breve{C}_{V(\Gamma)}(M_I)=(M_I \setminus \infty)^{2n} \setminus \{\mbox{all diagonals}\}=\mbox{Int}(C_{2n}(M_I)).$$
Fix $i \in N=\{1,2,\dots,\sharp N\}$, the contributions to $\Delta$ of the integrals of $\bigwedge_{e \in E}p_e^{\ast}(\omega(M_I))$ over 
$$\mbox{Int}(C_{2n}(M_I)) \cap (M_I \setminus C^i_I)^{2n}$$
are identical for $I=K$ and $I=K\cup\{i\}$ for any $K \subseteq N \setminus \{i\}$.
Since they enter the sum with opposite signs, they cancel each other.
This argument allows us to get rid of all the contributions of the integrals over
$$\bigcup_{i \in N}\left(\mbox{Int}(C_{2n}(M)) \cap (M_I \setminus C^i_I)^{2n}\right).$$
Thus, we are left with the contributions of the integrals over the subsets $P_I$ of $\mbox{Int}(C_{2n}(M_I)) \subset (M_I)^{2n}$ such that:\\
For any $i \in N$, any element of $P_I$ projects onto $C^i_I$ under at least one of the $(2n)$ projections onto $M_I$.
These subsets $P_I$ are clearly empty if $\sharp N >2n$, and the lemma is proved in this case. Otherwise, $P_I$ is equal to
$$\cup_{\sigma \in \Sigma_N}\prod_{i=1}^{2n}C_I^{\sigma(i)}$$
where $\Sigma_N$ is the set of permutations of $N$.
We get
$$\Delta=\sum_{\sigma \in \Sigma_N}\Delta_{\sigma}$$
with 
$$\Delta_{\sigma}=\sum_{I \subseteq N}(-1)^{\sharp I} \int_{\prod_{i=1}^{2n}C_I^{\sigma(i)}}\bigwedge_{e \in E}p_e^{\ast}(\omega(M_I)).$$
It is enough to prove that $\Delta_{\sigma}=\ell(D;\Gamma;\sigma)$. \index{NN}{lDGammasigma@$\ell(D;\Gamma;\sigma)$}
Recall that the vertices are numbered.
% by $\sigma:V(\Gamma) \longrightarrow \{1,2,\dot,2n\}$.
For any $i \in N$, $$p_i:\mbox{Int}(C_{2n}(M_I)) \longrightarrow M_I$$
denotes the projection onto the $i^{\mbox{\small th}}$ factor.
When $e$ is an oriented edge from the vertex $x(e) \in V(\Gamma)$ to $y(e) \in V(\Gamma)$.
$$p_e^{\ast}(\omega(M_I))_{|\prod_{i=1}^{2n}C_I^{\sigma(i)}}=$$
$$\sum_{\begin{array}{cc}j=1,\dots,g_{\sigma(x(e))}\\ \ell=1,\dots,g_{\sigma(y(e))} \end{array}}\ell(z^{\sigma(x(e))}_j,z^{\sigma(y(e))}_{\ell})p_{x(e)}^{\ast}(\eta(a^{\sigma(x(e))}_j)) \wedge p_{y(e)}^{\ast}(\eta(a^{\sigma(y(e))}_{\ell}))
$$
Recall that when $c$ is a half-edge, $v(c)$ denotes the label of the vertex contained in that half-edge, and that $H(\Gamma)$ denotes the set of half-edges. We shall also use the notation $x(e)$ and $y(e)$ for the corresponding halves of an edge $e$. Let $F$ denote the set of maps $f$
from $H(\Gamma)$ to $\NN$ such that for any $c \in H(\Gamma)$, ${f(c)} \in \{1,2,\dots,g_{\sigma(v(c))}\}$.
$$\Delta_{\sigma}=
\sum_{f \in F}\left(\prod_{e \in E}\ell(z^{\sigma(x(e))}_{f(x(e))},z^{\sigma(y(e))}_{f(y(e))})\right) I(f)$$
with
$$I(f)= \int_{\prod_{i=1}^{2n}(A^{\sigma(i)} \cup -B^{\sigma(i)})} \bigwedge_{c \in H} p_{v(c)}^{\ast}(\eta(a^{\sigma(v(c))}_{{f(c)}}))$$
Recall that the set of half-edges has been ordered in two equivalent ways (up to an even permutation)
in the beginning of 
\inconstsubapseckkt. 
%%% BLOC
The above exterior product must be computed with this order. In particular, thinking of this order as being given by the order of the vertices, we get:
$$I(f)=\prod_{i=1}^{2n} \int_{A^{\sigma(i)} \cup (-B^{\sigma(i)})} \bigwedge_{c \in H \cap v^{-1}(i)}\eta(a^{\sigma(i)}_{f(c)}).$$
%$$=\prod_{i=1}^{2n} \int_{A^{i} \cup (-B^{i})} \bigwedge_{c \in H \cap %v^{-1}(\sigma^{-1}(i)}\eta(a^{i}_{f(c)})$$
with the order of the three half-edges given by the vertex-orientation.
Now,
$$ \int_{A^{\sigma(i)} \cup (-B^{\sigma(i)})} \bigwedge_{c \in H \cap v^{-1}(i)}\eta(a^{\sigma(i)}_{f(c)})=\pm \CI({A^{\sigma(i)}},B^{\sigma(i)})(\bigotimes_{c \in H \cap v^{-1}(i)} a^{\sigma(i)}_{f(c)}),$$
where the sign $\pm$ does not depend on $i$.
It is easy to conclude that 
$\Delta_{\sigma}=\ell(D;\Gamma;\sigma)$.
\eop

\noindent{\sc Proof of Theorem~\ref{mainth}:}
This lemma easily implies
$$\begin{array}{llll}
\sum_{I \subseteq N}(-1)^{\sharp I}Z_n(M_I,\omega(M_I))&=& 0 &\mbox{if}\; \sharp N > 2n,\\
&=&\sum_{\Gamma}\frac{\ell(D;\Gamma)}{\sharp \mbox{\small Aut}(\Gamma)}[\Gamma] &\mbox{if}\; \sharp N=2n.\end{array}$$
and it suffices to deal with the framing corrections by proving that if $\sharp N \geq 2n$, then
$$\sum_{I \subseteq N}(-1)^{\sharp I}Z_n(M_I)=\sum_{I \subseteq N}(-1)^{\sharp I}(Z_n(M_I,\omega(M_I))$$
By definition, we have that
$$\sum_{I \subseteq N}(-1)^{\sharp I}Z_n(M_I)=$$
$$\sum_{I \subseteq N}(-1)^{\sharp I}\left(Z_n(M_I,\omega(M_I)) + \sum_{j<n}Z_j(M_I,\omega(M_I)) P_{n-j}(I)\right)$$
where $P_{n-j}(I)$ stands for an element of $\CA_{n-j}(\emptyset)$ that is a combination of $m[\Gamma]$ where the $m$ are monomials in $p_1(\tau(M))$ and in the $p(i)$ of degree
at most $(n-j)$, for degree $(n-j)$ Jacobi diagrams $\Gamma$. Furthermore, such an $m[\Gamma]$ appears in $P_{n-j}(I)$ if and only if $m$ is a monomial in the variables $p_1(\tau(M))$ and $p(i)$ for $i \in I$.
Therefore, we can rewrite the sum of the annoying terms by factoring out the $m[\Gamma]$. Let $K \subset N$ be the subset of the $i$ such that $p(i)$ appears in $m$.
$(\sharp K \leq n-j)$. The factor of $m[\Gamma]$ reads
$$\sum_{K \subseteq I \subseteq N}(-1)^{\sharp I}Z_j(M_I,\omega(M_I)).$$
Therefore, the sum runs over the subsets of $N \setminus K$ whose cardinality is at least $\sharp N +j - n$. Since $\sharp N \geq 2n$ and $j<n$, $\sharp N-n \geq n >j$, hence $\sharp N +j - n > 2j$ and the preceeding lemma ensures that the above sum is zero.
This concludes the reduction of the proof of Theorem~\ref{mainth} to the construction of {\em special admissible forms\/} in Subsection~\ref{subspecad},
and to the proofs of Proposition~\ref{propintpunspec} in Subsection~\ref{subproofpunspec} and Proposition~\ref{propnormasim} in Section~\ref{secnormap}.
\eop
\newpage
\section{Preliminaries to the simultaneous normalization of the forms}
\setcounter{equation}{0}
\label{secnormaprel}

In this section, Subsection~\ref{subprelprel} is useful even in the case
 when $\tau_{M|\partial B^j}$ extends to $B^j$ for any $j$, while
Subsections~\ref{subspecad} and~\ref{subproofpunspec} are useless in that case.

\subsection{Preliminaries}
\label{subprelprel}

Let $(e_1,e_2,e_3)$ denote the standard basis 
$((1,0,0),(0,1,0),(0,0,1))$ of $\RR^3$, 
and let $v_i:\RR^3 \longrightarrow \RR$ \index{NN}{vi@$v_i$} denote the 
$i^{\mbox{th}}$ coordinate with respect to this basis.
Let $R_{\theta}=R_{\theta,e_1}$ \index{NN}{Rtheta@$R_{\theta}$} denote the rotation of $\RR^3$ 
with axis directed by $e_1$ and with angle $\theta$.
Let $\omega_{S^2}$ \index{NN}{omegaS2@$\omega_{S^2}$} denote the homogeneous two-form on $S^2$
with total area $1$. 
When $\alpha \in S^2$, 
and when $v$ and $w$ are two tangent vectors of $S^2$ at $\alpha$, \index{NN}{omegaS2@$\omega_{S^2}$}
$$\omega_{S^2}(v \wedge w)=\frac{1}{4\pi}det(\alpha,v,w)$$
where $\alpha \wedge v \wedge w = 
det(\alpha,v,w) e_1 \wedge e_2 \wedge e_3$ in $\bigwedge^3\RR^3$.

\begin{lemma}
\label{lemprimtwist}
Let \index{NN}{Tcal@$\CT_k$} $$\begin{array}{llll}\CT_k:&\RR \times S^2 &\longrightarrow &S^2\\
& (\theta,\alpha) & \mapsto & R_{k\theta}(\alpha). \end{array}$$
Then with the notation above, \index{NN}{Tcal@$\CT_k$}
$$\CT_k^{\ast}(\omega_{S^2})=\CT_0^{\ast}(\omega_{S^2}) 
+ \frac{k}{4\pi} d\theta \wedge dv_1$$
\end{lemma}

\bp
Since $R_{\theta}$ preserves the area, the restrictions of 
$\CT_k^{\ast}(\omega_{S^2})$ and $\CT_0^{\ast}(\omega_{S^2})$ coincide on 
$\bigwedge^2T_{(\theta,\alpha)}(\{\theta \}\times S^2)$.
Therefore, we are left with the computation of 
$\left(\CT_k^{\ast}(\omega_{S^2})-\CT_0^{\ast}(\omega_{S^2})\right)(u \wedge v)$ 
when $u \in T_{(\theta,\alpha)}(\RR \times \{\alpha \})$
and $v \in T_{(\theta,\alpha)}(\{\theta \}\times S^2)$, where of course,
$\CT_0^{\ast}(\omega_{S^2})(u \wedge v)=0$, and by definition,
$$\CT_k^{\ast}(\omega_{S^2})(u \wedge v)=\frac{1}{4\pi}det(R_{k\theta}(\alpha),
D\CT_k(u),D\CT_k(v)).$$
Since $D\CT_k(v)=R_{k\theta}(v)$, 
and since $R_{k\theta}$ preserves the volume in $\RR^3$,
$$\CT_k^{\ast}(\omega_{S^2})(u \wedge v)=\frac{1}{4\pi}det(\alpha,
R_{-k\theta}(D\CT_k(u)),v).$$
Now, let $\alpha_i$ stand for $v_i(\alpha)$.
$D\CT_k(u)=kd\theta(u)R_{k\theta+\pi/2}(\alpha_2 e_2 +\alpha_3 e_3)$.
Therefore,
$\CT_k^{\ast}(\omega_{S^2})(u \wedge v)=\frac{kd\theta(u)}{4\pi}det(\alpha,
-\alpha_3 e_2 +\alpha_2 e_3,v),$ and, 
$$\begin{array}{ll}\CT_k^{\ast}(\omega_{S^2})(u \wedge .)
&=\frac{kd\theta(u)}{4\pi}
det\left(\begin{array}{ccc}\alpha_1 & 0 & dv_1\\
\alpha_2 & -\alpha_3 & dv_2\\
\alpha_3 & \alpha_2 & dv_3\end{array}\right)\\
&=\frac{kd\theta(u)}{4\pi}
\left((1-\alpha_1^2)dv_1 -\alpha_1\alpha_2 dv_2-
\alpha_1\alpha_3 dv_3\right)\\
&=\frac{kd\theta(u)}{4\pi} dv_1.\end{array}$$
\eop

\begin{definition}
\label{deftwicur}
Let $S$ be an oriented surface, let $T$ be an $\RR^3$-bundle over $S$, and let
$\tau: T\longrightarrow S \times \RR^3$ be a trivialisation of $T$.
Let $\varepsilon >0$ be a small positive real number.
Let \index{NN}{theta@$\theta$} $\theta : [-1,1] \longrightarrow [0,2 \pi]$ 
be a smooth map that maps $[-1,-1+\varepsilon]$ to $0$, 
that increases from $0$ to $2\pi$ on $[-1+\varepsilon, 1-\varepsilon]$, and such that $\theta(-x)+\theta(x)=2\pi$ for any $x \in [-1,1]$.
When $c \times [-1,1]$ is an oriented neighborhood of an oriented curve $c$
in $S$, 
let \index{NN}{thetac@$\theta(c)$} $$\begin{array}{llll}\theta(c):&S& \longrightarrow &\frac{[0,2 \pi]}{0 \sim 2\pi}\\
&x \notin c \times [-1,1]& \mapsto & 0\\ &(\gamma,u) \in c \times [-1,1]& \mapsto & \theta(u).\end{array}$$
A {\em twist of $\tau$\/} \index{TT}{twist of a trivialisation} across an oriented curve $c$ with oriented neighborhood $c \times [-1,1]$ is the trivialisation
\index{NN}{Tcalc@$\CT_c$}
$$\begin{array}{llll}\CT_c \circ \tau:&T &\longrightarrow &S \times \RR^3\\
&\tau^{-1}(x,v) & \mapsto & (x,R_{\theta(c)(x)}(v))
\end{array}$$
where $R_{\theta}$ is defined in the beginning of Subsection~\ref{subprelprel}.
\end{definition}

\begin{definition}
\label{defkiobs}
Let $S$ be a non-necessarily oriented compact surface with possible boundary, and let $v$ be a nonzero vector field of $TS$
on the boundary $\partial S$ of $S$.
Then the {\em \indexT{Euler number}\/} $\chi(TS; v_{|\partial S})$ \index{NN}{chiTSv@$\chi(TS; .)$} of $v$  is the obstruction to extend
the vector field $v$ of $TS_{|\partial S}$ to a nonzero section of $TS$.
More precisely, if $S$ is a disk then its unit tangent bundle $S(TS)$ is isomorphic to $S^1 \times S$, and
$\chi(TS; v_{|\partial S})$ is the degree of the composite map
$$ \partial S \;\hfl{v}\; S(TS)\cong S^1 \times S \;\hfl{\pi_{S^1}}\; S^1.$$
Note that the orientation of $S$ does not matter for this definition (because
changing it changes both orientations of $\partial S$ and of the fiber $S^1$).
When $S$ is connected, then the vector field  $v_{|\partial S}$ can be extended as a nonzero section $v$ outside the interior of a disk $D$ in the interior of $S$, and $\chi(TS; v_{|\partial S})=\chi(TD; v_{|\partial D})$. 
When $S$ is not connected, then 
$\chi(TS; v_{|\partial S})$ \index{NN}{chiTSv@$\chi(TS; .)$} is 
the sum of the 
numbers associated to the different components of $S$. 
\end{definition}
Note that when the boundary of $S$ is empty, $\chi(TS)$ is the Euler characteristic \index{NN}{chiS@$\chi(S)$} $\chi(S)$ of $S$. The following lemma is left to the reader.

\begin{lemma}
\label{lemcalchi}
Let $S$ be an oriented surface
and let $v$ be a nonzero vector field of $TS$
on $\partial S$. Use the two orthogonal unit vector fields on $\partial S$, "outward normal to $S$" $\vec{N}(\partial S)$ and "tangent vector to $\partial S$" $\vec{T}(\partial S)$ to trivialise $TS_{|\partial S}$ by mapping $(\vec{N}(\partial S),\vec{T}(\partial S))$
to the basis $(e_1,e_2)$ of $\RR^2$. Using this trivialisation,
$$S(TS)_{|\partial S}= S(\RR e_1 \oplus \RR e_2) \times \partial S=S^1 \times \partial S.$$
Let $d$ be the degree of the composite map
$$ \partial S \;\hfl{v}\; S(TS)_{|\partial S}=S^1 \times \partial S \;\hfl{\pi_{S^1}}\; S^1.$$
Then \index{NN}{chiTSv@$\chi(TS; .)$} $$\chi(TS; v_{|\partial S})=d+\chi(S).$$
\end{lemma}

\begin{lemma}
\label{lemintsur}
Let $S \times [-1,1]$ denote a collar of a connected oriented surface $S$ with possible boundary. Let $T$ denote the restriction of the tangent bundle
of $S \times [-1,1]$ to $S=S\times \{0\}$. Let $n$ denote a nonzero vector field of $T$ that is tangent to $\{x\} \times [-1,1]$ at $(x,0)$, and let $S(n)$ \index{NN}{Sn@$S(n)$} denote the corresponding section of the spherical bundle $S(T)$ over $S$.
Let $\tau : T \longrightarrow S \times \RR^3$ be a trivialisation that
maps $n(x)$ to $(x,e_1)$ for any $x \in \partial S$. Use $\tau$ to identify
$S(T)$ to $S \times S^2$, and let $\pi_{S^2}(\tau)$ denote the associated projection from $S(T)$ to $S^2$.
Then \index{NN}{chiTSv@$\chi(TS; .)$} $$\int_{S(n)}\pi_{S^2}(\tau)^{\ast}(\omega_{S^2}) = \frac12\chi(TS; \tau^{-1}(.,e_2)_{|\partial S}).$$
\end{lemma}
\bp
We first prove the result for a special trivialisation $\tau(S)$ of $T$.

If the boundary of $S$ is non-empty then $TS$ is trivialisable, and we fix $\tau(S)$ such that $\tau(S)(n(x))=(x,e_1)$ for $x \in S$. Then both sides of 
the equality to be shown vanish.

If the boundary of $S$ is empty, we trivially embed $S \times [-1,1]$ in 
$\RR^3$ and we pull-back the standard trivialisation of $\RR^3$ through this embedding. Then we compute the left-hand side as the degree
of the map "direction of the positive normal" from $S$ to $S^2$. This degree
is $(1-g(S)=\chi(S)/2)$. Therefore it coincides with the right-hand side that is  $(\chi(TS)/2=\chi(S)/2)$.

Of course, both sides of the equality are unchanged under a homotopy of $\tau$
such that $n=\tau^{-1}(e_1)$ on $\partial S$. Now, up to this kind of homotopy, any such trivialisation is obtained from $\tau(S)$ by a twist across a curve $c$ with possible boundary with a neighborhood $c \times [-1,1]$ properly embedded in $S$. Thus, it is enough to prove that both sides vary in the same way under such a twist.
It is clear that the right-hand side varies like half the degree of the map
$\tau^{-1}(.,e_2)$ from $\partial S$ to the $S^1$ fiber of
$S(TS)_{|\partial S}$ equipped with an arbitrary fixed trivialisation.
Therefore, the variation of the right-hand side is 
$$\frac12 \mbox{deg}\left(R_{-\theta(c)(.)}(e_2): \partial S \longrightarrow S^1=S(\RR e_2 \oplus \RR e_3)\right)=-\frac{\langle c,\partial S \rangle }{2}.$$
For the left-hand side, consider
$$S(n) \hookrightarrow S(TS) \;\hfl{\tau(S)}\; S \times S^2 \;\hfl{\theta(c) \times \mbox{Id}} \;\frac{[0,2\pi]}{0\sim 2\pi}  \times S^2 \;\hfl{\CT_1}\; S^2.$$
$\int_{S(n)}\left(\pi_{S^2}(\CT_c\circ\tau(S))^{\ast}\right)(\omega_{S^2})-
\int_{S(n)}\left(\pi_{S^2}(\tau(S))^{\ast}\right)(\omega_{S^2})=$
$$\int_{S(n)}(\theta(c) \times \mbox{Id})^{\ast}(\CT_1^{\ast}(\omega_{S^2})-\CT_0^{\ast}(\omega_{S^2}))$$
where according to Lemma~\ref{lemprimtwist},$$(\theta(c) \times \mbox{Id})^{\ast}(\CT_1^{\ast}(\omega_{S^2})-\CT_0^{\ast}(\omega_{S^2}))=-\frac{1}{4\pi}d(v_1 d \theta(c)).$$
Therefore, according to the Stokes theorem, the variation of the right-hand side is $$\int_{\partial S(n)} -\frac{1}{4\pi}v_1 d \theta(c)=-\frac{\langle c,\partial S \rangle }{2},$$
and we are done.
\eop

\subsection{Extensions of trivialisations on $3$-manifolds}
\label{subuseless}

This section is useless for the proofs. 
It only justifies why I could not avoid the following subsections, and some of the difficulties they contain.

Let $A$ be a compact oriented connected $3$-manifold 
with boundary $\partial A$.
Consider the {\em $\ZZ/2\ZZ$-Lagrangian\/} of $A$ \index{NN}{LAZ2@$\CL_A^{\ZZ/2\ZZ}$}
$$\CL_A^{\ZZ/2\ZZ}=\mbox{Ker}(H_1(\partial A;\ZZ/2\ZZ) \longrightarrow
H_1( A;\ZZ/2\ZZ))$$
This is a Lagrangian subspace of $(H_1(\partial A;\ZZ/2\ZZ);\langle ., . \rangle )$.

Let $K$ be a framed knot in an oriented $3$-manifold $M$, that is a knot equipped with a normal nonzero vector field $\vec{N}$, or equivalently with a parallel (up to homotopies). 
These data induce the direct trivialisation $\tau_K$ of $TM_{|K}$ (up to homotopy) such that $\tau_K^{-1}(e_1)=\vec{TK}$  and $\tau_K^{-1}(e_2)=\vec{N}$, where $\vec{TK}$ is a tangent vector of $K$ that is equipped with an arbitrary orientation.
The homotopy class of the trivialisation $\tau_K$ is well-defined and does not 
depend on the orientation of $K$.

Assume that $K$ bounds a possibly non-oriented compact surface $\Sigma$ that induces the given parallelisation of $K$.
Then if $\tau$ is a trivialisation of the tangent space of $M$ over $\Sigma$, the restriction of $\tau$ to $K$ is not homotopic to $\tau_K$. 
(This is clear when $\Sigma$ is a disk, and the trivialisation must be homotopic in the other cases, since the tangent bundle of an oriented $3$-manifold over a closed surface is trivialisable. Recall that $\pi_1(SO(3))=\ZZ/2\ZZ$.)

If $K$ is a framed knot in an oriented $3$-manifold $M$ and if $\tau$ is a trivialisation of the restriction of $TM$ to $K$, we shall say that
$K$ is {\em $\tau$-bounding\/} if $\tau$ is not homotopic to $\tau_K$.

%Let $A$ be a compact oriented connected $3$-manifold 
%with boundary $\partial A$. 
\begin{proposition}
\label{propuseless}
Let $\partial A$ be a connected oriented compact surface.
Let $\tau$ be a trivialisation of \/$T(\partial A \times [-2,2])$.
Then there exists a unique map 
%a $\ZZ/2\ZZ$-valued map from 
$$\phi_{\tau} : H_1(\partial A;\ZZ/2\ZZ) \longrightarrow \frac{\ZZ}{2\ZZ}$$
such that 
\begin{enumerate}
\item when $x$ is a connected curve of $\partial A=\partial A \times \{0\}$,
$\phi_{\tau}(x)=0$ if and only if $x$ (equipped with its parallelisation induced by $\partial A$) is $\tau$-bounding
and,
\item
$$\phi_{\tau}(x+y)=\phi_{\tau}(x) + \phi_{\tau}(y) + \langle x, y \rangle .$$
\end{enumerate}
Let $c$ be curve of $\partial A$ and let $\CT_c$ denote the twist across $c$, then
$$\phi_{\CT_c \circ \tau}(x)=\phi_{\tau}(x) + \langle x, c \rangle.$$
When $A$ is a compact oriented connected $3$-manifold 
with boundary $\partial A$, $\tau$ extends as a trivialisation over $A$ if and only if
$\phi_{\tau}(\CL_A^{\ZZ/2\ZZ})=\{0\}$.
\end{proposition}
\bp
%Let us check that the above definition is consistent.
Define $\phi_{\tau}$ for the embedded possibly non-connected curves $x$ in $\partial A$ by $\phi_{\tau}(x)=0$ if and only if $x$ is $\tau$-bounding, (that is such that $\tau$ would extend to a possibly non-oriented connected surface with framed boundary $x$).

With this definition that is consistent with the first part of the above definition, $\phi_{\tau}$ is additive under disjoint union. This is easy to see if one of the considered embedded curve is $\tau$-bounding, and this additivity  property
is preserved by a trivialisation change.

Then this definition only depends on the class of $x$ in $H_1(\partial A;\ZZ/2\ZZ)$.
Indeed let $x$ be an embedded (possibly non-connected) curve in $\partial A$
and let $y$ be another such in $\partial A \times \{-1\}$ that is homologous to $x$ mod $2$. Then there exists a framed (possibly non-orientable) cobordism between $x$ and $y$ in $\partial A \times [-1,1]$, and it is easy to see that 
$x$ is $\tau$-bounding if and only if $y$ is $\tau$-bounding.

Let us check that $\phi_{\tau}$ behaves as predicted under addition.
Because we are dealing with elements of $H_1(\partial A;\ZZ/2\ZZ)$, we can
consider representatives of $x$ and $y$ that are disjoint or that intersect once. The first case has already been treated.
Note that both sides of the equality to be proved vary in the same way under trivialisation changes. 
When $x$ and $y$ intersect once, consider the punctured torus neighborhood of $x \cup y$, and
a trivialisation $\tau$ that restricts to the punctured torus as the direct sum of a trivialisation of the torus and the normal vector to $\partial A$. Then 
$\phi_{\tau}(x+y)=\phi_{\tau}(x)= \phi_{\tau}(y)=1$.
The last two assertions are left to the reader.
\eop

\begin{example}
\label{exabadlag}
For any $\QQ$-handlebody $A$, there exists a Lagrangian subspace $\CL^{\ZZ}$ of $(H_1(\partial A;\ZZ);\langle ., . \rangle )$, such that
$\CL_A=\CL^{\ZZ} \otimes \QQ$. However, as the following example shows,
$\CL_A^{\ZZ/2\ZZ}$ \index{NN}{LAZ2@$\CL_A^{\ZZ/2\ZZ}$} is not necessarily equal to $\CL^{\ZZ} \otimes {\ZZ/2\ZZ}$.

The $\QQ$-handlebody $A$ will be the exterior of a knot $\partial M$ in 
$S^2 \times S^1$ described below. 
Consider a Moebius band $M$ embedded in the interior of a solid torus $D^2 \times S^1$ so that
the core of the solid torus is equal to the core of $M$. 
Embed $D^2 \times S^1$ into $S^2 \times S^1=D^2 \times S^1 \cup_{\partial D^2 \times S^1} (-D^2 \times S^1)$ as the first copy. 
Let $m$ be the meridian of the knot $\partial M$ that is oriented so that $\partial M$ pierces twice $S^2 \times 1$ positively, and let $\ell$ be the parallel of $\partial M$ induced by $M$.
Then $A$ is a $\QQ$-handlebody such that
$\CL_A^{\ZZ}=\ZZ[2m]$, $\CL_A=\QQ[m]$, and $\CL_A^{\ZZ/2\ZZ}=\ZZ/2\ZZ[\ell]$.
\end{example}

\subsection{Special admissible forms and their $p_1$.}
\label{subspecad}
\index{TT}{special admissible forms}

In this subsection, we fix
\begin{itemize} 
\item a rational homology sphere $M$,
\item an embedding of a neighborhood $N(c)=[a,b] \times c \times [-1,1]$  of a framed link $c$ in $B_M$ that respects the framing of $c$, 
\item a trivialisation $\tau$ of
$T(M\setminus (\infty \cup N(c))$ that is standard near $\infty$, and  
\item a trivialisation $\tau_b$ of $T(N(c))$ such that \index{NN}{taub@$\tau_b$}
$$\tau_b=\left\{\begin{array}{ll} \tau & \mbox{over}\; \partial([a,b] \times c \times [-1,1])\setminus (\{a\} \times c \times [-1,1]) \\
\CT_c^{-1} \circ \tau  & \mbox{over}\; \{a\} \times c \times [-1,1].
\end{array} \right.$$
\end{itemize}

With these data, we will associate: 
\begin{itemize}
\item{in Notation~\ref{notformad},} an admissible $2$-form $\omega(c;\tau,\tau_b)$ \index{NN}{omegactautaub@$\omega(c;\tau,\tau_b)$}
on  $\partial C_2(M)$ that reads
$p_M(\tau)^{\ast}(\omega_{S^2})$ over $ST(M \setminus (\infty \cup N(c))$, and then
\index{NN}{znctautaub@$z_n(c;\tau,\tau_b)$} $$z_n(c;\tau,\tau_b)=z_n(\omega(c;\tau,\tau_b))=\sum_{\Gamma \in {\cal E}_n} I_{{\Gamma}}(M;\bigwedge_{i=1}^{3n}p_i^{\ast}(\omega(c;\tau,\tau_b)))[{\Gamma}]$$
as in \rconstthkktfra,
\item{in Notation~\ref{nottrivspecad},} a trivialisation \index{NN}{tauCctautaub@$\tau_{\CC}(c;\tau,\tau_b)$} $\tau_{\CC}(c;\tau,\tau_b)$ of $T(M\setminus \infty) \otimes \CC$ that reads $\tau \otimes 1_{\CC}$ over $M \setminus (\infty \cup N(c))$, and its relative
Pontryagin class \index{NN}{ponectautaub@$p_1(c;\tau,\tau_b)$} $p_1(c;\tau,\tau_b)=p_1(\tau_{\CC}(c;\tau,\tau_b))$ that is defined like in the real case.
\end{itemize}
Then in Subsection~\ref{subproofpunspec}, we shall prove 
\begin{proposition}
\label{propintpunspec}
For any rational homology sphere $M$ equipped with a framed link $c$ in $B_M$, with a trivialisation $\tau$ of
$T(M\setminus \infty)$ that is standard near $\infty$ outside $N(c)$, 
and with a  trivialisation $\tau_b$ of $T(N(c))$ such that $$\tau_b=\left\{\begin{array}{ll} \tau & \mbox{over}\; \partial([a,b] \times c \times [-1,1])\setminus (\{a\} \times c \times [-1,1]) \\
\CT_c^{-1} \circ \tau  & \mbox{over}\; \{a\} \times c \times [-1,1].
\end{array} \right.$$
for any trivialisation $\tau_M$ of $M \setminus \infty$ that is standard near $\infty$,\index{NN}{znctautaub@$z_n(c;\tau,\tau_b)$}\index{NN}{ponectautaub@$p_1(c;\tau,\tau_b)$}
$$z_n(c;\tau,\tau_b) - z_n(\tau_M) = \frac{p_1(\tau_M)-p_1(c;\tau,\tau_b)}{4}\delta_n.$$
%%%SIGNE
In particular, \index{NN}{ZM@$Z(M)$} $Z(M)=Z(M;\omega(c;\tau,\tau_b)) \exp(\frac{p_1(\tau(c;\tau,\tau_b))}{4}\xi).$
%%%SIGNE
%for an integer $p_1(\omega_M)$ that is independent of $\tau_M$. 
\end{proposition}

The forms $\omega(c;\tau,\tau_b)$ are called {\em special admissible forms.\/}
We shall see that they satisfy Lemma~\ref{lemintsur} (see Lemma~\ref{lemintsurbis}). 

The two main ideas that make our constructions work are the following ones.
\begin{enumerate}
\item This neighborhood will be filled in by gadgets that factor through
the projection 
$p_c:[a,b] \times c \times [-1,1] \longrightarrow [a,b] \times [-1,1]$.
\item Over $N(c)$, $\omega(c;\tau,\tau_b)$ \index{NN}{omegactautaub@$\omega(c;\tau,\tau_b)$}will be the average of two forms coming from genuine trivialisations. 
\end{enumerate} 

\begin{notation}
\label{notformad}
Let $\varepsilon>0$ be a small positive number and let \index{NN}{F@$F$} $F$ be a smooth map such that
$$\begin{array}{llll}F:&[a,b] \times [-1,1] &\longrightarrow &SO(3)\\
& (t,u) & \mapsto & \left\{\begin{array}{ll}\mbox{Identity}&
 \mbox{if}\; |u|>1-\varepsilon\\
R_{\theta(u)} &  \mbox{if}\; t<a+\varepsilon\\
R_{-\theta(u)} &  \mbox{if}\; t>b-\varepsilon\end{array}\right.\end{array}$$
where $\theta$ has been defined in Definition~\ref{deftwicur}.
The map $F$ extends to $[a,b] \times [-1,1]$  because its restriction to the boundary
is trivial in $\pi_1(SO(3))$.

Let \index{NN}{Fctaub@${F}(c,\tau_b)$} ${F}(c,\tau_b)$ be defined on $ST(N(c)) \stackrel{\tau_b}{=}   [a,b] \times   c \times [-1,1] \times S^2$ as follows
$$\begin{array}{llll} {F}(c,\tau_b): &  [a,b] \times   c \times [-1,1] \times S^2 & \longrightarrow & S^2\\
&(t,\sigma,u, v) & \mapsto &  F(t,u)(v).\end{array}$$
Define the closed two-form \index{NN}{omegactaub@$\omega(c,\tau_b)$} $\omega(c,\tau_b)$ on $ST([a,b] \times c \times [-1,1])$ as \begin{equation}
\label{eqspecad}
\omega(c,\tau_b)=\frac{\pi_{S^2}(\CT_c \circ \tau_b)^{\ast}(\omega_{S^2}) +{F}(c,\tau_b)^{\ast}(\omega_{S^2})}2.
\end{equation}
Set \index{NN}{omegactautaub@$\omega(c;\tau,\tau_b)$}
$$\omega(c;\tau,\tau_b)=\left\{ \begin{array}{ll}p_M(\tau)^{\ast}(\omega_{S^2})& \mbox{on}\; ST(M \setminus (N(c) \cup \infty))\\
\omega(c,\tau_b)
&\mbox{on}\;ST(N(c)).
\end{array}\right.$$
\end{notation}
Extend $\omega(c,\tau_b)$ on $\partial C_2(M) \setminus ST(B_M)$ like in the case
of fundamental forms.

Observe that $\omega(c,\tau_b)$ is the average of two forms corresponding to trivialisations.
The definition of $\omega(c;\tau,\tau_b)$ is consistent because using Lemma~\ref{lemprimtwist}, we see that:  \index{NN}{omegactaub@$\omega(c,\tau_b)$}
$$\omega(c,\tau_b)= \left\{\begin{array}{ll} \pi_{S^2}( \tau_b=\tau)^{\ast}(\omega_{S^2}) &\mbox{on}\; ST(\{b\} \times c \times [-1,1])\\
\pi_{S^2}(\tau_b=\tau)^{\ast}(\omega_{S^2}) &\mbox{on}\; ST\left([a,b] \times c \times \{-1,1\}\right)\\
\pi_{S^2}(\CT_c \circ \tau_b=\tau)^{\ast}(\omega_{S^2}) &\mbox{on}\; ST(\{a\} \times c \times [-1,1])\end{array} \right.$$

Of course, $\omega(c;\tau,\tau_b)$ \index{NN}{omegactautaub@$\omega(c;\tau,\tau_b)$}depends on many choices. However, we shall see that
$z_n(\omega(c;\tau,\tau_b))$ only depends on the isotopy class of the framed link $c$, and on the homotopy classes of $\tau$ and $\tau_b$.

Bundle isomorphisms $\phi$ from $T(M \setminus \infty)$ to itself over a diffeomorphism $\overline{\phi}$ that is the
identity on $(M \setminus (B_M \cup \infty))$ induce isomorphisms of 
$ST(M \setminus \infty)$ that trivially extend to $\partial C_2(M)$.
The pull-back $\phi^{\ast}(\omega(c;\tau,\tau_b))$ of a special admissible form $\omega(c;\tau,\tau_b)$
under such an isomorphism is a special admissible form $\omega(\overline{\phi}^{-1}(c);\tau \circ \phi,\tau_b \circ \phi)$.

\begin{lemma}
\label{lemrelspecad}
If $\phi_1$ is a bundle isomorphism from $T(M \setminus \infty)$ to itself that is the
identity over $(M \setminus (B_M \cup \infty))$ and that is isotopic to the identity among these, then 
$z_n(\omega(\overline{\phi_1}^{-1}(c);\tau \circ \phi_1,\tau_b \circ \phi_1)) = z_n(\omega(c;\tau,\tau_b))$.\\
Furthermore,
$z_n(\omega(c;\tau,\tau_b))$\index{NN}{znctautaub@$z_n(c;\tau,\tau_b)$} is independent of the choices of $F$,
$c$ inside its homotopy class, $\tau_b$ and $\tau$ in their homotopy classes.
\end{lemma}
\bp For the first assertion, consider
the isotopy $$\phi: [0,1] \times ST(M) \longrightarrow ST(M).$$
Then \inconstpropkeytriv ~tells that 
$$z_n(\phi_1^{\ast}(\omega(c;\tau,\tau_b))) - z_n(\omega(c;\tau,\tau_b))= z_n([0,1] \times ST(B_M);\phi^{\ast}(\omega(c;\tau,\tau_b))).$$
The right-hand side vanishes thanks to \inconstlemtranspform.
This proves that when $\tau$ and $\tau_b$ vary in their homotopy classes so that they always satisfy 
$$\tau_b=\left\{\begin{array}{ll} \tau & \mbox{over}\; \partial([a,b] \times c \times [-1,1])\setminus (\{a\} \times c \times [-1,1]) \\
\CT_c^{-1} \circ \tau  & \mbox{over}\; \{a\} \times c \times [-1,1],
\end{array} \right.$$ $z_n(c;\tau,\tau_b)$ remains the same.

Similarly, the choice of the tubular neighborhood of $c$ in $\Sigma$ does not matter because it can be realised by an isotopy of $M$ that induces a bundle isomorphism isotopic to the identity on $ST(B_M)$.

Now, since $\pi_2(SO(3))$ is trivial, two maps $F=G_0$ and $G_1$ that satisfy
the hypotheses of Notation~\ref{notformad} are homotopic by a homotopy 
$$G:[0,1] \times [a,b] \times [-1,1] \longrightarrow SO(3).$$
Extend $\omega(G)$ that satisfies 
$$\omega(G)=\left\{\begin{array}{ll} \omega(G_1) & \mbox{on}\; \{1\} \times ST(B_M) \\
\omega(G_0) & \mbox{on}\; \{0\} \times ST(B_M) \\
\pi_{S^2}^{\ast}(\tau)(\omega_{S^2}) & \mbox{on}\; [0,1] \times ST(B_M \setminus
N(c)),
\end{array} \right.$$
by using Formula~\ref{eqspecad} on $[0,1] \times [a,b] \times c \times [-1,1] \times S^2$. Then according to \inconstpropkeytriv, $z_n(\omega(G_1)) -z_n(\omega(G_0))=z_n([0,1] \times ST(B_M);\omega(G))$, and since $\omega(G)$ pulls-back under a bundle morphism onto $[0,1] \times [a,b] \times [-1,1] \times S^2$ over $[0,1] \times N(c)$, then $$z_n([0,1] \times ST(B_M);\omega(G))=0,$$ thanks to \inconstlemtranspform. \eop

\begin{lemma}
\label{lemintsurbis} 
Let $S$
be an oriented surface of $M$ whose boundary does not meet $N(c)$, then
$\int_{S(n)}\omega(c;\tau,\tau_b)$ only depends on the topology of $S$ and on the restriction of $\tau$ on $\partial S$,
and if the positive normal $n$ of $S$ is the first vector of the trivialisation $\tau$ on $\partial S$, then \index{NN}{omegactautaub@$\omega(c;\tau,\tau_b)$}
\index{NN}{chiTSv@$\chi(TS; .)$}
$$\int_{S(n)}\omega(c;\tau,\tau_b)=\frac12\chi(TS; \tau^{-1}(.,e_2)_{|\partial S}).$$
\end{lemma}
\bp 
First isotope $S$ so that it meets $c$ along meridian squares $D(x_i)=[a,b] \times \{x_i\}\times[-1,1]$ of $c$.
Without changing the sides of the equality to be shown, perform a bundle isomorphism of $ST(M \setminus \infty)$ over $1_M$ that is isotopic to the identity and supported near $\partial D(x_i)$ so that the normal $n$ to the square $D(x_i)$ becomes the first vector of $\tau$ around $\partial D(x_i)$.
Let $D_i$ be a disk inside $D(x_i)$ with a smooth boundary outside $[a+\varepsilon,b-\varepsilon] \times \{x_i\}\times[-1+\varepsilon,1-\varepsilon]$ such that $\tau(n)=e_1$ on $\partial D_i$.
Set $\omega_M=\omega(c;\tau,\tau_b)$.
Then we have the following sublemma.
\begin{sublemma}
$\int_{D_i(n)}\omega_M=\frac12\chi(TS; \tau^{-1}(.,e_2)_{|\partial D_i}).$
\end{sublemma}
\bp Indeed, on $D_i$, $\omega_M$ is the average of two forms corresponding to 
the trivialisations $\CT_c \circ \tau_b$ and ${F}(c,\tau_b) \circ \tau_b$.
For both of these, $n$ is the first vector of the trivialisation on $\partial D_i$, and 
Lemma~\ref{lemintsur} tells us that $2\int_{D_i(n)}\omega_M$ is the average
of the obstructions to extend the second vector 
of these trivialisations to $D_i$. These obstructions vary like the degrees 
of the induced maps from $\partial D_i$ to the fiber $S^1$ of $S(TD_{i|\partial D_i})$ equipped with an arbitrary trivialisation.
We conclude because the degree of the map induced by $\tau$ is the average
of the degrees of the map induced by $\CT_c \circ \tau_b$ and the map induced by ${F}(c,\tau_b) \circ \tau_b$.
\eop

Back to the proof of Lemma~\ref{lemintsurbis},
use Lemma~\ref{lemintsur} and the additivity of both sides of the equality under gluing to treat the case when the positive normal $n$ of $S$ is the first vector of the trivialisation $\tau$ on $\partial S$.
For the other case, let $ [-1,0] \times \partial S$ be a collar of $\partial S =\partial S \times \{0\}$ inside $S$. After an arbitrary isotopy of $\tau$  around $ [-1,0] \times \partial S$ supported away from $\partial S$, the positive normal $n$ of $S$ is the first vector of the trivialisation $\tau$ on $\partial S \times \{-1\}$. Then, since $\int_{S(n)}\omega_M=\int_{(S \setminus ]-1,0] \times \partial S)(n)}\omega_M +\int_{([-1,0] \times \partial S)(n)}\omega_M$, according to the first case, $\int_{S(n)}\omega_M$ only depends on our arbitrary isotopy, on the topology of $S$ and on the restriction of $\tau$ on $\partial S$.
%After a bundle isomorphism $T(M \setminus \infty)$ 
%that is isotopic to the identity and that is the identity near
%$\partial S(n)$, we may assume that $n=\tau_b^{-1}(e_1)$
%on the boundaries of these squares.
%Then it is enough to prove the lemma when $S$ is one of these squares. 
\eop

\begin{notation}
\label{nottrivspecad}
Let $F_U$ \index{NN}{FU@$F_U$} be a smooth map such that
$$\begin{array}{llll}F_U:&[a,b] \times [-1,1] &\longrightarrow &SU(3)\\
& (t,u) & \mapsto & \left\{\begin{array}{ll}\mbox{Identity}&
 \mbox{if}\; |u|>1-\varepsilon\\
R_{\theta(u)} &  \mbox{if}\; t<a+\varepsilon\\
\mbox{Identity} &  \mbox{if}\; t>b-\varepsilon.\end{array}\right.\end{array}$$
$F_U$ extends to $[a,b] \times [-1,1]$  because $\pi_1(SU(3))$ is trivial.
Define the trivialisation $\tau_{\CC} =\tau_{\CC}(c;\tau,\tau_b)$ \index{NN}{tauCctautaub@$\tau_{\CC}(c;\tau,\tau_b)$} of $T(M\setminus \infty) \otimes \CC$ as follows.
\begin{itemize}
\item On $T(M\setminus (\infty \cup N(c)))$, $\tau_{\CC} =\tau \otimes 1_{\CC}$,
\item Over $[a,b] \times c \times [-1,1]$, 
$\tau_{\CC}(t,\gamma,u;v) =F_U(t,u)(\tau_b \otimes 1_{\CC})(t,\gamma,u;v)$.
\end{itemize}
(Here, as often $\tau_b$ that is valued in $N(c) \times \RR^3$ is identified
with $\pi_{\RR^3} \circ \tau_b$, and $\tau_{\CC}$ is identified with $\pi_{\CC^3} \circ \tau_{\CC}$.)
Since $\pi_2(SU(3))$ is trivial, the homotopy class of $\tau_{\CC}$
is well-defined. Since $\tau_{\CC}= \tau_M \otimes 1_{\CC}$ outside $B_M$, the definition of $p_1$ for trivialisations that are standard near $\infty$ extends to this kind of trivialisations, hence $\tau_{\CC}$ has a well-defined integral $p_1$ defined as in 
\inconstsubpont.
%%% BLOC
Set $$p_1(c;\tau,\tau_b)=p_1(\tau_{\CC}(c;\tau,\tau_b)).$$
\end{notation}
%\end{definition}

\subsection{Proof of Proposition~\ref{propintpunspec}.}
\label{subproofpunspec}
Of course it is enough to prove this proposition for some trivialisation $\tau_M$ of $M \setminus \infty$ that is standard near $\infty$,
thanks to \inconstproppontproptrois.
The last sentence of the proposition follows from the previous one by the arguments given
at the end of \inconstsubsketchpkkt.
%%2.2
\medskip

\noindent{\em Fixing the main notation for the proof.\/}\\
Let $A$ be a rational homology handlebody that meets $N(c)$ along
$\{a\} \times c \times [-1,1]$.
($A$ can be the genuine handlebody obtained by first thickening $N(c)$ to an embedding of 
$[a-5,b] \times c \times [-2,2]$ and then construct $A$ as a regular neighborhood of a connected graph whose loops are the connected components of $\{a-2\} \times c \times \{0\}$.)
Let $[a,b] \times \partial A$ be a collar of $\partial A$ in the closure of 
$B_M \setminus A$ such that the notation  $[a,b] \times c \times [-1,1]$
is consistent.
 
Let $D$ be a disk of $\partial A$ that does not meet the tubular neighborhood of $c$. Let $D^3$ be a topological ball of $M$ that contains $M \setminus B_M$, that intersects
$[a,b] \times \partial A$ as $[a,b] \times D$, and $A$ along
$\{a\} \times D$.
Choose a trivialisation $\tau_M$ of $M \setminus \infty$ that is standard near $\infty$ so that $\tau$ \index{NN}{tau@$\tau$} and $\tau_M$ \index{NN}{tauM@$\tau_M$} coincide on the complement
of $\{\infty\}$ in the ball $D^3$.
Set $B^D_M=M \setminus \mbox{Int}(D^3)$, \index{NN}{BMD@$B^D_M$} $\Sigma=\partial A \setminus \mbox{Int}(D)$, \index{NN}{Sigma@$\Sigma$} and
$$B = B^D_M \setminus (A \cup ([a,b[ \times \Sigma)) $$
Then $B^D_M$ is the union of three rational 
homology handlebodies $A$, $B$ and $[a,b] \times \Sigma$.

There exists a curve $c_B$ \index{NN}{cB@$c_B$} of $\Sigma$ that is transverse to $c$, such that the restriction of $\tau$ to $\{b\} \times \Sigma$ is homotopic to $\CT_{c_B} \circ \tau_{M|\{b\} \times\Sigma}$. Homotope $\tau_M$ so that \index{NN}{taub@$\tau_b$}
\begin{itemize}
\item $\tau=\CT_{c_B} \circ \tau_{M}$ on $[a,b] \times (\Sigma \setminus c \times [-1,1])$, and
\item $\tau_{b}=\CT_{c_B} \circ \tau_{M}$ on $[a,b] \times c \times [-1,1]$.
\end{itemize}
 
Let $c_A$ \index{NN}{cA@$c_A$} be the curve obtained from $c \cup c_B$ by replacing the neighborhoods \begin{pspicture}[.2](0,0)(.6,.45)
\psline{->}(0.05,0)(.45,.4)
\psline{->}(.45,0)(0.05,.4)
\end{pspicture} of the double points by \begin{pspicture}[.2](0,0)(.6,.45)
\pscurve{->}(0.45,0)(.3,.2)(.45,.4)
\pscurve{->}(.05,0)(.2,.2)(0.05,.4)
\end{pspicture} according to the orientation. Then $c_A$ represents $(c+c_B)$ in $H_1(\Sigma)$, and the restriction of $\tau_A=\tau_{|A}$ to $\{a\} \times \Sigma$ is homotopic to $\CT_{c_A} \circ \tau_{M|\{a\} \times \Sigma}$. 

Set $\omega_M=\omega(c;\tau,\tau_b)$. 

\medskip

\noindent{\em Sketch of the proof of Proposition~\ref{propintpunspec}.}\\
%of  
Consider
the part $$B_0=\partial ([0,1] \times B^D_M) \cup ([0,1] \times \{a,b\} \times \Sigma)$$ of 
the base $[0,1] \times B^D_M$ of the $S(\RR^3)$-bundle $E=[0,1] \times ST(B^D_M)$. We shall fix a closed two-form $\omega$ over $B_0$ 
such that $$\omega=\left\{ \begin{array}{ll} \omega_M & \mbox{on}\;\{0\} \times ST(B^D_M)\\
\omega(\tau_M)  & \mbox{on}\;\{1\} \times ST(B^D_M) 
\end{array} \right.$$
and such that any extension $\omega$ of $\omega$ to $E$ will satisfy 
$$z_n(\tau_M) - z_n(\omega_M)=z_n(E; \omega).$$
Simultaneously, we shall fix a related trivialisation $\tau_0$ of $T([0,1] \times B^D_M) \otimes \CC$ over $B_0$, that coincides with the trivialisation
associated to $\omega_M$ and $\tau_M$ on $\{0,1\} \times  B^D_M$.
Now, for each closure $C$ of a connected component of $([0,1] \times B^D_M) \setminus B_0$, that is for $C=[0,1] \times A$, $[0,1] \times B$ or  $[0,1] \times [a,b] \times \Sigma$, we are going to prove that:\\
$\omega$ extends as a closed form $\omega$ over $C$,
and $z_n(E_{|C};\omega) =-\frac{p_1(C;\tau_{0|\partial C})}{4} \delta_n$ 
%%% SIGNE ???
for each $C$.
Of course, once these goals are achieved (in Lemmas~\ref{lemcasa}, \ref{lemcasb} and \ref{lemdifsurf}), the proof will be finished.

\medskip

\noindent{\em Definition of $\omega$ over 
$B_0$.}\\
Define the bundle isomorphisms over the Identity of $[0,1] \times \{d\} \times \Sigma$ \index{NN}{Qcal@${\cal Q}$}
$$\begin{array}{llll}{\cal Q}: &ST([0,1] \times \{d\} \times \Sigma) & \longrightarrow & [0,1] \times ST(M)_{|\{d\} \times \Sigma}\\
&T((t,d) \times 1_{\Sigma})( v \in T_{\sigma}\Sigma) & \mapsto & (t;T(\{d\} \times 1_{\Sigma})(v)) \\
& \frac{\partial}{\partial u}(t+u,d,\sigma)_{u=0} & \mapsto & (t;\frac{\partial}{\partial u}(d+u,\sigma)_{u=0}).\end{array}$$
Define $\omega_{{\cal Q}}(c_B,\tau_M )$ \index{NN}{omegaQcalB@$\omega_{{\cal Q}}(c_B,\tau_M )$} on $[0,1] \times ST(M)_{|\{b\} \times \Sigma}$ and $\omega_{{\cal Q}}(c_A,\tau_M)$ \index{NN}{omegaQcalA@$\omega_{{\cal Q}}(c_A,\tau_M )$} on $[0,1] \times ST(M)_{|\{a\} \times \Sigma}$ by:
$$\omega_{{\cal Q}}(c_B,\tau_M)=\left\{ \begin{array}{ll} 
\left({\cal Q}^{-1}\right)^{\ast}
\left(\omega(c_B,\tau_M\circ {\cal Q}) \right)& \mbox{on}\;[0,1] \times ST(M)_{|\{b\} \times c_B \times [-1,1]}\\
\omega(\tau_M)& \mbox{on}\;[0,1] \times ST(M)_{|\{b\} \times (\Sigma \setminus c_B \times [-1,1])},
\end{array} \right.$$
$$\omega_{{\cal Q}}(c_A,\tau_M)=\left\{ \begin{array}{ll} 
\left({\cal Q}^{-1}\right)^{\ast}
\left(\omega(c_A,\tau_M\circ {\cal Q}) \right)& \mbox{on}\;[0,1] \times ST(M)_{|\{a\} \times c_A \times [-1,1]}\\
\omega(\tau_M)& \mbox{on}\;[0,1] \times ST(M)_{|\{a\} \times (\Sigma \setminus c_A \times [-1,1])},
\end{array} \right.$$
using the notation of Equation~\ref{eqspecad}. Set 
$$\omega=\left\{ \begin{array}{ll} \omega_M & \mbox{on}\;\{0\} \times ST(M \setminus \infty)\\
\omega(\tau_M)  & \mbox{on}\;\{1\} \times ST(M \setminus \infty) 
\\
\omega_{{\cal Q}}(c_B,\tau_M)& \mbox{on}\;[0,1] \times ST(M)_{|\{b\} \times \Sigma}\\
\omega_{{\cal Q}}(c_A,\tau_M)& \mbox{on}\;[0,1] \times ST(M)_{|\{a\} \times \Sigma}.
\end{array} \right.$$

\begin{pspicture}[.4](0,-1)(12,7) 
\psline(3,0)(9,0)(9,7)(3,7)(3,0)
\psline(3,2)(9,2)
\psline(3,5)(9,5)
\rput[b](6,5.1){$\omega_{{\cal Q}}(c_B,\tau_M)$}
\rput[t](6,1.9){$\omega_{{\cal Q}}(c_A,\tau_M)$}
\rput[l](9.1,1){$\omega(\tau_M)$}
\rput[l](9.1,3.5){$\omega(\tau_M)$}
\rput[l](9.1,6){$\omega(\tau_M)$}
\rput[r](2.9,6){$\omega(\tau_B=\tau_{|B})$}
\rput[r](2.9,1){$\omega(\tau_A=\tau_{|A})$}
\rput[r](2.9,3.5){$\tilde{\omega}(c,\tau_b=\CT_{c_B}\circ \tau_M)$}
\rput(6,1){$[0,1] \times ST(A)$}
\rput(6,6){$[0,1] \times ST(B)$}
\rput(6,3.5){$[0,1] \times ST([a,b] \times \Sigma)$}
%\psset{xunit=.5cm,yunit=.5cm} 
\end{pspicture}

Here, $\tilde{\omega}(c,\tau_b)$ is the following slight modification of ${\omega}(c;\tau,\tau_b)$ (allowed by Lemma~\ref{lemrelspecad}).
$$\omega=\left\{ \begin{array}{ll} {\omega}(c;\tau,\tau_b) & \mbox{on}\; ST([a+\varepsilon/2,b] \times \Sigma)\\
\pi_{S^2}(H \circ  \CT_{c_A} \circ \tau_M)^{\ast}(\omega_{S^2})  & \mbox{on}\;ST([a,a+\varepsilon/2] \times \Sigma)
\end{array} \right.$$
where $$H:[a,a+\varepsilon/2] \times \Sigma \longrightarrow SO(2) \subset SO(3)$$ is a homotopy supported near the intersection points of $c$ and $c_B$, valued in the subgroup $SO(2)$ of $SO(3)$ that fixes $e_1$,
such that $H(\{a\} \times \Sigma)=\{\mbox{Identity}\}$, and 
$H_{a+\varepsilon/2} \circ \CT_{c_A}=\CT_c\CT_{c_B}$.

\medskip

\noindent{\em Definition of the trivialisation $\tau_0$ of $T([0,1] \times B^D_M)\otimes \CC$
over $B_0$.}\\
Wherever $\omega$ is associated to a trivialisation $\tau$ of $T(B^D_M)$, $\tau_0$ is its natural stabilisation $S(\tau)$ obtained by mapping the unit tangent vector to $T([0,1] \times \{x\})$ to the first basis vector $E_1$ of $\CC^4$, and the $i^{\mbox{th}}$ basis vector of $\RR^3$ to the $(i+1)^{\mbox{th}}$ basis vector of $\CC^4$.
On the remaining parts of $B_0$ that are of the form $c \times D^2$, with
a trivialisation $\tau_b$ of $T(B_M)$ over $c \times D^2$ involved (possibly via ${\cal Q}$),
$$\tau_0(\gamma \in c,d \in D^2;v)= \tilde{F}_U(d)(S(\tau_b)(\gamma,d;v))$$
 where $\tilde{F}_U$ is a fixed map from $D^2$ to the stabilizer $SU(3)$ of $E_1$ in $SU(4)$.

\begin{lemma}
\label{lemexist}
The form $\omega$ extends as a closed two-form on $[0,1] \times ST(A)$, on $[0,1] \times ST(B)$ and on $[0,1] \times ST([a,b]\times \Sigma)$.
\end{lemma}
\bp
Let us first treat the case of $[0,1] \times ST(A)$.
We know $\omega$ on $\partial ([0,1] \times A) \times S^2$, and it is sufficient to prove that the integral of $\omega$ vanishes on the
kernel of the map induced by the inclusion $H_2(\partial ([0,1] \times A) \times S^2) \longrightarrow H_2([0,1] \times A \times S^2)$. This kernel is generated by the $\partial ([0,1] \times S(a_i)(n))$, for $i=1, \dots , g(\Sigma)$, where the $a_i$ are simple curves that generate $\CL_A$, $S(a_i)$ is a surface whose boundary is made of 
$k_i$ copies of $a_i$ and $S(a_i)(n)$ is the section of $ST(A)_{|S(a_i)}$ given by its positive normal that belongs to $T(\{x\} \times \partial A)$ on a neighborhood $[a-1,a+1] \times \partial A$ of $\partial A$. 

$$\int_{\partial ([0,1] \times S(a_i)(n))}\omega= \int_{S(a_i)(n)} \omega(\tau_M) - \int_{0 \times S(a_i)(n) \cup [0,1] \times \partial S(a_i)(n)} \omega .$$
Lemma~\ref{lemintsurbis} tells that the first integral is a well-determined function of the restriction of $\tau_M$
to $\partial S(a_i)(n)$ and of the topology of $S(a_i)(n)$.

Now, $\tilde{A}=(\{0\} \times A) \cup_{\partial A} ([0,1] \times \partial A)$ is equipped
with a smooth structure in a standard way so that the $S^2$-bundle over it 
$ST(\{x\} \times M)$ is identified with its unit tangent bundle with the help
of the bundle isomorphism ${\cal Q}$ over $[0,1] \times \partial A$. The smoothing
makes $$\tilde{S}(a_i)(n)=0 \times S(a_i)(n) \cup [0,1] \times \partial S(a_i)(n)$$ smooth and leaves the normal to the surface unchanged. On $\tilde{A}$, the form $\omega$ is a special admissible form with respect to
a trivialisation $\tilde{\tau}$ of $T(\tilde{A} \setminus N(c_A))$ and a trivialisation $\tilde{\tau}_b$ of $T(N(c_A))$. Now,  $(\tilde{A}, \tilde{S}(a_i)(n),\tilde{\tau},\tilde{\tau}_b)$ is isomorphic to some $(A, {S}(a_i)(n),\tau^{\prime},\tau^{\prime}_b)$ such that the restrictions of $\tau^{\prime}$ and $\tau_M$ to $\partial {S}(a_i)(n)$ coincide.  Therefore, Lemma~\ref{lemintsurbis} applies
to identify $\int_{\tilde{S}(a_i)(n)} \omega$ to $\int_{S(a_i)(n)} \omega(\tau_M)$.
Thus, $\omega$ extends as a closed 2-form on $[0,1]\times ST(A)$.

Similarly, for any surface $S(b_i)$ of $B$ whose boundary is in $\{b\} \times \Sigma$, Lemma~\ref{lemintsurbis} applies to identify 
$\int_{\{0\} \times S(b_i)(n)}\omega$
to $\int_{- [0,1] \times \partial S(b_i)(n) \cup \{1\} \times S(b_i)(n)}\omega$,
and allows us to prove that $\omega$ extends as a closed 2-form on $[0,1]\times ST(B)$.

Now, for any curve $\gamma$ of $\Sigma$, Lemma~\ref{lemintsurbis} also applies to prove 
$$\int_{\{0\} \times [a,b] \times \gamma \cup [0,1] \times \{b\} \times \gamma} \omega=\int_{[0,1] \times \{a\} \times \gamma \cup\{1\} \times [a,b] \times \gamma } \omega,$$
and we easily deduce from this fact that $\omega$ extends as a closed 2-form on $[0,1]\times ST([a,b] \times \Sigma)$.
\eop

\begin{lemma}
\label{lemsimphom}
Let $K$ be a rational homology handlebody, and let $\tau$ be a trivialisation
of $TK$ over $\partial ([0,1] \times K)$. Let $\omega$ denote the
associated two-form on $\partial([0,1] \times ST(K))$ and let $\tau_{\CC}$ be the associated trivialisation of $T([0,1] \times K) \otimes \CC$.
Then $\omega$ extends to $[0,1] \times ST(K)$ and for any closed extension of $\omega$ to $[0,1] \times ST(K)$,
$$z_n({[0,1] \times ST(K)};\omega)= -\frac14 p_1([0,1] \times K;\tau_{\CC|\partial ([0,1] \times K)} ) \delta_n.$$
%%%SIGNE ???
\end{lemma}
\bp
The existence of $\omega$ is shown as in the above proof of Lemma~\ref{lemexist}.
By \inconstlemdepboun, %%%,
since the restriction injects $H^2([0,1] \times K \times S^2)$ into $H^2(\partial([0,1] \times K) \times S^2)$ and maps $H^1([0,1] \times K \times S^2)$ onto $H^1(\partial([0,1] \times K) \times S^2)$, 
$z_n({[0,1] \times ST(K)};\omega)$ only depends on the values of $\omega$
on the boundary of $[0,1] \times ST(K)$.

Use the restriction $\tau_1$ of $\tau$ to $\{1\} \times K$ to identify
$ST(K)$ to $K \times S^2$. Then $ \tau = G \circ \tau_1$ for some map
$G:\partial ([0,1] \times K) \longrightarrow SO(3)$ that maps $\{1\} \times K$
to $1$. Then both sides only depend on the homotopy class of $G$ among these
maps. In particular, we may assume that $G$ maps $[0,1] \times \partial K$ 
to $1$. But in this case, it is enough to embed $K$ into a rational homology sphere $M$ where $\tau_1$ extends and to apply \inconstproppontproptroispropkeytriv. 
%%%that give the result in this case.
\eop

\begin{lemma}
\label{lemcasa}
$$z_n({[0,1] \times ST(A)};\omega)= -\frac14 p_1([0,1] \times A;\tau_{0|\partial [0,1] \times A} ) \delta_n.$$
\end{lemma}
\bp
There exists a bundle isomorphism $\psi: ST(A) \longrightarrow  ST(A)$ over the identity map of $A$ such that
$\tau_A= \tau_M \circ \psi $. \\
Let $\psi$ still denote
$\psi: \begin{array}{lll}[-1,0] \times ST(A) &\longrightarrow &[0,1] \times ST(A)\\
(t,y) & \mapsto & (t+1,\psi(y)). \end{array}$\\
Extend $\omega$ on $[-1,1] \times ST(A)$ by $\psi^{\ast}(\omega_{|[0,1] \times ST(A)})$ on $[-1,0] \times ST(A)$.
Then according to \inconstlemtranspform,
$$z_n({[-1,0] \times ST(A)};\omega)=z_n([0,1] \times ST(A);\omega)=\frac12 z_n({[-1,1] \times ST(A)};\omega).$$

Again, 
$z_n({[-1,1] \times ST(A)};\omega)$ only depends on the values of $\omega$
on the boundary of $[-1,1] \times ST(A)$ that are 
$$\omega=\left\{\begin{array}{ll}
\omega(\tau_M \circ \psi^2) & \mbox{on} \; \{-1\} \times ST(A)\\ \omega(\tau_M) & \mbox{on} \;\{1\} \times ST(A)\\
\omega_{\cal Q}(c_A,\tau_M) & \mbox{on} \;[0,1] \times  ST(A)_{|\partial A}\stackrel{\cal Q}{=}ST([0,1] \times \partial A)\\
\omega_{\cal Q}(c_A,\tau_M \circ \psi) & \mbox{on} \;[-1,0] \times  ST(A)_{|\partial A}\stackrel{\cal Q}{=}ST([-1,0] \times \partial A). \end{array}\right.$$

Furthermore, with the data on $\partial ([-1,1] \times A)$, we again
naturally associate a trivialisation $\tau_0$ of $T([-1,1] \times A) \otimes \CC$ over $\partial ([-1,1] \times A)$ such that 
$$p_1([-1,1] \times A;\tau_{0|\partial ([-1,1] \times A)})=2p_1([0,1] \times A;\tau_{0|\partial ([0,1] \times A)}).$$

Therefore, we are left with the proof of the following equality.
\begin{equation}
\label{eqomega}
z_n({[-1,1] \times ST(A)};\omega)=-\frac{p_1([-1,1] \times A;\tau_{0|\partial ([-1,1] \times A)})}4 \delta_n.
\end{equation}
%%%SIGNE

There exists
a smooth map $G$ for a small positive number $\varepsilon>0$
such that
$$\begin{array}{llll}G:&[-1,1] \times [-1,1] &\longrightarrow &SO(3)\\
& (t,u) & \mapsto & \left\{\begin{array}{ll}\mbox{Identity}&
 \mbox{if}\; |u|>1-\varepsilon\\
R_{2\theta(u)} &  \mbox{if}\; t<-1+\varepsilon\\
\mbox{Identity} &  \mbox{if}\; t>1-\varepsilon.\end{array}\right.\end{array}$$

Define $\tilde{G}$ as the Identity map over $[-1,1] \times (\partial A \setminus
c_A \times]-1,1[)  \times \RR^3$ and by
$$\tilde{G}((t,\sigma,u;v) \in [-1,1] \times
c_A \times]-1,1[  \times \RR^3)=(t,\sigma,u;G(t,u)(v)).$$

Thanks to Lemma~\ref{lemsimphom}, 
%(and to \inconstlemtranspform),
in order to prove Equality~\ref{eqomega}, it is enough 
to prove that changing the value of $\omega$ on $[-1,1] \times ST(\partial A)$ 
into $\omega(\tilde{G} \circ \tau_M)$ does not change the left-hand side of the equality
and that changing the trivialisation $\tau_0$ of $T([-1,1] \times A)\otimes \CC$ into $\tilde{G} \circ \tau_M \otimes 1_{\CC}$
over $[-1,1] \times \partial A$ does not change its right-hand side.

To do this, first assume without loss that over a collar $[-1,1] \times [a-1,a] \times
\partial A$ everything behaves as a product by $[a-1,a]$, and use this to  
extend both $\omega$ as $\pi_{[-1,1] \times
\partial A \times S^2}^{\ast}(\omega)$ on this collar, and the trivialisation $\tau_0$.
In particular, $\omega$ and $\tau_0$ are associated to the trivialisation $\tau_M$
over $[-1,1] \times [a-1,a] \times (\partial A \setminus (c_A \times [-1,1]))$,
and over $\{1\} \times A$. Over $[-1,1] \times  [a-1,a] \times c_A \times [-1,1]$, $\pi_{[-1,1] \times
\partial A \times S^2}^{\ast}(\omega)$ factors through $[-1,1] \times [-1,1]\times S^2$. We are going
to modify $\omega$ on $]-1,1[ \times  [a-1,a] \times c_A \times [-1,1] \times S^2$ by a form that 
still factors through the bundle projection onto $]-1,1[ \times  [a-1,a] \times [-1,1] \times S^2$, and therefore without changing $z_n([-1,1] \times A;\omega)$
(thanks to \inconstlemtranspform)
so that $\omega=\omega(\tilde{G} \circ \tau_M)$ over $[-1,1] \times \{a\} \times \partial A$. Then the associated trivialisation will read $\psi(K) \circ \tau_M$ on $\partial
([-1,1] \times  [a-1,a] \times [-1,1])$ for some fixed 
$$K:\partial([-1,1] \times  [a-1,a] \times [-1,1]) \longrightarrow SO(3)$$
Since $\pi_2(SU(4))=0$, $K$ will extend to $SU(4)$, and such an operation will not change $p_1([-1,1] \times A;\tau_0(\omega)_{\partial([-1,1] \times A)})$ either. 

Change $\omega$ over $[-1,1] \times  \{a\}\; (\times c_A) \times [-1,1]$
into $\omega(\tilde{G} \circ \tau_M)$.
In order to achieve our goal, it is enough to see that $\omega$ that is defined on the boundary
of $[-1,1] \times  [a-1,a] \times [-1,1] \times S^2$ extends to a closed form
on the whole space.
To do that, it is enough to prove that
$$\int_{\partial([-1,1] \times  [a-1,a] \times [-1,1]) \times v}\omega=0$$ for some $v \in S^2$. This integral equals $$\int_{\partial([-1,1] \times  [a-1,a]) \times [-1,1] \times v}\omega.$$
It  depends neither on $c_A$ nor on $\tau_M$. Assume that there exists a closed curve $d$ that intersects $c_A \times [-1,1]$ along $ \{\sigma\} \times  [-1,1]$.
%, and that $v=\tau_M(n)$ on this curve where $n$ is the normal to $d$. 
Then 
$$\int_{\partial([-1,1] \times  [a-1,a]) \times [-1,1] \times v}\omega
=\int_{\partial([-1,1] \times  [a-1,a]) \times d \times v}\omega.$$
and this integral vanishes as in the proof of Lemma~\ref{lemexist}.

Therefore, Lemma~\ref{lemsimphom} (together with \inconstlemtranspform) 
applies to conclude the proof. 
\eop

Similarly, we get the following lemma.
\begin{lemma}
\label{lemcasb}
$$z_n({[0,1] \times ST(B)};\omega)= -\frac14 p_1([0,1] \times B;\tau_{0|\partial [0,1] \times B} ) \delta_n$$
\end{lemma}
\bp
We can also deduce this lemma from the previous one by gluing a rational
homology handlebody $A^{\prime}$ along $\{b\} \times \Sigma$ to $B$ so that both $\tau_B$ and $\CT_{c_B} \circ \tau_B$ extend to $A^{\prime}$, and 
$A^{\prime} \cup B$ is rational homology ball. 
\eop

Unfortunately, the same methods do not apply to prove the following
similar lemma over $[0,1] \times [a,b] \times \Sigma$.

\begin{lemma}
\label{lemdifsurf}
Let $\Sigma$ be an oriented compact surface with one boundary component.
Let $\tau_M$ be a trivialisation of $ST([a,b] \times \Sigma)$, let $c$, $c_A$, $c_B$ be three non-necessarily connected curves of $\Sigma$, let
$$\omega=\left\{ \begin{array}{ll} \tilde{\omega}(c,\CT_{c_B} \circ \tau_M) & \mbox{on}\;\{0\} \times ST([a,b] \times \Sigma)\\
\omega(\tau_M)  & \mbox{on}\;\{1\} \times ST([a,b] \times \Sigma) 
\\
\omega_{\cal Q}(c_B,\tau_M)& \mbox{on}\;[0,1] \times ST([a,b] \times \Sigma)_{|\{b\} \times \Sigma}\\
\omega_{\cal Q}(c_A,\tau_M )& \mbox{on}\;[0,1] \times ST([a,b] \times \Sigma)_{|\{a\} \times \Sigma}.
\end{array} \right.$$
$$z_n({[0,1] \times ST([a,b] \times \Sigma)};\omega)= -\frac14 p_1([0,1] \times [a,b] \times \Sigma;\tau_{0|\partial ([0,1] \times [a,b] \times \Sigma)} ) \delta_n .$$
\end{lemma}

\begin{remark}
Here, we have three genuine trivialisations, namely $\tau_M$, $\CT_{c_A} \circ \tau_M$
and $\CT_{c_B} \circ \tau_M$, and it may be possible that the three of them
cannot simultaneously extend to a rational homology handlebody whose boundary
contains $\Sigma$. (Indeed, according to Proposition~\ref{propuseless}, if these three trivialisations extend to a $\QQ$-handlebody $K$, $\phi_{\tau_M}(\CL_K^{\ZZ/2\ZZ})=\phi_{\CT_{c_A} \circ \tau_M}(\CL_K^{\ZZ/2\ZZ})=\{0\}$, this implies that $\langle c_A,.\rangle$ vanishes
on $\CL_K^{\ZZ/2\ZZ}$ and hence that $c_A$ belongs to $\CL_K^{\ZZ/2\ZZ}$. Similarly $c_B$ must be in $\CL_K^{\ZZ/2\ZZ}$. Thus, if the intersection of $c_A$ and $c_B$ mod 2 does not vanish, such a $K$ cannot exist.) Therefore, we have to give further arguments to prove Lemma~\ref{lemdifsurf} that will conclude the proof of Proposition~\ref{propintpunspec}.
\end{remark}

\noindent{\sc Proof of Lemma~\ref{lemdifsurf}:}
Fix an arbitrary trivialisation $\tau_M$ of $T([a,b] \times \Sigma)$.
Any other trivialisation is obtained from $\tau_M$ by a bundle isomorphism.
By \inconstlemtranspform,   
the left-hand side does not depend on $\tau_M$. The right-hand side does not depend on $\tau_M$ either for the same reason.
Hence, if Lemma~\ref{lemdifsurf} is true for $\tau_M$, it is true for any other trivialisation of $T([a,b] \times \Sigma)$.

Then for simplicity, we shall assume that $\tau_M(\frac{\partial}{\partial u}((u, \sigma) \in [a,b] \times \Sigma))=e_1$ and that $\tau_M \circ T(\Sigma \hookrightarrow \{c\} \times \Sigma)$ is independent of $c$.
In particular, $\tau_M$ makes unambiguous sense on all $[\alpha, \beta] \times 
\Sigma$.

\begin{lemma}
\label{lemdessin}
Under the hypotheses above and the assumptions of Lemma~\ref{lemdifsurf}, equip $[0,1] \times ([2,4] \times \Sigma)$ with the trivialisation $\tau_M$ of $ [2,4] \times \Sigma$.
Let $\omega_1$ be a closed form on $[0,1] \times ST( [2,4] \times \Sigma)$
such that
$$\omega_1=\left\{ \begin{array}{ll} \omega(\tau_M)  & \mbox{on}\; [0,1] \times ST([2,4] \times \Sigma)_{|\{4\} \times \Sigma} 
\\
\omega(\CT_{c_A} \circ \tau_M)  & \mbox{on}\; [0,1] \times ST([2,4] \times \Sigma)_{|\{2\} \times \Sigma} 
\\
\omega(\tau_M)  & \mbox{on}\; \{1\} \times ST([3,4] \times \Sigma)\\
\tilde{\omega}(c,\CT_{c_B} \circ \tau_M) & \mbox{on}\;\{0\} \times ST([2,3] \times \Sigma)\\
\omega(c_B,\tau_M)& \mbox{on}\;\{0\} \times ST([3,4] \times \Sigma)\\
\omega(c_A,\tau_M)& \mbox{on}\;\{1\} \times ST([2,3] \times \Sigma).
\end{array} \right.$$
Then $z_n({[0,1] \times ST([2,4] \times \Sigma)};\omega_1)=
z_n({[0,1] \times ST([a,b] \times \Sigma)};\omega)$.
Let $W=[0,1] \times [2,4] \times \Sigma$.
Let $\tau_1$ be the trivialisation of $TW \otimes \CC$
over $\partial W$ corresponding to the given trivialisations.
Then
$$p_1(W;\tau_{1|\partial W} ) =p_1([0,1] \times [a,b] \times \Sigma;\tau_{0|\partial ([0,1] \times [a,b] \times \Sigma)} )
\delta_n.$$
\end{lemma}
\begin{pspicture}[.4](1,0)(15,6) 
\psline(3,2)(7,2)(7,4)(3,4)(3,2)
\rput[b](5,4.1){$\omega_{\cal Q}(c_B,\tau_M)$}
\rput[t](5,1.9){$\omega_{\cal Q}(c_A,\tau_M)$}
\rput[l](7.1,3){$\omega(\tau_M)$}
\rput[r](2.9,3){$\tilde{\omega}(c,\CT_{c_B}\circ \tau_M)$}
\rput(5,3){$[0,1] \times [a,b] \times \Sigma \times S^2$}
\psline(10,1)(13,3)(13,5)(10,3)(10,1)
\pspolygon[fillstyle=hlines,hatchangle=0](10,1)(13,1)(13,3)
\pspolygon[fillstyle=hlines,hatchangle=0](10,3)(10,5)(13,5)
\rput[r](9.9,4){$\omega(c_B,\tau_M)$}
\rput[r](9.9,1){$(0,2)$}
\rput[r](9.9,5){$(0,4)$}
\rput[l](13.1,2){$\omega(c_A,\tau_M)$}
\rput[l](13.1,4){$\omega(\tau_M)$}
\rput[l](13.1,1){$(1,2)$}
\rput[l](13.1,5){$(1,4)$}
\rput[b](11.5,5.1){$\omega(\tau_M)$}
\rput[r](9.9,2){$\tilde{\omega}(c,\CT_{c_B}\circ \tau_M)$}
\rput[t](11.5,.9){${\omega}(\CT_{c_A}\circ \tau_M)$}
\rput(5,3){$[0,1] \times [a,b] \times \Sigma \times S^2$}
\end{pspicture}

\bp
Let us first treat the $p_1$ case, there is a map $G_U$ from $\partial W$ to $SU(3)$
such that $\tau_1$ is the stabilisation of $G_U \circ \tau_M$ over $\partial W$,
$$G_U=\begin{array}{ll}1 & \mbox{on}\;[0,1] \times \{4\} \times \Sigma \\
1 & \mbox{on}\;[0,1] \times \{2\} \times (\Sigma \setminus c_A \times [-1,1])\\
1 & \mbox{on}\;\{1\} \times [3,4] \times \Sigma\\
1 & \mbox{on}\;\{1\} \times [2,4] \times (\Sigma \setminus c_A \times [-1,1])\\
1 & \mbox{on}\;\{0\} \times [3,4] \times (\Sigma \setminus c_B \times [-1,1])\\
1 & \mbox{on}\;\{0\} \times [2,4] \times (\Sigma \setminus (c \cup c_B \times [-1,1]))\\
\CT_{c_A} & \mbox{on}\;[0,1] \times \{2\} \times c_A \times [-1,1]\\
\CT_{c_B} & \mbox{on}\;\{0\} \times \{3\} \times c_B \times [-1,1]\\
\CT_{c_B} & \mbox{on}\;\{0\} \times [2,3] \times (c_B \times [-1,1] \setminus
c \times [-1,1])\\
%T_{c} & \mbox{on}\;\{0\} \times [2,3] \times (c \times [-1,1] \setminus
%c_B \times [-1,1])
\end{array} $$
Furthermore, on $\{1\} \times [2,3] \times c_A \times [-1,1]$, $G_U$ factors through the natural projection onto $[2,3] \times [-1,1]$,
on $\{0\} \times [3,4] \times c_B \times [-1,1]$, $G_U$ factors through the natural projection onto $[3,4] \times [-1,1]$, and on 
$\{0\} \times [2,3] \times (c \times [-1,1] \setminus
c_B \times [-1,1])$, $G_U$ factors through the natural projection onto $[2,3] \times [-1,1]$.

$p_1(W;\tau_{1|\partial W})$ only depends on the homotopy class of $G_U$ with respect to $\{(1,4)\} \times \Sigma$. Since $p_1([0,1] \times [a,b] \times \Sigma;\tau_{0|\partial ([0,1] \times [a,b] \times \Sigma)} )$ is defined by a homotopic map. The two $p_1$ coincide.

The coincidence of the integrals can be seen on the
picture where the projection of $ST([0,1] \times [a,b] \times \Sigma) \stackrel{\tau_M}{=} 
[0,1] \times [a,b] \times \Sigma \times S^2$ on $[0,1] \times [a,b]$ is represented on the left-handside, and the projection of $ST([0,1] \times [2,4] \times \Sigma) \stackrel{\tau_M}{=} 
[0,1] \times [2,4] \times \Sigma \times S^2$ onto $[0,1] \times [2,4]$ is represented on the right-hand side. Next $\omega_1$ can be defined on the
(products by $\Sigma \times S^2$ of the) two hatched triangles by pulling-back $\omega(c_B,\tau_M)$ and $\omega(c_A,\tau_M)$ using the pictured horizontal projections. On the remaining part, $\omega_1$ can be filled in by the pull-back
of $\omega$ under the obvious diffeomorphism from this remaining part to
$[0,1] \times [a,b] \times \Sigma \times S^2$.
Now, the $z_n$ corresponding to the hatched triangles will vanish
while the $z_n$ corresponding to the remaining part 
equals $z_n({[0,1] \times ST([a,b] \times \Sigma)};\omega)$ thanks to \inconstlemtranspform.
\eop

\begin{lemma}
\label{lemdessindeux}
Let $c^{\prime}_A$ be a curve with the same class as $c_A$ in $H_1(\Sigma;\ZZ/2\ZZ)$.
Let $$H:  [0,1] \times \{2\} \times \Sigma  \longrightarrow  SO(3)$$
satisfy
$$\begin{array}{l} H([0,1] \times \{2\} \times \partial \Sigma)=\{1\}\\
H=\CT_{c_A}\;\;\mbox{ on }\;\{(0,2)\} \times \Sigma,\\ 
H=\CT_{c^{\prime}_A}\;\;\mbox{ on }\;\{(1,2)\} \times \Sigma.\end{array}$$
Let $\omega_3$ be a closed form on $ST([0,1] \times [2,3] \times \Sigma)$
such that
$$\omega_3=\left\{ \begin{array}{ll} \omega(\tau_M)  & \mbox{on}\; [0,1] \times ST([2,3] \times \Sigma)_{|\{3\} \times \Sigma} 
\\
\omega(H \circ \tau_M)  & \mbox{on}\; [0,1] \times ST([2,3] \times \Sigma)_{|\{2\} \times \Sigma}\\ 
\omega(c_A,\tau_M)& \mbox{on}\;\{0\} \times ST([2,3] \times \Sigma)\\
\omega(c^{\prime}_A,\tau_M)& \mbox{on}\;\{1\} \times ST([2,3] \times \Sigma).
\end{array} \right.$$
Let $\tau_3$ be the trivialisation of $T([0,1] \times [2,3] \times \Sigma) \otimes \CC$
corresponding to the trivialisations used to define $\omega_3$.
Then $$z_n({[0,1] \times ST([2,3] \times \Sigma)};\omega_3)=-\frac{
p_1([0,1] \times [2,3] \times \Sigma;\tau_{3|\partial ([0,1] \times [2,3] \times \Sigma)} )}{4}
\delta_n.$$
%%% SIGNE
\end{lemma}
\bp By an obvious modification of Lemma~\ref{lemdessin} above, we would not change either side of the equality
to be shown by setting rather
$$\omega_3=\left\{ \begin{array}{ll}\omega_{\cal Q}(c_A,\tau_M) & \mbox{on}\; [0,1] \times ST([2,3] \times \Sigma)_{|\{3\} \times \Sigma} 
\\
\omega_{\cal Q}(c^{\prime}_A,\tau_M)  & \mbox{on}\; [0,1] \times ST([2,3] \times \Sigma)_{|\{2\} \times \Sigma} \\
\omega(H \circ \tau_M)  & \mbox{on}\;\{0\} \times ST([2,3] \times \Sigma)\\
\omega(\tau_M) & \mbox{on}\;\{1\} \times ST([2,3] \times \Sigma).
\end{array} \right.$$
Now, both $\omega_3$ and $\tau_3$ trivially extend over $[0,1] \times [2,3] \times \partial A$ and can be glued along $[0,1] \times \{2\}$ with the implicit picture
of Lemma~\ref{lemcasa} with $c^{\prime}_A$ instead of $c_A$, and therefore
$\omega^{\prime}$ instead of $\omega$ and $\tau^{\prime}_0$ instead of $\tau_0$.
Then Lemma~\ref{lemcasa} applied to $A \cup_{\partial A} [2,3] \times \partial A$
instead of $A$ tells us that
$$4 z_n(\omega_3)+ 4 z_n([0,1] \times A;\omega^{\prime})=-p_1([0,1] \times [2,3] \times \partial A;\tau_3)\delta_n-p_1([0,1] \times A;\tau^{\prime}_0)\delta_n,$$
%%% SIGNE
and that
$$4 z_n([0,1] \times A;\omega^{\prime})=- p_1([0,1] \times A;\tau^{\prime}_0)\delta_n.$$
%%% SIGNE
\eop

This lemma has the following corollary whose proof is similar and therefore left 
to the reader.

\begin{lemma}
\label{lemmoddeux}
If there exist curves $c^{\prime}$, $c^{\prime}_A$ and $c^{\prime}_B$ that
are homologous modulo 2 to $c$, $c_A$ and $c_B$, respectively such that
Lemma~\ref{lemdifsurf} is true when replacing $(c, c_A, c_B)$ by $(c^{\prime}, c^{\prime}_A, c^{\prime}_B)$, then Lemma~\ref{lemdifsurf} is true for $(c, c_A, c_B)$.
\end{lemma}

\begin{lemma}
\label{lemnonintge}
Lemma~\ref{lemdifsurf} is true when $c$ and $c_B$ do not intersect.
\end{lemma}
\bp 
According to Lemma~\ref{lemdessin}, it is enough to prove that in this case,
$$z_n({[0,1] \times ST([2,4] \times \Sigma)};\omega_1)=-\frac14 p_1(W;\tau_{1|\partial W}) \delta_n.$$
We may assume that $$c_A \times [-1,1]= c  \times [-1,1] \coprod
c_B \times [-1,1].$$ 
Use an isotopy of $\{0\} \times [2,4] \times \Sigma$ supported away from 
$\{0\} \times [2,4] \times c  \times [-1,1]$ to lower 
lower $[3+\varepsilon,4-\varepsilon] \times c_B \times [-1,1]$ to 
$[2+\varepsilon,3-\varepsilon] \times c_B \times [-1,1]$ without changing either side of the equality.
After this isotopy, the trivialisation and the form over $\{0\} \times [2,4] \times \Sigma$ coincide with the trivialisation and the form over $\{1\} \times [2,4] \times \Sigma$. Therefore both sides of the equality to be shown vanish.
\eop

As a direct corollary of these two lemmas, we get the following lemma.
\begin{lemma}
\label{lemnonintal}
Lemma~\ref{lemdifsurf} is true when the algebraic intersection $\langle c,c_B \rangle $ of $c$ and $c_B$ in $\Sigma$ is even.
\end{lemma}
\eop

\begin{lemma}
\label{lemintpos}
There exists an integer $p^+$ and an element $I_n \in \CA_n(\emptyset)$ such
that if, under the assumptions of Lemma~\ref{lemdifsurf}, $c$ and $c_B$ are transverse and the sign of all their intersection points
is positive, then\\
$z_n({[0,1] \times ST([a,b] \times \Sigma)};\omega)= \langle c,c_B \rangle I_n$ and \\
$p_1([0,1] \times [a,b] \times \Sigma;\tau_{0|\partial ([0,1] \times [a,b] \times \Sigma)} )=\langle c,c_B \rangle p^+$.
\end{lemma}

Once this lemma is proved, using Lemma~\ref{lemnonintal} for two curves $c$ and $c_B$ as in Lemma~\ref{lemintpos} such that $\langle c,c_B \rangle =2$ shows that
$I_n=-\frac{p^+}{4}\delta_n$. Therefore, Lemma~\ref{lemdifsurf} is true for any
two curves that only intersect positively. Now, since it is easy to change $c_B$
without changing its class in $H_1(\Sigma;{\ZZ}/2\ZZ)$ so that $c$ and $c_B$
only intersect positively,  Lemma~\ref{lemdifsurf} will be proved right after
Lemma~\ref{lemintpos} is proved.

\noindent{\sc Proof of Lemma~\ref{lemintpos}:}
First isolate the intersection points of $c$ and $c_B$ inside
boxes $[-2,2] \times [-2,2]$ that $c \times [-1,1]$ intersects as $[-2,2] \times
[-1,1]$ and $c_B \times [-1,1]$ as $-[-1,1] \times
[-2,2]$.
Then lower $[3+\varepsilon,4-\varepsilon] \times c_B \times [-1,1]$ to 
$[2+\varepsilon,3-\varepsilon] \times c_B \times [-1,1]$ in $\{0\} \times [2,4] \times \Sigma$  except on the cubes $[2,4] \times [-2,2]^2$ by an isotopy of the framed link $c \cup c_B$ supported away from $[2,4] \times [-1,1]^2$.
This does not change either side of the equality and after this isotopy, we
may assume that $$\omega_1=\pi_{\{1\} \times ST([2,4] \times \Sigma)}^{\ast}\left(\omega_{1|\{1\} \times ST( [2,4] \times \Sigma)}\right)$$
except over the cubes $[0,1] \times [2,4] \times [-2,2]^2$.
But on the boundaries of the products by $S^2$ of these cubes, the value
of $\omega$ is always the same and since $H^2(\partial ([0,1] \times [2,4] \times [-2,2]^2 \times S^2))=H^2([0,1] \times [2,4] \times [-2,2]^2 \times S^2)$, $\omega$ extends as a closed form there, and we may choose the same extension for all the cubes that are the only ones to produce nonzero integrals according to \inconstlemtranspform ~and that all produce the same integral $I_n$.
Similarly, $p^+$ is the obstruction to extend the trivialisation associated to $\omega_1$
on the boundary of such a cube.
\eop
\newpage
\section{Simultaneous normalization of the fundamental forms}
\setcounter{equation}{0}
\label{secnormap}

This section is devoted to the proof of Proposition~\ref{propnormasim}.
We use real coefficients for homology and cohomology.

\subsection{Sketch}

First note that the closed 2-forms 
$\omega(M_I)$ \index{NN}{omegaMI@$\omega(M_I)$} on $\partial C_2(M_I)$ 
defined in the beginning of Section~\ref{secskp} (after Remark~\ref{rkfirstap}) are antisymmetric and extend as closed antisymmetric 2-forms on $C_2(M_I)$ because of \inconstlemexisfunap.

Now, we wish to arrange the forms $\omega_{M_I}=\omega(M_I)$ as in Proposition~\ref{propnormasim}.
To do that, we shall first show how to make $\omega_M$ explicit in some part of 
$C_2(M)$.

Recall that $[-4,4] \times \partial A^i$ denotes a regular neighborhood of $\partial A^i$
embedded in $M$, that intersects $A^i$ as $ [-4,0] \times \partial A^i$.
All the neighborhoods $[-4,4] \times \partial A^i$ are disjoint from each other.
Throughout this paragraph, we shall use the corresponding  coordinates on the image of this implicit embedding.

For $t \in [-4,4]$, set \index{NN}{Ait@$A^i_t$}
$$A^i_t=\left\{\begin{array}{ll} A^i \cup  ([0,t] \times \partial A^i) & \mbox{if}\;\; t \geq 0\\
 A^i \setminus ( ]t,0] \times \partial A^i) & \mbox{if} \;\;t \leq 0\end{array} \right.$$
$$ \partial A^i_t =\{t\} \times \partial A^i.$$
Choose curves $(b^i_j)_{j=1,\dots,g_i}$,
and $(y^i_j)_{j=1,\dots,g_i}$ of $\partial A^i$ such that
\begin{itemize}
\item the homology classes of the $(b^i_j)_{j=1,\dots,g_i}$ form a basis of
$\CL(M \setminus \mbox{Int}(A^i))$,
\item and $\langle y_j^i, [b_k^i]\rangle_{\partial A^i}= \delta_{jk}$ (thus, the homology classes of the $(y^i_j)$ form a basis of 
$H_1(M \setminus A^i)$).
\end{itemize}
Choose a basepoint $p^i$ in $\partial A^i$
outside the neighborhoods $a^i_j \times [-1,1]$ of the $a^i_j$ and outside neighborhoods $b^i_j \times [-1,1]$ the $b^i_j$.
Fix a path $[p^i,(0,0,1)]$ from $p^i$ to $(0,0,1)$ in 
$$B_M(1) \setminus \left(\mbox{Int}(A^i)\cup_{k,k\neq i} A^k_4\right),$$ that is extended into a path $[p^i,\infty(v)]$ \index{NN}{piinfty@$[p^i,\infty(v)]$} by the vertical line $(0,0) \times [1,\infty[ \cup \{\infty(v)\}$, where $\infty(v)$ is the intersection with $\partial C_1(M)$ of the closure in $C_1(M)$ of the vertical half-line $(0,0) \times [1,\infty[$. Fix a closed two-form $\omega(p^i)$ on $(M \setminus \mbox{Int}(A^i))$ such that
\begin{itemize}
\item the integral of $\omega(p^i)$ \index{NN}{omegapi@$\omega(p^i)$} along a closed surface of 
$(M \setminus \mbox{Int}(A^i))$ is its algebraic intersection with $[p^i,\infty(v)]$,
\item the support of $\omega(p^i)$ intersects $(B_M \setminus \mbox{Int}(A^i)) $
inside a tubular neighborhood of $[p^i,\infty(v)]$ disjoint from $$\left( \cup_{k,k\neq i} A^k_4\right) \cup \left([-4,4] \times \left( \cup_{j=1}^{g_i} ( (a^i_j \times [-1,1]) \cup (b^i_j \times [-1,1]) )\right) \right).$$
\item $\omega(p^i)$ restricts as the usual volume form on $\partial C_1(M)=S^2$.
\end{itemize}

Here, a {\em two-chain\/} is a linear rational combination of smooth compact oriented surfaces with boundaries. The {\em integral\/} of a 2-form along such a chain is the corresponding linear rational combination of its integral along the surfaces. The {\em support\/} of such a 2-chain is the union of the involved surfaces.
A {\em 2-cycle\/} is a two-chain with empty (or null) boundary.

For any $i \in N$, and for any $j \in \{1,2, \dots, g_i\}$,
extend $\eta(a^i_j)$ \index{NN}{etaaij@$\eta(a^i_j)$} on $[-4,4] \times \partial A^i$ into a closed one-form $\eta(a^i_j)$ supported on $[-4,4] \times a^i_j \times [-1,1]$ where $\eta(a^i_j)$ is again given by the formula.
$$\eta(a^i_j)=\pi_{[-1,1]}^{\ast}(\eta_{[-1,1]}).$$
Let $S(a^i_j)$ be a 2-chain in $A^i_4$ with boundary $4 \times a^i_j \times 0$ that intersects 
$[-4,4] \times \partial A^i$ along $[-4,4] \times a^i_j \times 0$. 
%We assume that $S(a^i_j)$ does not meet the support of the $\omega(p^i)$.

Let $i \in N$, and let $j \in \{1,2, \dots, g_i\}$.
Choose a $2$-chain
$S(b^i_j)$ \index{NN}{Sbij@$S(b^i_j)$} in $(B_M \setminus \mbox{Int}(A^i))$ that is bounded by $b^i_j$ whose support is disjoint from all the supports of the $\omega(p^k)$, and that intersects 
$A^k_4$ as a combination of $S(a^k_{\ell})$ for all $k \in N$ such that $k \neq i$.
Define a   
closed one-form \index{NN}{etabij@$\eta(b^i_j)$}
$\eta(b^i_j)$ on $(M \setminus \mbox{Int}(A^i))$
such that 
\begin{itemize}
\item the integral of $\eta(b^i_j)$ along a closed curve of 
$(M \setminus \mbox{Int}(A^i))$ is its algebraic intersection with $S(b^i_j)$,
\item the support of $\eta(b^i_j)$ is in a neighborhood of $S(b^i_j)$ disjoint from the support of the 
$\omega(p^k)$, for all $k \neq i$.
\item for all $k \neq i$, the restriction of $\eta(b^i_j)$ to $A^k_4$, $k\neq i$ is a rational combination of the $\eta(a^k_{\ell})$, $\ell=1,\dots,g_k$ (the coefficients are linking numbers determined by the first condition).
\end{itemize}

In Subsection~\ref{proofnorma}, we shall prove that these forms can be used as follows to make $\omega_M$ explicit in some parts of $C_2(M)$.

\begin{proposition}
\label{propnorma}
With the above notations, we can choose $\omega_M$ \index{NN}{omegaM@$\omega_M$} so that:
\begin{enumerate} 
\item for every $i \in N$, the restriction of $\omega_M$ to 
$$A^i \times (C_1(M) \setminus A^i_3) \subset C_2(M)$$
equals
$$\sum_{(j,k) \in \{1, \dots, g_i\}^2} \ell(z_j^{i}, (4 \times y_k^{i})) p_1^{\ast}(\eta(a^i_j)) \wedge p_2^{\ast}(\eta(b^i_k)) + p_2^{\ast}(\omega(p^i))$$
where $p_1$ and $p_2$ denote the first and the second projection of $C_2(M)$ onto $C_1(M)$,
respectively;
\item for every $i$, for any $j \in \{1,2, \dots, g_i\}$, $$\int_{S(a^i_j) \times p^i} \omega_M=0;$$
\item $\omega_M$ is fundamental with respect to $\tau_M$.
\end{enumerate}
\end{proposition}

Assume that Proposition~\ref{propnorma} is proved. This is the goal of Subsection~\ref{proofnorma}.
When changing some $A^i$ into some $B^i$ with the same Lagrangian, it is easy to change the restrictions of $\omega_M$ inside the parts mentioned in the first paragraph of the statement of Proposition~\ref{propnorma}
(and inside their symmetric parts under $\iota$ that are also determined by the statement). Indeed, all the forms $\eta(a^i_j)$, $\eta(b^i_j)$ and $\omega(p^i)$ can be defined on the parts of the $M_I$ where they are needed so that these forms coincide with each other whenever it makes sense, and so that they have the properties that were required for $M$.
(Recall that the $\eta(a^i_j)$ are defined both in $A^i$ and $B^i$ and that they are identical near $\partial A^i$ and $\partial B^i$ while $\omega(p^i)$ is supported in $\left(M \setminus \left(\cup_{k \in N} \mbox{Int}(A_k)\right)\right)$ and  while  the $\eta(b^i_j)$ restrict to the $A^k$ as a (fixed by the clover data) combination of $\eta(a^k_{\ell}$).
Define $\omega_0(M_I)$ on \\ 
$ D(\omega_0(M_I))=$ \index{NN}{Domegazero@$D(\omega_0(M_I))$}
$$\left(C_2(M_I) \setminus \left( \cup_{i \in I} p_{12}^{-1}\left((B^i_{-1} \times B^i_3) \cup (B^i_3 \times B^i_{-1}) \right) \right)  \right) \cup p_{12}^{-1}(\mbox{diag}(M \setminus \{\infty\})) $$
so that
\begin{enumerate}
\item $\omega_0(M_I)=\omega_M$ on $C_2\left(M \setminus (\cup_{i \in I}B^i_{-1})\right)$,
\item $$\omega_0(M_I)=\sum_{(j,k) \in \{1, \dots g_i\}^2} \ell(z_j^{i}, (4 \times y_k^{i})) p_1^{\ast}(\eta(a^i_j)) \wedge p_2^{\ast}(\eta(b^i_k)) + p_2^{\ast}(\omega(p^i))$$ on $p_{12}^{-1}(B^i \times (M \setminus B^i_3) )$ when $i \in I$,
\item $\omega_0(M_I)=-\iota^{\ast}(\omega_0(M_I))$ on $p_{12}^{-1}((M \setminus B^i_3) \times B^i )$  when $i \in I$,

\item $\omega_0(M_I)=
\omega(M_I)$ on $\partial C_2(M_I)$. (See the definition after Remark~\ref{rkfirstap}.)
\end{enumerate}
Note that this definition is consistent. 
%Also note that $D(\omega_0(M_i=M_{\{i\}}))$
%deformation retracts onto $(C_2(M_i) \setminus C_2(B^i_3)) \cup ST(B^i_4)$.

\begin{lemma}
\label{lemker}
With the above notation, for any $i \in N$, 
$\omega_0(M_i=M_{\{i\}})$ vanishes on the kernel
 of the map induced by the inclusion
$$H_2\left(D(\omega_0(M_i)\right) \longrightarrow H_2(C_2(M_i)).$$
\end{lemma}

This lemma is surprisingly difficult to prove for me. It will be proved in Subsection~\ref{subpropclef}. Assume it for the moment. Then (the cohomology class of) $\omega_0(M_i)$  is in the image of the natural map
$$H^2(C_2(M_i)) \longrightarrow H^2(D(\omega_0(M_i))).$$
Therefore  $\omega_0(M_i)$ extends to a closed form $\omega_1(M_i)$, and
\index{NN}{omegaMi@$\omega(M_i)$} $$\omega(M_i)=\frac{\omega_1(M_i) - \iota^{\ast}(\omega_1(M_i))}{2}$$
is an admissible form for $M_i$ .
%with respect to $\tau(M_i)$.

Now, for any $I \subset N$, we may define
\index{NN}{omegaMI@$\omega(M_I)$}
$$\omega(M_I)= \left\{ \begin{array}{ll} \omega_{0}(M_I) \; &\mbox{on}\; C_2(M_I) \setminus \left( \cup_{i \in I} p_{12}^{-1}\left((B^i_{-1} \times B^i_4) \cup (B^i_4 \times B^i_{-1}) \right) \right)\\
\omega(M_i) &\mbox{on}\;C_2(B^i_4) \; \mbox{for}\; i \in I
\end{array}\right.$$ 
since the $C_2(B^i_4)$ do not intersect.
These forms $\omega(M_I)$ satisfy the conclusions of
Proposition~\ref{propnormasim} that will be proved once Proposition~\ref{propnorma}
and Lemma~\ref{lemker} are proved. Their proofs will occupy the next two subsections.

\subsection{Proof of Proposition~\ref{propnorma}}
\label{proofnorma}
The homology classes of the $(z_j^i \times (4 \times y_k^i))_{(j,k) \in \{1,\dots,g_i\}^2}$ and $(p^i \times \partial C_1(M))$ form a basis of
$$H_2\left(A^i \times (C_1(M) \setminus A^i_3)\right)=(H_1(A^i) \otimes H_1(M \setminus A^i)) \oplus H_2(C_1(M) \setminus A^i).$$
The evaluation of $L_M$ (defined after \inconstlemhctwo)
along these classes is $\ell(z_j^i, (4 \times y_k^i))$ for the first ones
and $1$ for the last one. 
In particular the form of the statement integrates correctly on this basis.

Let us first prove Proposition~\ref{propnorma} when $N=\{1\}$. Set $A^1=A$, and forget about the 
superfluous superscripts $1$. Let $\omega_0$ be a 2-form fundamental with respect to $\tau_M$ given by \inconstlemexisfunap, and let $\omega$ be the closed 2-form defined on 
$\left(A_1 \times (C_1(M) \setminus \mbox{Int}(A_2))\right)$ by the statement (naturally extended). Since this form $\omega$ integrates correctly on $H_2\left(A_1 \times (C_1(M) \setminus \mbox{Int}(A_2))\right)$, there exists a one-form $\eta$ on
$\left(A_1 \times (C_1(M) \setminus \mbox{Int}(A_2))\right)$ such that
$\omega=\omega_0 +d \eta$.

This form $\eta$ is closed on $A_1 \times \partial C_1(M)$.
Since $H^1\left(A_1 \times (C_1(M) \setminus \mbox{Int}(A_2))\right)$ maps surjectively to 
$H^1(A_1 \times \partial C_1(M))$, we may extend $\eta$ to a closed one-form $\tilde{\eta}$
on $\left(A_1 \times (C_1(M) \setminus \mbox{Int}(A_2))\right)$. Changing $\eta$ into 
$(\eta-\tilde{\eta})$, turns $\eta$ into a primitive of $(\omega-\omega_0)$
that vanishes on $A_1 \times \partial C_1(M)$.

Let $\chi$ be a smooth function on $C_2(M)$ supported in 
$\left(A_1 \times (C_1(M) \setminus \mbox{Int}(A_2))\right)$, and constant with the value 
1 on $(A \times C_1(M) \setminus A_3)$.

Set $$\omega_a=\omega_0 +d \chi \eta.$$

Then $\omega_a$ is a closed form that has the required form on $(A \times (C_1(M) \setminus A_3))$. Furthermore, the restrictions of 
$\omega_a$ and $\omega_0$ agree on $\partial C_2(M)$ since $d \chi \eta$
vanishes there (because $\eta$ vanishes on $A_1 \times \partial C_1(M)$).

Adding to $\eta$ a combination $\eta_c$ of the closed forms $p_2^{\ast}(\eta(b_j))$
that vanish on $A_1 \times \partial C_1(M)$ does not change the above properties,
but adds 
$$\int_{p \times ( [2,3] \times a_j)} d (\chi \eta_c)=\int_{p \times (3 \times a_j)} \eta_c$$
to
$\int_{p \times S(a_j)} \omega_a$.
Therefore since the $p_2^{\ast}(\eta(b_j))$ generate the dual of $\CL_A$, we may choose $\eta_c$ so that all the $\int_{p \times S(a_j)} \omega_a$ vanish.
After this step, $\omega_a$ is a closed form that takes the prescribed values on 
$$PS_a=\partial C_2(M) \cup \left(A \times (C_1(M) \setminus A_3)\right)$$
and such that all the $\int_{p \times S(a_j)} \omega_a$ vanish.
In order to make $\omega_a$ antisymmetric with respect to $\iota^{\ast}$, 
we apply 
similar modifications to $\omega_a$ on the symmetric part $(C_1(M) \setminus \mbox{Int}(A_2)) \times A_1$. The support of these modifications is disjoint
from the support of the previous ones. Thus, they do not interfer and transform 
$\omega_a$ into a closed form $\omega_b$ with the additional properties:
\begin{itemize}
\item $\omega_b$ has the prescribed form on $(C_1(M) \setminus A_3) \times A$, (It is prescribed there because of the prescribed form on $\left( A \times (C_1(M) \setminus A_3) \right)$ and because of the prescribed antisymmetry with respect to $\iota^{\ast}$.)
\item $\int_{S(a_j) \times p } \omega_b =0$, for all $j=1, \dots g_1$.
\end{itemize}

Now, the form $\omega_M=\frac{\omega_b - \iota^{\ast}(\omega_b)}{2}$
has all the required properties, and the proposition is proved for $N=\{1\}$.

We now proceed by induction on $\sharp N=i$.
We start with a $2$-form $\omega_0$ that satisfies all the hypotheses on $\{1,\dots,i-1\}$ instead of $N=\{1,\dots,i\}$, and by the first
step, we also assume that we have a $2$-form $\omega$ that satisfies all the hypotheses on $\{i\}$ instead of $N=\{1,\dots,i\}$, with $A_1^i$ replacing $A^i$. 

Now, we proceed similarly.
There exists a one-form $\eta$ on
$C_2(M)$ such that
$\omega=\omega_0 +d \eta$.
The exact sequence
$$ 0=H^1(C_2(M)) \longrightarrow H^1(\partial C_2(M)) \longrightarrow H^2(C_2(M), \partial C_2(M)) \cong H_4(C_2(M))=0 $$
shows that $H^1(\partial C_2(M))$ is trivial. Therefore, $\eta$ is exact
on $\partial C_2(M)$, we can assume that $\eta$ vanishes on  $\partial C_2(M)$, and we do assume so.\\ Let $\chi$ be a smooth function on $C_2(M)$ supported in 
$\left(A^i_1 \times (C_1(M) \setminus \mbox{Int}(A^i_2))\right)$, and constant with the value 
1 on $\left(A^i \times (C_1(M) \setminus A^i_3)\right)$.
Again, we are going to modify $\eta$ by some closed forms
so that $$\omega_a=\omega_0 +d \chi \eta$$
has the prescribed value on 
$$PS_a=\partial C_2(M) \cup \left(\cup_{k \in N} \left(A^k \times (C_1(M) \setminus A^k_3)\right) \right) \cup \left(\cup_{k \in N \setminus \{i\}}\left( (C_1(M) \setminus A^k_3) \times A^k\right)\right).$$

Our form $\omega_a$ is as required anywhere except possibly in
$$\left(A^i_1 \times (C_1(M) \setminus \mbox{Int}(A^i_2))\right) \setminus
\left(A^i \times (C_1(M) \setminus A^i_3)\right)$$
and in particular in the intersection of this domain with the domains 
where it was previously normalized, that are included
in $$\left(A^i_1 \times (\partial C_1(M) \cup (\cup_{k=1}^{i-1}A^k))\right).$$

Recall that $\eta$ vanishes on $A^i_1 \times \partial C_1(M)$.
Our assumptions also imply that $\eta$ is closed on $A^i_1 \times A^k$.
Let us prove that they imply that $\eta$ is exact on $A^i_1 \times A^k$ for any $k <i$.
To do that it suffices to check that:
\begin{enumerate}
\item For any $j =1,\dots,g_i,$
$\int_{b^i_j \times p^k} \eta=0.$
\item For any $j =1,\dots,g_k,$
$\int_{p^i \times b^k_j} \eta=0.$
\end{enumerate}
Let us prove the first assertion.
Since $\int_{b^i_j \times \infty(v)} \eta=0$,
$$\int_{b^i_j \times p^k} \eta= \int_{\partial (b^i_j \times [p^k,\infty(v)])}
\eta=\int_{ b^i_j \times [p^k,\infty(v)]}
(\omega-\omega_0).$$
where $\int_{ b^i_j \times [p^k,\infty(v)]}
\omega=0$ because the supports of the $\eta(b^i_{\ell})$ do not intersect $[p^k,\infty(v)]$. 
Now, $$\int_{ b^i_j \times [p^k,\infty(v)]}
\omega_0= - \int_{ S(b^i_j) \times \partial [p^k,\infty(v)]} \omega_0
=\int_{ S(b^i_j) \times \{p^k\}} \omega_0.$$

The latter integral vanishes because
\begin{enumerate}
\item $S(b^i_j)$ intersects $A^k_4$ as copies of $S(a^k_{\ell})$,
\item $\int_{ S(a^k_{\ell}) \times p^k}\omega_0=0$ (that is the second condition of Proposition~\ref{propnorma}), and,
\item the integral of $\omega_0$ also vanishes on the remaining part of $S(b^i_j) \times p^k$ because $\omega_0$ is determined on $((C_1(M) \setminus A^k_4) \times A^k)$ and because the support of $\omega(p^k)$ is disjoint from $S(b^i_j)$.
\end{enumerate}

Let us prove the second assertion.
Again, since $\eta$ vanishes on $\partial C_2(M)$, $\int_{\infty(v) \times b^k_j} \eta=0$ and therefore
$$\int_{p^i \times b^k_j} \eta=- \int_{ [p^i,\infty(v)] \times b^k_j }
(\omega-\omega_0).$$
$\int_{ [p^i,\infty(v)] \times b^k_j }
\omega_0=0$ because of the form of $\omega_0$ on $(C_1(M) \setminus A^k_4) \times A^k$.
$$\int_{ [p^i,\infty(v)] \times b^k_j }
\omega= \int_{ \partial [p^i,\infty(v)] \times S(b^k_j) }
\omega= - \int_{ \{p^i\} \times S(b^k_j) } \omega.$$
Again, we know that this integral is zero along the intersection 
of $\{p^i\} \times S(b^k_j) $ with $A^i \times (C_1(M) \setminus A^i_4)$ because 
 $S(b^k_j)$ does not meet the support of $\omega(p^i)$, and we conclude
because $\int_{ \{p^i\} \times S(a^i_{\ell}) } \omega=0$ and because $S(b^k_j)$ intersects $A^i_4$ along copies of $ S(a^i_{\ell})$.

Since $\eta$ is exact on the annoying parts, we can assume that it identically vanishes there.

Thus, $\omega_a$ takes the prescribed values on $A^i \times (C_1(M) \setminus A^i_4)$, $\omega_a$ coincides with $\omega_0$ where $\omega_0$ was prescribed and $\omega_a$ integrates
correctly along the $S(a^k_{\ell}) \times p^k$ and their symmetric with respect to $\iota$, for $k\neq i$.
Let us now modify $\eta$ by adding a linear combination of $p_2^{\ast}(\eta(b^i_j))$ that vanishes on the $A^i_1 \times A^k$, for $k<i$, and thus without changing the above properties so that 
the integrals of $\omega_a$ along the $\{p^i\} \times S(a^i_{\ell})$ vanish, for $\ell=1, \dots, g_i$, too.
Let $f: H_1(M \setminus \mbox{Int}(A^i)) \longrightarrow \RR$ be the linear map 
defined by
$$f(a^i_{\ell})=-\int_{\{p^i\} \times S(a^i_{\ell})}\omega_a.$$
There exists a combination $\eta_c$ of $p_2^{\ast}(\eta(b^i_j))$
such that for any $x \in \CL_A$, $f(x)=\int_{p^i \times x}\eta_c$.

Observe that 
$$\int_{\{p^i\} \times S(b^k_j)} \omega_a=\int_{\{\infty(v)\} \times S(b^k_j)} \omega_a - \int_{[p^i,\infty(v)] \times b^k_j} \omega_a=0.$$
This implies that $f(\mbox{Im}(H_1(A^k) \longrightarrow H_1(M \setminus \mbox{Int}(A^i))))=0$. Thus, $\eta_c$ vanishes on $A^i_1 \times A_k$.
Changing $\eta$ into $(\eta+\eta_c)$ does not change $\omega_a$ on the prescribed set but adds $\int_{\{p^i\} \times (4 \times  a^i_{\ell})}\eta_c=f(a^i_{\ell})$
to $\int_{\{p^i\} \times  S(a^i_{\ell})} \omega_a$ that becomes 0. 

After this step, $\omega_a$ is a closed form that takes the prescribed values on $PS_a$ such that the integrals of $\omega_a$ along the $(\{p^i\} \times S(a^i_{\ell}))$ vanish, for $\ell=1, \dots, g_i$.
In order to make $\omega_a$ antisymmetric with respect to $\iota^{\ast}$, 
we apply 
similar modifications to $\omega_a$ on the symmetric part $(C_1(M) \setminus \mbox{Int}(A^i_2)) \times A^i_1$. Again, the support of these modifications is disjoint
from the support of the previous ones. Thus, they do not interfer and transform 
$\omega_a$ into a closed form $\omega_b$ with the additional properties:
\begin{itemize}
\item $\omega_b$ has the prescribed form on $(C_1(M) \setminus A^i_3) \times A^i$,
\item $\int_{S(a^i_j) \times p^i } \omega_b =0$, for all $j=1, \dots g_i$.
\end{itemize}

Now, the form $\omega_M=\frac{\omega_b - \iota^{\ast}(\omega_b)}{2}$
has all the required properties, and Proposition~\ref{propnorma} is proved.

\eop

\subsection{Proof of Lemma~\ref{lemker}}
\label{subpropclef}

We need some more notation before stating the key proposition that will lead to  the proof of  Lemma~\ref{lemker}.

We assume that our trivialisation $\tau_M$ (fixed since the beginning of Section~\ref{secskp}) maps the unit tangent vector of $( \{s\} \times x \times [0,1])$ at $(s,x,t)$ to the first basis vector $e_1 \in \RR^3$ for any $(s,x,t) \in [-4,0] \times a^i_j \times [0,1]$.
When $X$ is a unit vector field on the image of some chain $F(P)$ of $M \setminus \infty$, then
$\mbox{diag}(X)(F(P))$ denotes the chain of the blow-up of the diagonal that is the image of $P$ under the map
$\left( p \mapsto (F(p),X(F(p)) )\right)$.

Fix a basepoint $p(a^i_j)$ \index{NN}{paij@$p(a^i_j)$} on every curve $a^i_j$.

We shall construct several $2$-chains. 
%All the manifolds and chains are oriented. 
%The boundary of an oriented manifold is oriented with the outward normal first 
%convention. 
When chains are presented as products, the orientation is the product orientation with respect to the order of the factors from left to right. Similarly, when chains are described with coordinates varying inside oriented manifolds, they are oriented as the image of the product of the oriented manifolds ordered by the order used to write the coordinates. The field $\RR$ is given its standard orientation. 
%$a$ is oriented.
A minus sign reverses the orientation.

Let $a$ denote a curve $a^i_j=0 \times a^i_j \times 0$ that lives inside
the fixed neighborhood $([-4,4] \times (a \times [-1,1] \subset \partial A^i))$ in $M$. Fix $I \subset N$.
We are about to explicitly construct a 2-dimensional cycle $F(a)$ in $C_2(M_I)$ of the form \index{NN}{Fa@$F(a)$}\index{NN}{Ca@$C(a)$}
$$F(a)=C(a) \cup e(S_0(a))(S^2_{\mbox{\small diag}}=ST(M)_{|p(a)})$$
$$-\left(S_0(a) \times (4 \times p(a))\right)  
\cup -\left((4 \times p(a)) \times S_0(a) \right) \cup 
\mbox{diag}({n})(S_0(a))$$
for a rational two-chain $S_0(a)$ in $C^i_I$ ($=A^i$ or $B^i$) whose boundary
is $(0 \times a \times 0)$, equipped with a vector field ${n}$ and with a rational number $e(S_0(a))$,
and for a two-chain $C(a)$ in $C_2([0,4] \times \partial A^i)$ described below.

Let us first describe $S_0(a)$, ${n}$ and $e(S_0(a))$.

There exists a minimal positive nonzero integer $k$ such that $ka=0$ in $H_1(C^i_I;\ZZ)$.
Then there exists a connected surface $S_{2,-4}$ \index{NN}{Stwomf@$S_{2,-4}$} embedded in $C^i_{I,-4}$ whose
boundary is $$ \partial S_{2,-4} = (\{-4\} \times a \times \{0,\frac2{2k-1}, \frac4{2k-1}, \dots, \frac{2k-2}{2k-1}\})$$ 
and whose normal is $\tau_M^{-1}(.,e_1)$ on $ \partial S_{2,-4}$.
It can be easily extended to an immersed connected surface $S_2$\index{NN}{Stwo@$S_2$} with boundary $k$
copies of $\{0\} \times a \times \{0\}$
such that
$S_2 \cap \left([-2,0] \times \partial A\right)$ is made of $k$ copies of
$[-2,0] \times a  \times \{0\}$, and such that
$S_2$ is obtained from $S_{2,-4}$ by gluing $k$ embedded annuli transverse to the vector field $\tau_M^{-1}(.,e_1)$.
Then $S_0(a)=\frac1kS_2$,\index{NN}{Szeroa@$S_0(a)$} ${n}$ is the normal vector to $S_2$, (that is homotopic to $\tau_M^{-1}(.,e_1)$ on $S_2 \setminus S_{2,-4}$,) and \index{NN}{eSzeroa@$e(S_0(a))$}
$$ e(S_0(a))=\frac{g(S_2)+k-1}{k}=\frac12 -\frac{\chi(S_2)}{2k}.$$

Let us now describe the wanted two-chain $C(a)$ \index{NN}{Ca@$C(a)$}in $C_2([0,4] \times \partial A^i)$ with boundary 
$$\partial C(a)=\left(\{0\} \times a \times \{0\} \times (4 \times p(a))\right) $$
$$\cup \left((4 \times p(a)) \times \{0\} \times a \times \{0\}) \right) \cup 
-\mbox{diag}(e_1)(\{0\} \times a \times \{0\}).$$

For two given based parametrized closed curves $\{x(v);v \in [0,1]/(0 \sim 1)\}$ and $\{y(v);v \in [0,1]/(0 \sim 1)\}$ with respective basepoints $x(0)$ and $y(0)$, we fix a cobordism $T(x,y)$  \index{NN}{Txy@$T(x,y)$} in the torus $x \times y$ between the diagonal $\{(x(v),y(v)); v \in [0,1]\}$ and $(x(1)\times y) \cup (x \times y(0))$. (The notation $(x(1)\times y)$ stands for $(\{x(1)\}\times y)$ for lightness).
The fixed cobordism $T(x,y)$ is
the image of the triangle $\{(v,w) \in (0,1)^2; v \geq w\}$ by the map $((v,w) \mapsto (x(v),y(w))$ in the torus $x \times y$.
%$a \times \{(t,0)\} \times a \times \{(u,0)\} \subset C_2(M)$.
$$ \begin{pspicture}[.2](0,0)(2.5,2.5)
\psset{xunit=.5cm,yunit=.5cm}
\rput[tr](3.8,.8){$v$}
\rput[tr](.8,3.8){$ w$}
\rput(3,2){$T(x,y)$}
%\rput[br](3.8,4.2){$a \times \{(t,0)\} \times p(a) \times \{0\}$}
%\rput[tl](4.2,3.8){$ p(a) \times \{0\} \times a \times \{(u,0)\}$}
\pspolygon*[linecolor=lightgray](1,1)(4,1)(4,4)
\psline{->}(1,1)(4,1)
\psline{->}(1,1)(1,4)
\psline{->}(1,4)(4,4)
\psline{->}(4,1)(4,4)
\psline{->}(4,4)(1,1)
\end{pspicture}$$

Let $t$ and $u$ be such that $0 \leq t < u <1$. Let $s \in [-4,4]$. 
Define \index{NN}{Atus@$A(t,u;s)$}
$$A(t,u;s)=\overline{\{((s,x,t),(s,x,t+\lambda(u-t))) ;\lambda \in ]0,1], x \in a\} } \subset C_2(M)$$
$$=\{((s,x,t),(s,x,t+\lambda(u-t))) ; \lambda \in ]0,1],x \in a\} \cup \mbox{diag}(e_1)(s \times a \times t).$$

Let $C(a)$ \index{NN}{Ca@$C(a)$} denote the sum of the following $2$-chains in $C_2(M)$:

\begin{enumerate}
\item $T( (0 \times a \times 0), (0 \times a \times 1))$
\item $A(0,1;0)$
\item $ (0 \times a \times 0) \times  \left[-([0,4] \times p(a) \times 1) \cup 
(4 \times p(a) \times [0,1]) \right]$
\item $\left((4 \times p(a) \times 0) \times (0 \times a \times [0,1])\right)
\cup \left(([0,4] \times p(a) \times 0) \times (0 \times a \times 1)\right)$
%\item $ ([0,4] \times a \times t) \times  (4 \times p(a) \times u)$
%\item $  (4 \times p(a) \times t) \times ([0,4] \times a \times u)$
%\item $ \pm $
%\item $-\{(x,t,\lambda,x,u,\lambda) ;x \in a, \lambda \in [0,4]\}$.\\
%(Here, the minus means that the implicit orientation induced by $a \times %[0,4]$ is reversed.)\\
%\item $A(t,u;4)=-\overline{\{((4,x,t),(4,x,t+\lambda(u-t)) ;x \in a, \lambda %\in ]0,1]\}} \subset C_2(M)$
\end{enumerate}

Now,  $F(a)$ \index{NN}{Fa@$F(a)$} is completely defined as a 2-dimensional cycle in $C_2(M_I)$, for all $I \subset N$.

We postpone the proof of the most difficult lemma to the end of this subsection.

\begin{lemma}
\label{lemkey}
With the above definition, $F(a)=F(a^i_j)$ \index{NN}{Fa@$F(a)$} is null-homologous in $C_2(M)$.
\end{lemma}

Assuming this lemma, the proof of Lemma~\ref{lemker} goes as follows.

We first prove:

\begin{lemma}
\label{lemcalker}
For any $i \in N$, the homology classes of the $F(a^i_j)$ for $j=1, \dots g_i$ generate 
the kernel of
$$H_2(D(\omega_0(M_i)) \longrightarrow H_2(C_2(M_i)).$$
\end{lemma}
\bp
Since the inclusion from $D(\omega_0(M_i)$ to $(C_2(M_i) \setminus C_2(B^i_{-1})) \cup ST(B^i)$ is a homotopy equivalence, it is the same to prove that the $F(a^i_j)$ generate 
the kernel of
$$H_2((C_2(M_i) \setminus C_2(B^i_{-1})) \cup ST(B^i)) \longrightarrow H_2(C_2(M_i)).$$
It is also the same to prove that the $F(a^i_j)$ generate the kernel of
$$H_2((C_2(M) \setminus C_2(A^i_{-1})) \cup ST(A^i)) \longrightarrow H_2(C_2(M)).$$
Set $C=C_1(M)$ and $A=A^i$. Then 
$C_2(M) \setminus C_2(A)$ has the homotopy type of 
$$ \breve{C}_2(C) \setminus \breve{C}_2(A) \stackrel{\mbox{\small def}}{=} \left( C^2 \setminus A^2 \right) \setminus \mbox{diag}(C \setminus A) .$$
Let us compute the real homology of $C_2(M) \setminus C_2(A)$ in degrees 1 and 2. Recall that $C$ has the homology of a point, that $H_2(C \setminus A)=\RR[\partial C]$, and that $H_1(C \setminus A)$ is isomorphic to $\CL_A$. Therefore, in the Mayer-Vietoris sequence
$$ H_{i}((C \setminus A)^2) \hfl{\alpha_{i}} H_i(C \times (C \setminus A)) \oplus H_i((C \setminus A) \times C)$$ $$ \longrightarrow H_i \left((C \times (C \setminus A)) \cup_{(C \setminus A)^2} ((C \setminus A) \times C) \right) \hfl{\partial_i} H_{i-1}((C \setminus A)^2), $$ the map $\alpha_{i}$ is onto for $i \geq 1$, and
the map $\partial_i$ is an injection into the kernel of $\alpha_{i-1}$ that is an injection, when $i=1$ and $2$.
This shows that 
$$H_i \left(C^2 \setminus A^2=(C \times (C \setminus A)) \cup ((C \setminus A) \times C) \right)=\{0\}\;\mbox{ for} \; i \in \{1,2\}.$$

Let us now compute the effect of removing $\mbox{diag}(C \setminus A)$, by using the long exact sequence:
$$ \longrightarrow H_{i+1}(C^2 \setminus A^2,\breve{C}_2(C) \setminus \breve{C}_2(A)) \longrightarrow  H_i(\breve{C}_2(C) \setminus \breve{C}_2(A)) \longrightarrow H_i(C^2 \setminus A^2)  \longrightarrow... $$
%H_{i}(C^2 \setminus A^2,\breve{C}_2(C) \setminus \breve{C}_2(A)).$$
By excision, $H_{i+1}(C^2 \setminus A^2,\breve{C}_2(C) \setminus \breve{C}_2(A))$ is isomorphic to $$H_{i+1}\left((C \setminus A) \times \RR^3, (C \setminus A) \times (\RR^3\setminus \{0\})\right)\cong H_{i-2}(C \setminus A)
\otimes H_2(S^2).$$
This shows that $H_1(C_2(M) \setminus C_2(A))=H_1(\breve{C}_2(C) \setminus \breve{C}_2(A))$ is trivial, and that $H_2(C_2(M) \setminus C_2(A))$ is generated by the homology class of a fiber of the unit tangent bundle of $(C \setminus A)$. Since this fiber is not null homologous
in $C_2(M)$, we conclude that $H_2(C_2(M) \setminus C_2(A))=\RR[S^2_{\mbox{\small diag}}].$

We end the computation of $H_2((C_2(M) \setminus C_2(A)) \cup ST(A))$ by gluing
$(ST(C)\cong C \times S^2)$ to $(C_2(M) \setminus C_2(A))$ along $(ST(C \setminus A)\cong (C \setminus A) \times S^2)$, and by using the Mayer-Vietoris sequence that yields the exact sequence
$$0 \longrightarrow \RR[S^2_{\mbox{\small diag}}] \longrightarrow H_2((C_2(M) \setminus C_2(A_{-1})) \cup ST(A)) $$ $$ \hfl{\partial_{MV}} H_1((C \setminus A_{-1}) \times S^2)
\hfl{\gamma_1} H_1(C \times S^2).$$
The kernel of $\gamma_1$ is freely generated by the curves 
$\mbox{diag}(e_1)(a_j)$ for $j=1, \dots, g_i$, and the assumptions on the $F(a^i_j)$ 
%in the statement of Proposition~\ref{propclef} 
ensure that $\partial_{MV}$ maps $F(a^i_j)$
to $\pm \mbox{diag}(e_1)(a_j)$. Therefore, 
$$H_2((C_2(M) \setminus C_2(A_{-1})) \cup ST(A)) = \oplus_{j=1}^{g_i}[F(a^i_j)] \oplus \RR[S^2_{\mbox{\small diag}}].$$
and Lemma~\ref{lemkey} ensure that the kernel of $$H_2(C_2(M) \setminus C_2(A_{-1}) \cup ST(A)) \longrightarrow H_2(C_2(M))$$ is generated by the $[F(a^i_j)]$.
\eop

Thus, in order to prove Lemma~\ref{lemker}, it is enough to prove that, for any $j=1,\dots,g_i$, 
$$\int_{F(a^i_j)}\omega_0(M_i)=0.$$
Set $a=a^i_j$.
Since $C(a)$ only depends on $\CL_{A^i}$ and since $C(a)$ lives inside the part $C_2\left([0,4]\times \partial A^i\right)$ of $C_2(M)$ or $C_2(M_i)$ where $\omega_M$ and $\omega_0(M_i)$ coincide, 
$$\int_{C(a)}\omega_M=\int_{C(a)}\omega_0(M_i).$$
The normalizations also imply that the integrals of the forms vanish on 
$$-\left(S_0(a) \times (4 \times p(a))\right)  
\cup -\left((4 \times p(a)) \times S_0(a) \right).$$
Now, since $F(a)$ is null-homologous in $C_2(M)$, $\int_{F(a)}\omega_M=0$, and we are left with the proof that
$$\int_{\mbox{ \small diag}({n})(S_0(a))}\omega_M + e(S_0(a))= \int_{\mbox{ \small diag}({n})(S_{0,i}(a))}\omega_0(M_i) + e(S_{0,i}(a)).$$ 
$$\int_{\mbox{ \small diag}({n})(S_0(a))}\omega_M=\frac1k \int_{\mbox{ \small diag}({n})(S_{2})} \omega_M
=\frac1{2k}\chi(TS_2; \tau_M^{-1}(.,e_2)_{|\partial S_2})$$
according to Lemma~\ref{lemintsur}. 
Then  Lemma~\ref{lemcalchi} implies that
$$\int_{\mbox{ \small diag}({n})(S_0(a))}\omega_M=\frac{kd + \chi(S_2)}{2k}$$
where $d$ is the degree of the map $\tau_M^{-1}(.,e_2)$ from $a$ to $S(\RR \vec{N}(a) \oplus \RR \vec{T}(a))$. Thus,
%When $\tau_M^{-1}(.,e_2)$ extends to a disk $D_a$ 
%that frames $a$ like $\partial A$,
%$$\int_{\mbox{ \small diag}({n})(S_0(a))}\omega_M = %\frac{2-2g(S_2)}{2k}=1-e(S_0(a)),$$ and in general 
%$$\int_{\mbox{ \small diag}({n})(S_0(a))}\omega_M + e(S_0(a))=1 + %\frac1{2}\chi(D_a; \tau_M^{-1}(.,e_2)_{|a}).$$
$$\int_{\mbox{ \small diag}({n})(S_0(a))}\omega_M + e(S_0(a))=\frac{d+1}2.$$
Since similar equalities hold when $(M_i,\omega_0(M_i), S_{0,i}(a))$ replaces $(M,\omega_M,S_0(a))$ according to Lemma~\ref{lemintsurbis},
%$$\int_{\mbox{ \small diag}({n})(S_{0,i}(a))}\omega_0(M_i) + e(S_{0,i}(a))
%=1 + \frac1{2}\chi(D_a; \tau_M^{-1}(.,e_2)_{|a}),$$
$$\int_{\mbox{ \small diag}({n})(S_{0,i}(a))}\omega_0(M_i) + e(S_{0,i}(a))
=\frac{d+1}2,$$
and this concludes the proof of Lemma~\ref{lemker}
%We are now left with the (surprisingly difficult ???) proof of 
up to Lemma~\ref{lemkey}. \eop

\noindent{\sc Proof of Lemma~\ref{lemkey}}
We shall first replace $F(a)$ by an integral cycle.

Let $S=S_{2,-4} \cup ([-4,0] \times a  \times \{0,\frac2{2k-1}, \frac4{2k-1}, \dots, \frac{2k-2}{2k-1}\})$ be the connected surface embedded in $A$ whose
boundary is $$ \partial S = \{0\} \times a \times \{0,\frac2{2k-1}, \frac4{2k-1}, \dots, \frac{2k-2}{2k-1}\},$$ such that
$$S \cap \left([-4,0] \times \partial A\right) =[-4,0] \times a  \times \{0,\frac2{2k-1}, \frac4{2k-1}, \dots, \frac{2k-2}{2k-1}\}.$$
Choose a tubular neighborhood $S  \times [-\frac1{2k},\frac1{2k-1}]$ of $S$
such that, when $v=(t,(x,\frac{2(\ell -1)}{2k-1}))$ belongs to the part $([-4,0] \times a  \times \{0,\frac2{2k-1}, \frac4{2k-1}, \dots, \frac{2k-2}{2k-1}\})$
of $S$, the element $(v,u)$ of this neighborhood reads
$(t,x,\frac{2(\ell -1)}{2k-1} + u)$ in $[-4,0] \times (a \times [-1,1])$.

We use $S$ to find an integral cycle that is homologous to $k F(a)$.

Define two basepoints on $(a \times [0,1] \subset \partial A \subset M)$, one left one and one right one: 
$$p_{\ell}=(0,p(a),0) \;\; \mbox{and} \;\; p_r=(0,p(a),1)$$

\begin{sublemma}
\label{sublemG}
$k F(a)$ is homologous to the integral cycle $G(a)$, \index{NN}{Ga@$G(a)$} where
$G(a)$ is the sum of the seven following 
2-chains.

\begin{enumerate}
\item $\cup_{\ell=1}^{k}T((0 \times a \times \{\frac{2(\ell-1)}{2k-1}\}) \times (0\times a \times \{\frac{2\ell-1}{2k-1}\}))$
\item $A(\partial S)=\cup_{\ell=1}^{k}A(\frac{2(\ell-1)}{2k-1},\frac{2\ell-1}{2k-1};0)$
\item $A_r(0)=-\cup_{\ell=1}^{k} (0 \times  p(a) \times [0,\frac{2(\ell-1)}{2k-1}]) \times (0 \times a \times \{\frac{2\ell-1}{2k-1}\}) $
\item $A_{\ell}(0)=-\cup_{\ell=1}^{k} (0 \times a \times \{\frac{2(\ell-1)}{2k-1}\} ) \times (0 \times  p(a) \times 
[\frac{2\ell-1}{2k-1},1])$
\item $-S \times   p_r$
\item $-p_{\ell}  \times (S \times \{\frac1{2k-1}\})$
\item $\mbox{diag}({n})(S) + (g(S)+k-1)S^2_{\mbox{diag}}$
\end{enumerate}
\end{sublemma}
\bp
The cycle
$\left(k\mbox{diag}({n})S_0(a)-\left(\mbox{diag}({n})(S) \cup_{\ell=1}^{k} \mbox{diag}({e_1})(a \times [0,\frac{2(\ell-1)}{2k-1}])\right)\right)$
is homologous to 
$$\mbox{diag}({e_1})\left( (S_2 \setminus S_{2,-4})
-\left(\cup_{\ell=1}^{k} ([-4,0] \times a  \times \{\frac{2(\ell-1)}{2k-1}\}\cup a \times [0,\frac{2(\ell-1)}{2k-1}])\right)\right)$$
and is therefore null-homologous because $H_2(A;\QQ)=0$.

Therefore, in the definition of  $k F(a)$,
we may change  
$k\mbox{diag}(n)(S_0(a))$
into
$$\mbox{diag}(n)(S) \cup \cup_{\ell=1}^{k} \mbox{diag}(e_1)\left(a \times [0,\frac{2(\ell-1)}{2k-1}]\right)$$ and stay in the same homology class.

We may also change
$$-k\left(S_0(a) \times (4 \times p(a))\right) \cup $$
$$\cup
k(0 \times a \times 0) \times  \left(-([0,4] \times p(a) \times 1) \cup 
(4 \times p(a) \times [0,1]) \right)$$
into

$$-\left(S \cup_{\ell=1}^{k} (a \times [0,\frac{2(\ell-1)}{2k-1}])\right) \times p_r. $$
Indeed these chains have the same boundary and they live inside
$$\left(A \setminus (0 \times a \times 1) \right) \times \left( (-[0,4] \times p(a) \times 1) \cup 
(4 \times p(a) \times [0,1]) \right)$$
that has the homotopy type of $A$, and that therefore has a trivial $H_2$.

Similarly, we may change 
$$-k\left((4 \times p(a)) \times S_0(a) \right) \cup $$
$$\cup k\left((4 \times p(a) \times 0) \times (0 \times a \times [0,1])\right)
\cup k\left(([0,4] \times p(a) \times 0) \times (0 \times a \times 1)\right)$$
into

$$-p_{\ell} \times \left((S \times \frac1{2k-1}) \cup( 0 \times -a \times (\cup_{\ell=1}^{k} ([\frac{2\ell-1}{2k-1},1]))\right)$$
inside $$\left( ([0,4] \times p(a) \times [0,1]) \times A\right) \setminus (\mbox{diag}(0 \times p(a) \times [0,1]))$$
that has the homotopy type of $A$.

Afterwards, it is enough to check that, for any $\ell =1, \dots k$,
$$T( (0 \times a \times 0) \times  (0 \times a \times 1)) \cup A(0,1;0)$$
may be replaced by the chain $C(a)(\frac{2(\ell-1)}{2k-1},\frac{2\ell-1}{2k-1})$
with boundary
$$ p_{\ell} \times  (0 \times a \times 1) \cup (0 \times a \times 0) \times p_r
\cup -\mbox{diag}(e_1)(0 \times a \times 0)$$
where, for any $(t,u) \in [0,1]^2$ such that $0 \leq t <u \leq 1$, $C(a)(t,u)$ is the sum of the following chains:
\begin{enumerate}
\item
$T((0 \times a \times \{t\}) , (0\times a \times \{u\}))$
\item $A(t,u;0)$
\item $-(0 \times  p(a) \times [0,t]) \times (0 \times a \times \{u\}) $
\item$-(0 \times a \times \{t\} ) \times (0 \times  p(a) \times 
[u,1])$
\item $-\mbox{diag}(e_1)\left(a \times [0,t]\right)$
\item $(a \times [0,t]) \times p_r$
\item $-p_{\ell} \times ( 0 \times a \times [u,1])$.
\end{enumerate}

Since $C(a)(0,1)$ is homologous to $T( (0 \times a \times 0) \times  (0 \times a \times 1)) \cup A(0,1;0)$ and since the $C(a)(t,u)$ form a continuous family
of chains with the same boundary indexed by a connected set, we are done.
\eop

Then to prove Lemma~\ref{lemkey}, it is sufficient (and necessary) to prove that $G(a)$ represents $0$ in $(H_2(C_2(M);\RR)=\RR[S^2_{\mbox{\tiny diag}}])$.
To do that, we shall describe some homotopies explicitly, and we need some more notation.
All our homotopies will take place inside $C_2\left((S \times [0,\frac{1}{2k-1}]) \cup (\{0\} \times a \times [0,1])\right)$. We shall 
%write $ [-4,0] \times a \times [0,1] $ rather than $a \times [0,1] \times %[-4,0]$ and  
use the implicit coordinates there.

Let $g(S)$ be the genus of $S$. We define a Morse function $h_S$ from $S$ to $[-6 -3g(S),0]$ that is the height function of $S$ with respect to an embedding such as in the following picture
of $S$, and such that
$h_S$ coincides with the projection on $[-4,0]$ on $S \cap ([-4,0] \times \partial A)$.
In particular, $h_S$ is maximal and constant on the boundary of $S$, $h_S$ has a unique minimum at the height $(-6 -3g(S))$, $h_S$ has $(2g(S)+(k-1))$ index one critical points, and $h_S$ has no other critical points. Furthermore, $h_S^{-1}(-6-3j)$ is a circle for $j=0, \dots, g(S)-1$ and there are two critical points between $h_S^{-1}(-6-3j)$ and $h_S^{-1}(-6-3(j+1))$ for $j=0, \dots, g(S)-1$.

We construct a connected (compact) graph $\Gamma$ on $S$ that intersects every connected component of every height level of $S$
exactly once, and that intersects $S \cap ([-4,0] \times \partial A)$ as $\cup_{i=0}^{k-1} \left([-4,0] \times p(a) \times \{\frac{2i}{2k-1}\} \right)$. 

\begin{center}
\begin{pspicture}[.4](0,-.2)(5,6) 
\psset{xunit=.5cm,yunit=.5cm}   
\psecurve[linewidth=2pt]{->}(0,11)(.9,11.3)(1.8,11)(.9,10.7)(0,11)(.9,11.3)(1.8,11)
\rput[br](.9,11.4){$a \times \{0\}$}
\psecurve[linewidth=2pt]{->}(3,11)(3.9,11.3)(4.8,11)(3.9,10.7)(3,11)(3.9,11.3)(4.8,11)
\rput[b](3.9,11.4){$a \times \{\frac{2}{2k-1}\}$}
\psecurve[linewidth=2pt]{->}(6,11)(6.9,11.3)(7.8,11)(6.9,10.7)(6,11)(6.9,11.3)(7.8,11)
\rput[bl](6.9,11.4){$a \times \{\frac{2k-2}{2k-1}\}$}
\pscurve[linewidth=2pt](1.8,11)(2.4,7.4)(3,11)
\pscurve[linewidth=2pt]{-}(4.8,11)(5.4,6.8)(6,11)
\psccurve[linewidth=2pt](4.5,5)(4.8,4.25)(4.5,3.5)(4.2,4.25)
\psccurve[linewidth=2pt](4.5,2.5)(4.8,1.75)(4.5,1)(4.2,1.75)
\pscurve[linewidth=2pt](0,11)(.6,7.4)(2.8,5)(3,4.25)(3,3)(3,1.75)(3.2,1)(4.5,0)(5.8,1)(6,1.75)(6,3)(6,4.25)
(6.1,5)(7.2,6.8)(7.8,11)
%Les huit
\psecurve[linewidth=.5pt](3.4,5.2)(2.8,5)(3.4,4.8)(4.5,5)(4.5,5)
\psecurve[linestyle=dashed,dash=3pt 2pt](3.4,4.8)(2.8,5)(3.4,5.2)(4.5,5)(4.5,5)
\psecurve(5.5,5.2)(6.1,5)(5.5,4.8)(4.5,5)(4.5,5)
\psecurve[linestyle=dashed,dash=3pt 2pt](5.5,4.8)(6.1,5)(5.5,5.2)(4.5,5)(4.5,5)
\psecurve(3.8,1.15)(3.2,1)(3.8,.85)(4.5,1)(4.5,1)
\psecurve[linestyle=dashed,dash=3pt 2pt](3.8,.85)(3.2,1)(3.8,1.15)(4.5,1)(4.5,1)
\psecurve(5.2,1.15)(5.8,1)(5.2,.85)(4.5,1)(4.5,1)
\psecurve[linestyle=dashed,dash=3pt 2pt](5.2,.85)(5.8,1)(5.2,1.15)(4.5,1)(4.5,1)
% Les autres courbes meridiennes
\psecurve(4.5,3.3)(3,3)(4.5,2.7)(6,3)(4.5,3.3)
\psecurve[linestyle=dashed,dash=3pt 2pt](4.5,2.7)(3,3)(4.5,3.3)(6,3)(4.5,2.7)
\psecurve(3.6,4.4)(3,4.25)(3.6,4.1)(4.2,4.25)(3.6,4.4)
\psecurve[linestyle=dashed,dash=3pt 2pt](3.6,4.1)(3,4.25)(3.6,4.4)(4.2,4.25)(3.6,4.1)
\psecurve(5.4,1.9)(6,1.75)(5.4,1.6)(4.8,1.75)(5.4,1.9)
\psecurve[linestyle=dashed,dash=3pt 2pt](5.4,1.6)(6,1.75)(5.4,1.9)(4.8,1.75)(5.4,1.6)
\rput(4,6.1){$S$}
\end{pspicture}
\begin{pspicture}[.4](-1,-.2)(5,6) 
\psset{xunit=.5cm,yunit=.5cm} 
\psecurve{->}(0,11)(.9,11.3)(1.8,11)(.9,10.7)(0,11)(.9,11.3)(1.8,11)
%\rput[br](.9,11.4){$a \times \{0\}$}
\psecurve{->}(3,11)(3.9,11.3)(4.8,11)(3.9,10.7)(3,11)(3.9,11.3)(4.8,11)
%\rput[b](3.9,11.4){$a \times \{\frac{2}{2k-1}\}$}
\psecurve{->}(6,11)(6.9,11.3)(7.8,11)(6.9,10.7)(6,11)(6.9,11.3)(7.8,11)
\pscurve(1.8,11)(2.4,7.4)(3,11)
\pscurve{-}(4.8,11)(5.4,6.8)(6,11)
\psccurve(4.5,5)(4.8,4.25)(4.5,3.5)(4.2,4.25)
\psccurve(4.5,2.5)(4.8,1.75)(4.5,1)(4.2,1.75)
\pscurve(0,11)(.6,7.4)(2.8,5)(3,4.25)(3,3)(3,1.75)(3.2,1)(4.5,0)(5.8,1)(6,1.75)(6,3)(6,4.25)
(6.1,5)(7.2,6.8)(7.8,11)
%Gamma
\psecurve[linewidth=2pt]{*-*}(.9,10.7)(.9,10.7)(2.4,7.4)(3.9,10.7)
\psecurve[linewidth=2pt]{*-}(3.9,10.7)(3.9,10.7)(2.4,7.4)(.9,10.7)
\pscurve[linewidth=2pt]{-*}(2.4,7.4)(5.4,6.8)(6,7.1)(6.9,10.7)
\psline[linewidth=2pt]{*-*}(5.4,6.8)(4.5,5)
\psline[linewidth=2pt]{*-*}(4.5,3.5)(4.5,2.5)
\psline[linewidth=2pt]{*-*}(4.5,1)(4.5,0)
\psccurve[linewidth=2pt](4.5,5)(5.1,4.25)(4.5,3.5)(3.9,4.25)
\psccurve[linewidth=2pt](4.5,2.5)(5.1,1.75)(4.5,1)(3.9,1.75)
\rput(4.5,6.1){$\Gamma$}
\end{pspicture}

\end{center}

We also extract two connected subgraphs $\Gamma_{\ell}$ and  $\Gamma_r$ of $\Gamma$  that intersect every height level of $S$
exactly once, and such that:
\begin{itemize}
\item $\Gamma_{\ell} \cap h_S^{-1}([-6 -3g(S), -6])=\Gamma_r  \cap h_S^{-1}([-6 -3g(S), -6])$,
\item $\Gamma_{\ell} \cap ([-4,0] \times \partial A) = [-4,0] \times p(a) \times \{0\} $, and,
\item $\Gamma_{r} \cap ([-4,0] \times \partial A) =[-4,0] \times  p(a) \times \{\frac{2k-2}{2k-1}\}$.
\end{itemize}

\begin{center}
\begin{pspicture}[.4](0,-.2)(5,6) 
\psset{xunit=.5cm,yunit=.5cm} 
\psecurve{->}(0,11)(.9,11.3)(1.8,11)(.9,10.7)(0,11)(.9,11.3)(1.8,11)
%\rput[br](.9,11.4){$a \times \{0\}$}
\psecurve{->}(3,11)(3.9,11.3)(4.8,11)(3.9,10.7)(3,11)(3.9,11.3)(4.8,11)
%\rput[b](3.9,11.4){$a \times \{\frac{2}{2k-1}\}$}
\psecurve{->}(6,11)(6.9,11.3)(7.8,11)(6.9,10.7)(6,11)(6.9,11.3)(7.8,11)
\pscurve(1.8,11)(2.4,7.4)(3,11)
\pscurve{-}(4.8,11)(5.4,6.8)(6,11)
\psccurve(4.5,5)(4.8,4.25)(4.5,3.5)(4.2,4.25)
\psccurve(4.5,2.5)(4.8,1.75)(4.5,1)(4.2,1.75)
\pscurve(0,11)(.6,7.4)(2.8,5)(3,4.25)(3,3)(3,1.75)(3.2,1)(4.5,0)(5.8,1)(6,1.75)(6,3)(6,4.25)
(6.1,5)(7.2,6.8)(7.8,11)
%Gamma
\psecurve[linewidth=2pt]{*-*}(.9,10.7)(.9,10.7)(2.4,7.4)(3.9,10.7)
%\psecurve[linewidth=2pt]{*-}(3.9,10.7)(3.9,10.7)(2.4,7.4)(.9,10.7)
\psecurve[linewidth=2pt]{-*}(2.4,7.4)(2.4,7.4)(5.4,6.8)(6.9,10.7)
\psline[linewidth=2pt]{*-*}(5.4,6.8)(4.5,5)
\psline[linewidth=2pt]{*-*}(4.5,3.5)(4.5,2.5)
\psline[linewidth=2pt]{*-*}(4.5,1)(4.5,0)
\psecurve[linewidth=2pt](5.1,4.25)(4.5,3.5)(3.9,4.25)(4.5,5)(5.1,4.25)
\psecurve[linewidth=2pt](5.1,1.75)(4.5,1)(3.9,1.75)(4.5,2.5)(5.1,1.75)
\rput(4.5,6.1){$\Gamma_{\ell}$}
\end{pspicture}
\begin{pspicture}[.4](-2,-.2)(5,6.8) 
\psset{xunit=.5cm,yunit=.5cm} 
\psecurve{->}(0,11)(.9,11.3)(1.8,11)(.9,10.7)(0,11)(.9,11.3)(1.8,11)
%\rput[br](.9,11.4){$a \times \{0\}$}
\psecurve{->}(3,11)(3.9,11.3)(4.8,11)(3.9,10.7)(3,11)(3.9,11.3)(4.8,11)
%\rput[b](3.9,11.4){$a \times \{\frac{2}{2k-1}\}$}
\psecurve{->}(6,11)(6.9,11.3)(7.8,11)(6.9,10.7)(6,11)(6.9,11.3)(7.8,11)
\pscurve(1.8,11)(2.4,7.4)(3,11)
\pscurve{-}(4.8,11)(5.4,6.8)(6,11)
\psccurve(4.5,5)(4.8,4.25)(4.5,3.5)(4.2,4.25)
\psccurve(4.5,2.5)(4.8,1.75)(4.5,1)(4.2,1.75)
\pscurve(0,11)(.6,7.4)(2.8,5)(3,4.25)(3,3)(3,1.75)(3.2,1)(4.5,0)(5.8,1)(6,1.75)(6,3)(6,4.25)
(6.1,5)(7.2,6.8)(7.8,11)
\psecurve[border=1pt](4.8,11.1)(3.9,10.7)(3,10.3)(2.7,11)(3.9,11.6)(7.5,11.6)(6.9,10.7)(6,10.3)
\psecurve{->}(1.8,11.1)(.9,10.7)(0,10.3)(-.3,11)(1.2,11.9)(4.5,11.6)
\psecurve[border=1pt](-.3,11)(1.2,11.9)(4.5,11.6)(3.9,10.7)(3,10.3)
\rput[br](1.1,12){$\{0\} \times p(a) \times [0,\frac{2k-2}{2k-1}]$}
%Gamma
%\psecurve[linewidth=2pt]{*-*}(.9,10.7)(.9,10.7)(2.4,7.4)(3.9,10.7)
%\psecurve[linewidth=2pt]{*-}(3.9,10.7)(3.9,10.7)(2.4,7.4)(.9,10.7)
\psecurve[linewidth=2pt]{-*}(2.4,7.4)(5.4,6.8)(6,7.1)(6.9,10.7)(6.9,10.7)
\psline[linewidth=2pt]{*-*}(5.4,6.8)(4.5,5)
\psline[linewidth=2pt]{*-*}(4.5,3.5)(4.5,2.5)
\psline[linewidth=2pt]{*-*}(4.5,1)(4.5,0)
\psecurve[linewidth=2pt](5.1,4.25)(4.5,3.5)(3.9,4.25)(4.5,5)(5.1,4.25)
\psecurve[linewidth=2pt](5.1,1.75)(4.5,1)(3.9,1.75)(4.5,2.5)(5.1,1.75)
\rput(4.5,6.1){$\Gamma_{r}$}
\psecurve(3.9,11.3)(4.8,11)(3.9,10.7)(3,11)(3.9,11.3)
\end{pspicture}
\end{center}

Then, define three associated projections $p$, $p_{\ell}$ and $p_r$ from $S$ to $\Gamma$ such that
for any element $s$ of $S$, 
\begin{itemize}
\item $p(s)$ is the intersection point of $\Gamma$ and the connected component of $s$ in $h_S^{-1}(h_S(s))$,
\item $\{p_{\ell}(s)\}=\Gamma_{\ell} \cap h_S^{-1}(h_S(s))$,
\item $\{p_{r}(s)\}=\Gamma_{r} \cap h_S^{-1}(h_S(s))$.
\end{itemize}
The maps $p_{\ell}$ and $p_r$ factor through $h_S$, and the quotient maps will still be denoted by $p_{\ell}$ and $p_r$.

Our Morse function $h_S$ is such that $\Gamma_{\ell}$ contains all the critical points of $h_S$, and such that the shortest path in $\Gamma$ from $p_{\ell}(-6)$
to $0 \times p(a) \times \{\frac{2i}{2k-1}\}$ 
\begin{itemize}
\item is $p_{\ell}[-6,0]$, if $i=0$,
\item is $p_{r}[-6,0]$, if $i=k-1$,
\item contains exactly $k-i$ critical points for $0<i<k$.
\end{itemize}

For any $p \in \Gamma \setminus \Gamma_{\ell}$, define the path $[p,p_{\ell}(h_S(p))]$ as the injective path from  $p$ to $p_{\ell}(h_S(p))$
in $$\left(\Gamma \cap h_S^{-1}( [h_S(p),0]) \right) \cup (\{0\} \times p(a) \times [0,\frac{2k-2}{2k-1}]) $$
%that is over the height of $p$ in $\Gamma$ and 
that is in $\Gamma_{\ell}$ on its way down (that is when $h_S$ is decreasing). 
%and whose image intersects $(\{0\} \times p(a) \times [0,\frac{2k-2}{2k-1}])$ %as little as possible. 
Define $[p,p_{\ell}(h_S(p))]$ as the constant path for
$p \in \Gamma_{\ell}$.

For any $p \in \Gamma \setminus \Gamma_{r}$,
define the path $[p,p_{r}(h_S(p))]$ as the injective path from  $p$ to $p_r(h_S(p))$
in $$\left(\Gamma \cap h_S^{-1}( [h_S(p),0]) \right) \cup (\{0\} \times p(a) \times [0,\frac{2k-2}{2k-1}]) $$ 
%that is over the height of $p$ in $\Gamma$ and 
whose image intersects $(\{0\} \times p(a) \times [0,\frac{2k-2}{2k-1}])$ as little as possible, and define $[p,p_{r}(h_S(p))]$ as the constant path for
$p \in \Gamma_{r}$.

$$\begin{pspicture}[.4](0,-.2)(6,7) 
\psset{xunit=.5cm,yunit=.5cm} 
\psecurve{->}(0,11)(.9,11.3)(1.8,11)(.9,10.7)(0,11)(.9,11.3)(1.8,11)
%\rput[br](.9,11.4){$a \times \{0\}$}
\psecurve{->}(3,11)(3.9,11.3)(4.8,11)(3.9,10.7)(3,11)(3.9,11.3)(4.8,11)
%\rput[b](3.9,11.4){$a \times \{\frac{2}{2k-1}\}$}
\psecurve{->}(6,11)(6.9,11.3)(7.8,11)(6.9,10.7)(6,11)(6.9,11.3)(7.8,11)
\pscurve(1.8,11)(2.4,7.4)(3,11)
\pscurve{-}(4.8,11)(5.4,6.8)(6,11)
\psccurve(4.5,5)(4.8,4.25)(4.5,3.5)(4.2,4.25)
\psccurve(4.5,2.5)(4.8,1.75)(4.5,1)(4.2,1.75)
%la courbe qui fait tout le tour
\pscurve(0,11)(.6,7.4)(2.8,5)(3,4.25)(3,3)(3,1.75)(3.2,1)(4.5,0)(5.8,1)(6,1.75)(6,3)(6,4.25)
(6.1,5)(7.2,6.8)(7.8,11)
\psecurve[linewidth=2pt, border=1pt](4.8,11.1)(3.9,10.7)(3,10.3)(2.7,11)(3.9,11.6)(7.5,11.6)(6.9,10.7)(6,10.3)
\psecurve[linewidth=2pt]{-}(1.8,11.1)(.9,10.7)(0,10.3)(-.3,11)(1.2,11.9)(4.5,11.6)
\psecurve[linewidth=2pt, border=1pt](-.3,11)(1.2,11.9)(4.5,11.6)(3.9,10.7)(3,10.3)
%%la grasse a gauche qui descend
\psecurve[linewidth=2pt]{*-*}(.9,10.7)(.9,10.7)(2.4,7.4)(3.9,10.7)
\psecurve[linestyle=dashed,dash=3pt 2pt]{*-}(3.9,10.7)(3.9,10.7)(2.4,7.4)(.9,10.7)
\pscurve[linestyle=dashed,dash=3pt 2pt]{-*}(2.4,7.4)(5.4,6.8)(6,7.1)(6.9,10.7)
%ligne qui descend du point de connexion le plus bas
\psline[linestyle=dashed,dash=3pt 2pt]{*-*}(5.4,6.8)(4.5,5)
\psline[linestyle=dashed,dash=3pt 2pt]{*-*}(4.5,3.5)(4.5,2.5)
\psline[linestyle=dashed,dash=3pt 2pt]{*-*}(4.5,1)(4.5,0)
\psccurve[linestyle=dashed,dash=3pt 2pt](4.5,5)(5.1,4.25)(4.5,3.5)(3.9,4.25)
\psccurve[linestyle=dashed,dash=3pt 2pt](4.5,2.5)(5.1,1.75)(4.5,1)(3.9,1.75)
%la grasse a droite qui monte
\psecurve[linewidth=2pt](6.9,10.7)(6.9,10.7)(6,7.1)(5.4,6.8)
\psecurve[linewidth=2pt]{->}(2.4,7.4)(2.4,7.4)(3.5,7.1)(5.4,6.8)
\rput[t](3.3,7){$[p_1,p_{\ell}(p_1)]$}
\rput(8,7.1){$p_1$}
\psecurve[linewidth=2pt]{->}(4.5,3.5)(5.1,4.25)(4.5,5)(3.9,4.25)(4.5,3.5)
\rput[l](5.2,4.25){$p_2$}
\psframe*[linecolor=white](0,3.5)(3.8,5)
\rput(1.9,4.25){$[p_2,p_{\ell}(p_2)]$}
%\rput[r](2.9,4.25){$[p_2,p_{\ell}(p_2)]$}
%\rput(4.5,6.1){$\Gamma$}
\psecurve(3.9,11.3)(4.8,11)(3.9,10.7)(3,11)(3.9,11.3)
\end{pspicture}
\begin{pspicture}[.4](-1,-.2)(5,7) 
\psset{xunit=.5cm,yunit=.5cm} 
\psecurve{->}(0,11)(.9,11.3)(1.8,11)(.9,10.7)(0,11)(.9,11.3)(1.8,11)
%\rput[br](.9,11.4){$a \times \{0\}$}
\psecurve{->}(3,11)(3.9,11.3)(4.8,11)(3.9,10.7)(3,11)(3.9,11.3)(4.8,11)
%\rput[b](3.9,11.4){$a \times \{\frac{2}{2k-1}\}$}
\psecurve{->}(6,11)(6.9,11.3)(7.8,11)(6.9,10.7)(6,11)(6.9,11.3)(7.8,11)
\pscurve(1.8,11)(2.4,7.4)(3,11)
\pscurve{-}(4.8,11)(5.4,6.8)(6,11)
\psccurve(4.5,5)(4.8,4.25)(4.5,3.5)(4.2,4.25)
\psccurve(4.5,2.5)(4.8,1.75)(4.5,1)(4.2,1.75)
%la courbe qui fait tout le tour
\pscurve(0,11)(.6,7.4)(2.8,5)(3,4.25)(3,3)(3,1.75)(3.2,1)(4.5,0)(5.8,1)(6,1.75)(6,3)(6,4.25)
(6.1,5)(7.2,6.8)(7.8,11)
\psecurve[linewidth=2pt, border=1pt](4.8,11.1)(3.9,10.7)(3,10.3)(2.7,11)(3.9,11.6)(7.5,11.6)(6.9,10.7)(6,10.3)
\psecurve[linewidth=2pt]{-}(1.8,11.1)(.9,10.7)(0,10.3)(-.3,11)(1.2,11.9)(4.5,11.6)
\psecurve[linewidth=2pt, border=1pt](-.3,11)(1.2,11.9)(4.5,11.6)(3.9,10.7)(3,10.3)
\psecurve[linestyle=dashed,dash=3pt 2pt]{*-*}(.9,10.7)(.9,10.7)(2.4,7.4)(3.9,10.7)
%la grasse a gauche qui monte
\psecurve[linewidth=2pt]{*->}(.9,10.7)(.9,10.7)(1.2,9.1)(2.4,7.4)
\psecurve[linestyle=dashed,dash=3pt 2pt]{*-}(3.9,10.7)(3.9,10.7)(2.4,7.4)(.9,10.7)
\pscurve[linestyle=dashed,dash=3pt 2pt]{-*}(2.4,7.4)(5.4,6.8)(6,7.1)(6.9,10.7)
%ligne qui descend du point de connexion le plus bas
\psline[linestyle=dashed,dash=3pt 2pt]{*-*}(5.4,6.8)(4.5,5)
\psline[linestyle=dashed,dash=3pt 2pt]{*-*}(4.5,3.5)(4.5,2.5)
\psline[linestyle=dashed,dash=3pt 2pt]{*-*}(4.5,1)(4.5,0)
\psccurve[linestyle=dashed,dash=3pt 2pt](4.5,5)(5.1,4.25)(4.5,3.5)(3.9,4.25)
\psccurve[linestyle=dashed,dash=3pt 2pt](4.5,2.5)(5.1,1.75)(4.5,1)(3.9,1.75)
%la grasse a droite qui descend
\psecurve[linewidth=2pt]{-}(6.9,10.7)(6.9,10.7)(6.6,9.1)(5.4,6.8)
%\psecurve[linewidth=2pt](2.4,7.4)(2.4,7.4)(3.9,7.1)(5.4,6.8)
\psframe*[linecolor=white](-.8,8.35)(.6,9.85)
\rput[r](.8,9.1){$[p_3,p_{\ell}(p_3)]$}
\rput(8.5,9.1){$p_3$}
\psecurve[linewidth=2pt]{->}(4.5,3.5)(5.1,4.25)(4.5,5)(3.9,4.25)(4.5,3.5)
\rput[l](5.2,4.25){$p_2$}
\psframe*[linecolor=white](0,3.5)(3.8,5)
\rput(1.9,4.25){$[p_2,p_{\ell}(p_2)]$}
%\rput(4.5,6.1){$\Gamma$}
\psecurve(3.9,11.3)(4.8,11)(3.9,10.7)(3,11)(3.9,11.3)
\end{pspicture}$$

$$\begin{pspicture}[.4](0,-.2)(6,7) 
\psset{xunit=.5cm,yunit=.5cm} 
\psecurve{->}(0,11)(.9,11.3)(1.8,11)(.9,10.7)(0,11)(.9,11.3)(1.8,11)
%\rput[br](.9,11.4){$a \times \{0\}$}
\psecurve{->}(3,11)(3.9,11.3)(4.8,11)(3.9,10.7)(3,11)(3.9,11.3)(4.8,11)
%\rput[b](3.9,11.4){$a \times \{\frac{2}{2k-1}\}$}
\psecurve{->}(6,11)(6.9,11.3)(7.8,11)(6.9,10.7)(6,11)(6.9,11.3)(7.8,11)
\pscurve(1.8,11)(2.4,7.4)(3,11)
\pscurve{-}(4.8,11)(5.4,6.8)(6,11)
\psccurve(4.5,5)(4.8,4.25)(4.5,3.5)(4.2,4.25)
\psccurve(4.5,2.5)(4.8,1.75)(4.5,1)(4.2,1.75)
%la courbe qui fait tout le tour
\pscurve(0,11)(.6,7.4)(2.8,5)(3,4.25)(3,3)(3,1.75)(3.2,1)(4.5,0)(5.8,1)(6,1.75)(6,3)(6,4.25)
(6.1,5)(7.2,6.8)(7.8,11)
\psecurve[linewidth=2pt, border=1pt](4.8,11.1)(3.9,10.7)(3,10.3)(2.7,11)(3.9,11.6)(7.5,11.6)(6.9,10.7)(6,10.3)
\psecurve[linestyle=dashed,dash=3pt 2pt]{->}(1.8,11.1)(.9,10.7)(0,10.3)(-.3,11)(1.2,11.9)(4.5,11.6)
\psecurve[linestyle=dashed,dash=3pt 2pt, border=1pt](-.3,11)(1.2,11.9)(4.5,11.6)(3.9,10.7)(3,10.3)
\psecurve[linestyle=dashed,dash=3pt 2pt]{*-*}(.9,10.7)(.9,10.7)(2.4,7.4)(3.9,10.7)
%la grasse a gauche qui monte
\psecurve[linewidth=2pt]{*-}(3.9,10.7)(3.9,10.7)(2.4,7.4)(.9,10.7)
\pscurve[linestyle=dashed,dash=3pt 2pt]{-*}(2.4,7.4)(5.4,6.8)(6,7.1)(6.9,10.7)
%ligne qui descend du point de connexion le plus bas
\psline[linestyle=dashed,dash=3pt 2pt]{*-*}(5.4,6.8)(4.5,5)
\psline[linestyle=dashed,dash=3pt 2pt]{*-*}(4.5,3.5)(4.5,2.5)
\psline[linestyle=dashed,dash=3pt 2pt]{*-*}(4.5,1)(4.5,0)
\psccurve[linestyle=dashed,dash=3pt 2pt](4.5,5)(5.1,4.25)(4.5,3.5)(3.9,4.25)
\psccurve[linestyle=dashed,dash=3pt 2pt](4.5,2.5)(5.1,1.75)(4.5,1)(3.9,1.75)
%la grasse a droite qui descend
\psecurve[linewidth=2pt]{->}(6.9,10.7)(6.9,10.7)(6,7.1)(5.4,6.8)
\psecurve[linewidth=2pt](2.4,7.4)(2.4,7.4)(3.5,7.1)(5.4,6.8)
\rput[t](3.5,7){$p_4$}
\psframe*[linecolor=white](6.7,6.6)(9.9,7.6)
\rput(8.1,7.1){$[p_4,p_r(p_4)]$}
\psecurve[linewidth=2pt]{->}(4.5,3.5)(5.1,4.25)(4.5,5)(3.9,4.25)(4.5,3.5)
\rput[l](5.2,4.25){$p_2$}
\psframe*[linecolor=white](0,3.65)(3.8,4.85)
\rput(1.9,4.25){$[p_2,p_r(p_2)]$}
%\rput[r](2.9,4.25){$[p_2,p_r(p_2)]$}
%\rput(4.5,6.1){$\Gamma$}
\psecurve(3.9,11.3)(4.8,11)(3.9,10.7)(3,11)(3.9,11.3)
\end{pspicture}
\begin{pspicture}[.4](-1,-.2)(5,7) 
\psset{xunit=.5cm,yunit=.5cm} 
\psecurve{->}(0,11)(.9,11.3)(1.8,11)(.9,10.7)(0,11)(.9,11.3)(1.8,11)
%\rput[br](.9,11.4){$a \times \{0\}$}
\psecurve{->}(3,11)(3.9,11.3)(4.8,11)(3.9,10.7)(3,11)(3.9,11.3)(4.8,11)
%\rput[b](3.9,11.4){$a \times \{\frac{2}{2k-1}\}$}
\psecurve{->}(6,11)(6.9,11.3)(7.8,11)(6.9,10.7)(6,11)(6.9,11.3)(7.8,11)
\pscurve(1.8,11)(2.4,7.4)(3,11)
\pscurve{-}(4.8,11)(5.4,6.8)(6,11)
\psccurve(4.5,5)(4.8,4.25)(4.5,3.5)(4.2,4.25)
\psccurve(4.5,2.5)(4.8,1.75)(4.5,1)(4.2,1.75)
%la courbe qui fait tout le tour
\pscurve(0,11)(.6,7.4)(2.8,5)(3,4.25)(3,3)(3,1.75)(3.2,1)(4.5,0)(5.8,1)(6,1.75)(6,3)(6,4.25)
(6.1,5)(7.2,6.8)(7.8,11)
\psecurve[linewidth=2pt, border=1pt](4.8,11.1)(3.9,10.7)(3,10.3)(2.7,11)(3.9,11.6)(7.5,11.6)(6.9,10.7)(6,10.3)
\psecurve[linewidth=2pt]{->}(1.8,11.1)(.9,10.7)(0,10.3)(-.3,11)(1.2,11.9)(4.5,11.6)
\psecurve[linewidth=2pt, border=1pt](-.3,11)(1.2,11.9)(4.5,11.6)(3.9,10.7)(3,10.3)
\psecurve[linestyle=dashed,dash=3pt 2pt]{*-*}(.9,10.7)(.9,10.7)(2.4,7.4)(3.9,10.7)
%la grasse a gauche qui monte
\psecurve[linewidth=2pt]{*-}(.9,10.7)(.9,10.7)(1.2,9.1)(2.4,7.4)
\psecurve[linestyle=dashed,dash=3pt 2pt]{*-}(3.9,10.7)(3.9,10.7)(2.4,7.4)(.9,10.7)
\pscurve[linestyle=dashed,dash=3pt 2pt]{-*}(2.4,7.4)(5.4,6.8)(6,7.1)(6.9,10.7)
%ligne qui descend du point de connexion le plus bas
\psline[linestyle=dashed,dash=3pt 2pt]{*-*}(5.4,6.8)(4.5,5)
\psline[linestyle=dashed,dash=3pt 2pt]{*-*}(4.5,3.5)(4.5,2.5)
\psline[linestyle=dashed,dash=3pt 2pt]{*-*}(4.5,1)(4.5,0)
\psccurve[linestyle=dashed,dash=3pt 2pt](4.5,5)(5.1,4.25)(4.5,3.5)(3.9,4.25)
\psccurve[linestyle=dashed,dash=3pt 2pt](4.5,2.5)(5.1,1.75)(4.5,1)(3.9,1.75)
%la grasse a droite qui descend
\psecurve[linewidth=2pt]{->}(6.9,10.7)(6.9,10.7)(6.6,9.1)(5.4,6.8)
%\psecurve[linewidth=2pt](2.4,7.4)(2.4,7.4)(3.9,7.1)(5.4,6.8)
\rput[r](1.1,9.1){$p_5$}
\psframe*[linecolor=white](7.1,8.6)(10.4,9.6)
\rput(8.6,9.1){$[p_5,p_r(p_5)]$}
\psecurve[linewidth=2pt]{->}(4.5,3.5)(5.1,4.25)(4.5,5)(3.9,4.25)(4.5,3.5)
\rput[l](5.2,4.25){$p_2$}
\psframe*[linecolor=white](0,3.65)(3.8,4.85)
\rput(1.9,4.25){$[p_2,p_r(p_2)]$}
%\rput(4.5,6.1){$\Gamma$}
\psecurve(3.9,11.3)(4.8,11)(3.9,10.7)(3,11)(3.9,11.3)
\end{pspicture}$$

For $t \in [-6 -3g(S),0]$, set $$S_t=h_S^{-1}( [-6 -3g(S),t]) \subseteq S.$$
$$S_{\ell}(t)=\{(x,0,p(x),\frac1{2k-1}); x \in S_t\},$$
$$A_{\ell}(t)=-\cup_{c \; \mbox{\small component of $\partial S_t$}} c \times 
[p(c), p_r(t)] \times  \{\frac1{2k-1}\},$$
and define the cycle
$$C_{\ell}(t)= -\left( S_t \times \{p_r(t)\} \times \{\frac1{2k-1}\} \right) \cup  S_{\ell}(t) \cup A_{\ell}(t).$$
In a symmetric way, set

$$S_{r}(t)=\{(p(x),0,x \times \{\frac1{2k-1}\}); x \in S_t\},$$
$$A_{r}(t)=\cup_{c \; \mbox{\small component of $\partial S_t$}} [p(c), p_{\ell}(t)] \times   c \times \{\frac1{2k-1}\},$$
and define the cycle
$$C_{r}(t)= - \left( {p_{\ell}(t)} \times S_t \times \{\frac1{2k-1}\} \right) \cup S_{r}(t) \cup A_{r}(t).$$
Let $n$ denote the unit normal positive vector field of $S$.
Set
$$C(0)=-S_{\ell}(0) \cup -S_r(0) \cup \mbox{diag}(n)(S) \cup A(\partial S)$$
$$ \cup \cup_{c \;\mbox{\tiny component of $\partial S$}}T(c, c \times \frac1{2k-1}).$$

Note that the homology class of $G(a)$ that has been defined in Sublemma~\ref{sublemG} is
$$[G(a)]=[C(0)] + (g(S)+k-1) [S^2_{\mbox{\tiny diag}}] + 
[C_{\ell}(0)] + [C_r(0)].$$
Therefore, the proof of Lemma~\ref{lemkey} reduces to the proof of the two following lemmas. 
\begin{lemma}
\label{lemcyc1}
$[C(0)]=0 \;\; \mbox{in} \;\; H_2(C_2(M)).$
\end{lemma}
\begin{lemma}
\label{lemcyc3}
$[C_{\ell}(0)] + [C_r(0)]=-(g(S)+k-1) [S^2_{\mbox{\tiny diag}}].$
\end{lemma}

\noindent{\sc Proof of Lemma~\ref{lemcyc1}:}
Set $$\mbox{diag}(\frac1{2k-1})(S_t)=\{(x,0,x,\frac1{2k-1}); x \in S_t\}$$
and $$\tilde{E}(t)=\cup_{c \; \mbox{\small component of $\partial S_t$}}T(c, c \times \frac1{2k-1}).$$
Then $C(0)$ is homologous to $\tilde{C}(0)$ with
$$\tilde{C}(t)=-S_{\ell}(t) \cup -S_r(t) \cup \mbox{diag}(\frac1{2k-1})(S_t) \cup \tilde{E}(t).$$
We shall now replace $\tilde{E}(t)$ by a chain $E(t)$ such that
\begin{enumerate} 
\item $\partial E(t)= \partial \tilde{E}(t)=\partial S_{\ell}(t) \cup \partial S_r(t) \cup -\partial \mbox{diag}(\frac1{2k-1})(S_t)$\\
\item
$(E(0)-\tilde{E}(0))$ is null-homologous in $C_2(M)$,\\
\item $E(t)$ is defined continuously enough so that if
$$D(t)=-S_{\ell}(t) \cup -S_r(t) \cup \mbox{diag}(\frac1{2k-1})(S_t) \cup E(t),$$
then $\cup_{t \in [-6-3g,0]} D(t)$ may be considered as a cobordism from $D(0)$ to $0$, and used to show that $D(0)$ is null-homologous.
\end{enumerate}
Since the first two conditions imply that $D(0)$ is homologous to $C(0)$, this will be enough to finish proving the lemma.

In order to define $E(t)$, we use our graph $\Gamma$ on $S$.
This graph equips each component $c$ of $\partial S_t$ with the basepoint $p(c)$. We furthermore equip the set of connected components of $\partial S_t$
with the total order from left to right in the picture. 

Then we set 
$$E(t)=\tilde{E}(t) \cup \bigcup_{\{(c,c^{\prime}); c,c^{\prime} \;\mbox{\small components of}\; \partial S_t;   c^{\prime}<c\}} c \times\{0\}
\times c^{\prime} \times\{\frac1{2k-1}\}$$

Since the linking number of $a$ and a curve parallel to $a$ on $\partial A$
is zero, the homology classes of the tori $a \times \{\frac{2i}{2k-1}\} \times a \times \{\frac{2j+1}{2k-1}\}$
are null in $C_2(M)$. Therefore, $(E(0)-\tilde{E}(0))$ is null-homologous in $C_2(M)$.
It is also clear that adding tori has not changed the boundary of $\tilde{E}(t)$.

Since $D(-6-3g)$ is supported in a point, it is null-homologous.
To conclude, we prove that for any subinterval $[h_1,h_2]$ of $[-6-3g,0]$,
$D(h_1)$ is homologous to $D(h_2)$. 
This is clear when $[h_1,h_2]$ does not contain any critical value of $h_S$.
Let us prove that this is still true when $[h_1,h_2]$ contains exactly one critical value $h_c$ in its interior.
%$\cup_{t \in [-6-3g,0]} D(t)$ provides a cobordism from $D(0)$ to $D(-6-3g)$ %that is reduced to a point. 
%It is clear that if a subinterval $[h_1,h_2]$ of $[-6-3g,0]$ 
%does not contain any critical value of $h_S$, 
%then $\cup_{t \in [h_1,h_2]} D(t)$ 
%is a cobordism between $D(h_1)$ and $D(h_2)$. 

There are two cases. Either the number of components of 
$h_S^{-1}(]h_c,h_2])$ is greater than the number of components of 
$h_S^{-1}([h_1,h_c[)$, or it is smaller. In the first case, the corresponding index one
critical point will be called a {\em positive saddle point\/}, in the second case, the critical point  will be called a {\em negative saddle point\/}. 
Let us treat the case of a positive saddle point.
%Because of the symmetry we study only the case of a positive saddle point, 
%and we work under these hypotheses from now on.

Let $r$ be the number of components of $D(h_2)$, we are going to use a continuous map $$f: \sqcup_{i=0}^{r-1}[2i,2i+1] \times [h_1, h_2]
\longrightarrow S,$$ to parametrise
$h_S^{-1}([h_1,h_2])$.
Our parametrisation $f$ has the following properties:
\begin{enumerate}
\item $h_S \circ f(x,t)=t$ on $\sqcup_{i=0}^{r-1}[2i,2i+1] \times [h_1, h_2]$,
\item $f(.,t)$ provides a homeomorphism
from $\frac{[2i,2i+1]}{2i \sim 2i+1}$ to the $(i+1)^{\mbox{\small th}}$ component from left to right of $\partial S_t$, for any $(i,t) \in \{0,1, \dots, r-1\} \times ]h_c, h_2]$,
\item $f(2i,t)=f(2i+1,t) \in \Gamma$ for any $(i,t) \in \{2, 3, \dots, r-1\} \times [h_1, h_2]$, and for any $(i,t) \in \{0,1, \dots, r-1\} \times [h_c, h_2]$,
\item $f(0,t)=f(3,t) \in \Gamma$ and $f(1,t)=f(2,t)$ for any $t \in [h_1,h_c]$,
\item $f(.,t)$ provides a homeomorphism
from the circle $\frac{[0,1] \coprod [2,3]}{0 \sim 3, 1 \sim 2}$ to the first component of $\partial S_t$ for any $t \in [h_1,h_c[$,
$f(.,t)$ provides a homeomorphism
from $\frac{[2i,2i+1]}{2i \sim 2i+1}$ to the $i^{\mbox{\small th}}$ component from left to right of $\partial S_t$ for any $t \in [h_1,h_c]$, for any $i \in \{2,3,\dots,r-1\}$.
\end{enumerate}
Let $T_r$ be the following part of $\RR^2$:
$$T_r =\{(v,w) \in (\sqcup_{i=0}^{r-1}[2i,2i+1])^2; v \geq w\}.$$

\begin{center}
$$ \begin{pspicture}[.2](0,0)(2.5,2.5)
%\psset{xunit=.5cm,yunit=.5cm}
%\rput[tr](3.8,.8){$v$}
%\rput[tr](.8,3.8){$ w$}
\rput(1,2){$T_3$}
\pspolygon*[linecolor=lightgray](0,0)(.5,0)(.5,.5)
\pspolygon*[linecolor=lightgray](1,1)(1.5,1)(1.5,1.5)
\pspolygon*[linecolor=lightgray](2,2)(2.5,2)(2.5,2.5)
\pspolygon*[linecolor=lightgray](1,0)(1.5,0)(1.5,.5)(1,.5)
\pspolygon*[linecolor=lightgray](2,0)(2.5,0)(2.5,.5)(2,.5)
\pspolygon*[linecolor=lightgray](2,1)(2.5,1)(2.5,1.5)(2,1.5)
\end{pspicture}$$
\end{center}

Then $\cup_{t \in [h_1,h_2]} E(t)$ is the image of the continuous map
$$\begin{array}{llll}F:&T_r \times [h_1,h_2] & \longrightarrow & S \times \{0\} \times S \times \{\frac1{2k-1}\}\\&(v,w,t) & \mapsto & (f(v,t),0,f(w,t),\frac1{2k-1}) \end{array}$$
that may be extended to provide a cobordism between $D(h_1)$ and $D(h_2)$.
The case of a negative saddle point is symmetric.
\eop

\noindent{\sc Proof of Lemma~\ref{lemcyc3}:}
The graph $\Gamma$ is as in the figure.
We observe that $\cup_{t \in [-6-3g(S),0]} C_{\ell}(t)$ may not be considered as a cobordism between
$C_{\ell}(0)$ and $C_{\ell}(-6-3g(S))$ (that is a point) because the chains $A_{\ell}(t)$ are not continuously defined with respect to $t$. More precisely, the jumps occur exactly at the heights of the positive saddle points, because the paths
$[p(c),p_r(t)]$ are not continuously defined near a positive saddle point.
Let $h_1$, $h_2$, \dots, $h_{g+k-1}$ be the heights of the critical values ordered
from top to bottom (in decreasing order) corresponding to the positive saddle points $p_{\ell}(h_1)$, $p_{\ell}(h_2)$, \dots, $p_{\ell}(h_{g+k-1})$. The 
figure 8 horizontal curve $c_{8}(h_i)$ that contains $p_{\ell}(h_i)$ is the union of two topological circles $c_1(h_i)$ and  $c_2(h_i)$. Let $c_1(h_i)$ be the part of $c_{8}(h_i)$ where $[p(c),p_r(t)]$ changes. 
%($c_1(h_i)$ is the left-hand side part of $c_{8}(h_i)$.) 
When approaching $c_1(h_i)$ from above $h_i$, this path approaches some path $[p^+(c_1(h_i)),p_r(h_i)]$ homotopic to the path composition of a loop $\gamma(h_i)$ based at $p_{\ell}(h_i)$ and the path $[p_{\ell}(h_i),p_r(h_i)]$.
The loops $c_1(h_i)$ and  $\gamma(h_i)$ are shown in the following picture.

\begin{pspicture}[.4](0,-.2)(5,7) 
\psset{xunit=.5cm,yunit=.5cm} 
\psecurve{->}(0,11)(.9,11.3)(1.8,11)(.9,10.7)(0,11)(.9,11.3)(1.8,11)
%\rput[br](.9,11.4){$a \times \{0\}$}
\psecurve{->}(3,11)(3.9,11.3)(4.8,11)(3.9,10.7)(3,11)(3.9,11.3)(4.8,11)
%\rput[b](3.9,11.4){$a \times \{\frac{2}{2k-1}\}$}
\psecurve{->}(6,11)(6.9,11.3)(7.8,11)(6.9,10.7)(6,11)(6.9,11.3)(7.8,11)
%Les deux contours interieurs de la surface et ses trous
\pscurve(1.8,11)(2.4,8)(3,11)
\pscurve{-}(4.8,11)(5.4,6.8)(6,11)
\psccurve(4.5,5)(4.8,4.25)(4.5,3.5)(4.2,4.25)
\psccurve(4.5,2.5)(4.8,1.75)(4.5,1)(4.2,1.75)
%Le tour de la surface
\pscurve(0,11)(.5,8)(2,6.8)(3,4.25)(3,3.5)(3,1.75)(3.2,1)(4.5,0)(5.8,1)(6,1.75)(6,3.5)(6,4.25)
(6.1,5)(7.2,6.8)(7.8,11)
%Le chemin en zig-zag
%\psecurve[border=1pt](4.8,11.1)(3.9,10.7)(3,10.3)(2.7,11)(3.9,11.6)(7.5,11.6)(6.9,10.7)(6,10.3)
%\psecurve{->}(1.8,11.1)(.9,10.7)(0,10.3)(-.3,11)(1.2,11.9)(4.5,11.6)
%\psecurve[border=1pt](-.3,11)(1.2,11.9)(4.5,11.6)(3.9,10.7)(3,10.3)
%Gamma
%\psecurve[linewidth=2pt]{*-*}(.9,10.7)(.9,10.7)(2.4,8)(3.9,10.7)
%\psecurve[linewidth=2pt]{*-}(3.9,10.7)(3.9,10.7)(2.4,8)(.9,10.7)
%\pscurve[linewidth=2pt]{-*}(2.4,8)(5.4,6.8)(6,7.1)(6.9,10.7)
%\psline[linewidth=2pt]{*-*}(5.4,6.8)(4.5,5)
%\psline[linewidth=2pt]{*-*}(4.5,3.5)(4.5,2.5)
%\psline[linewidth=2pt]{*-*}(4.5,1)(4.5,0)
%\psccurve[linewidth=2pt](4.5,5)(5.1,4.25)(4.5,3.5)(3.9,4.25)
%\psccurve[linewidth=2pt](4.5,2.5)(5.1,1.75)(4.5,1)(3.9,1.75)
%\rput(4.5,6.1){$\Gamma$}
%Le huit du haut
%\psecurve[linewidth=.5pt](3.4,5.2)(2.8,5)(3.4,4.8)(4.5,5)(4.5,5)
%\psecurve[linestyle=dashed,dash=3pt 2pt](3.4,4.8)(2.8,5)(3.4,5.2)(4.5,5)(4.5,5)
%\psecurve(5.5,5.2)(6.1,5)(5.5,4.8)(4.5,5)(4.5,5)
%\psecurve[linestyle=dashed,dash=3pt 2pt](5.5,4.8)(6.1,5)(5.5,5.2)(4.5,5)(4.5,5)
%\gamma(h_2)
\rput[b](6,12.4){$\gamma(h_2)$}
\psecurve[linestyle=dotted, border=1pt]{->}(4.8,11.1)(3.9,10.7)(3,10.3)(2.7,11)(3.2,11.9)(6,12.2)(7.5,11.6)
\psecurve[linestyle=dotted, border=1pt](3.2,11.9)(6,12.2)(7.5,11.6)(6.9,10.7)(6,10.3)
\pscurve[linestyle=dotted](2.4,8)(5.4,6.8)(6,7.1)(6.9,10.7)
\psecurve[linestyle=dotted](.9,10.7)(2.4,8)(3.9,10.7)(3.9,10.7)
%$\gamma(h_1)
\rput[br](1.7,12){$\gamma(h_1)$}
\psecurve[linestyle=dashed,dash=2pt 3pt]{->}(1.8,11.1)(.9,10.7)(0,10.3)(-.3,11)(1.2,11.9)(4.5,11.6)
\psecurve[linestyle=dashed,dash=2pt 3pt, border=1pt](-.3,11)(1.2,11.9)(4.5,11.6)(3.9,10.7)(3,10.3)
\pscurve[linestyle=dashed,dash=2pt 3pt](.9,10.7)(2.4,8)(3.9,10.7)
%Pour h_1 (2.4,8)
\rput[r](.5,8){$c_1(h_1)$}
\psecurve(1.4,8.2)(.5,8)(1.4,7.8)(2.4,8)
\psecurve{->}(2.4,8)(2.4,8)(1.4,7.8)(.5,8)
\psecurve[linestyle=dashed,dash=3pt 2pt](1.4,7.8)(.5,8)(1.4,8.2)(2.4,8)(2.4,8)
%Pour h_2 (5.4,6.8)
\rput[t](3.4,6.45){$c_1(h_2)$}
\psecurve(3.4,7.05)(2,6.8)(3.4,6.55)(5.4,6.8)
\psecurve{->}(5.4,6.8)(5.4,6.8)(3.4,6.55)(2,6.8)
\psecurve[linestyle=dashed,dash=3pt 2pt](3.4,6.55)(2,6.8)(3.4,7.05)(5.4,6.8)(5.4,6.8)
%Pour h_3
\psframe*[linecolor=white](2,3.85)(3.8,4.65)
\rput[r](3.7,4.25){$\gamma(h_3)$}
\psecurve{->}(4.5,5)(3.9,4.25)(4.5,3.5)(5.1,4.25)(4.5,5)(3.9,4.25)(4.5,3.5)
\rput[t](5.2,3.25){$c_1(h_3)$}
\psecurve{->}(5.2,3.65)(6,3.5)(5.2,3.35)(4.5,3.5)
\psecurve(6,3.5)(5.2,3.35)(4.5,3.5)(4.5,3.5)
\psecurve[linestyle=dashed,dash=3pt 2pt](5.2,3.35)(6,3.5)(5.2,3.65)(4.5,3.5)(4.5,3.5)
%Le $\gamma(h_4)$
\psframe*[linecolor=white](2,1.35)(3.8,2.15)
\rput[r](3.7,1.75){$\gamma(h_4)$}
\psecurve{->}(4.5,2.5)(3.9,1.75)(4.5,1)(5.1,1.75)(4.5,2.5)(3.9,1.75)(4.5,1)
%Le huit du bas
\rput[l](5.9,1){$c_8(h_4)$}
\psecurve{->}(4.5,1)(4.5,1)(3.8,.85)(3.2,1)
\psecurve(4.5,1)(3.8,.85)(3.2,1)(3.8,1.15)
\psecurve[linestyle=dashed,dash=3pt 2pt](3.8,.85)(3.2,1)(3.8,1.15)(4.5,1)(4.5,1)
\psecurve{->}(5.2,1.15)(5.8,1)(5.2,.85)(4.5,1)
\psecurve(5.8,1)(5.2,.85)(4.5,1)(4.5,1)
\psecurve[linestyle=dashed,dash=3pt 2pt](5.2,.85)(5.8,1)(5.2,1.15)(4.5,1)(4.5,1)
% Grande meridienne
%\psecurve(4.5,3.3)(3,3)(4.5,2.7)(6,3)(4.5,3.3)
%\psecurve[linestyle=dashed,dash=3pt 2pt](4.5,2.7)(3,3)(4.5,3.3)(6,3)(4.5,2.7)
% Les autres courbes meridiennes%
%\psecurve(3.6,4.4)(3,4.25)(3.6,4.1)(4.2,4.25)(3.6,4.4)
%\psecurve[linestyle=dashed,dash=3pt 2pt](3.6,4.1)(3,4.25)(3.6,4.4)(4.2,4.25)(3.6,4.1)
%\psecurve(5.4,1.9)(6,1.75)(5.4,1.6)(4.8,1.75)(5.4,1.9)
%\psecurve[linestyle=dashed,dash=3pt 2pt](5.4,1.6)(6,1.75)(5.4,1.9)(4.8,1.75)(5.4,1.6)
\psecurve(3.9,11.3)(4.8,11)(3.9,10.7)(3,11)(3.9,11.3)
\end{pspicture}
\begin{pspicture}[.4](0,-.2)(5,7) 
\psset{xunit=.5cm,yunit=.5cm} 
\psecurve{->}(0,11)(.9,11.3)(1.8,11)(.9,10.7)(0,11)(.9,11.3)(1.8,11)
%\rput[br](.9,11.4){$a \times \{0\}$}
\psecurve{->}(3,11)(3.9,11.3)(4.8,11)(3.9,10.7)(3,11)(3.9,11.3)(4.8,11)
%\rput[b](3.9,11.4){$a \times \{\frac{2}{2k-1}\}$}
\psecurve{->}(6,11)(6.9,11.3)(7.8,11)(6.9,10.7)(6,11)(6.9,11.3)(7.8,11)
%Les deux contours interieurs de la surface et ses trous
\pscurve(1.8,11)(2.4,8)(3,11)
\pscurve{-}(4.8,11)(4.8,8)(5.4,6.8)(6,11)
\psccurve(4.5,5)(4.8,4.25)(4.5,3.5)(4.2,4.25)
\psccurve(4.5,2.5)(4.8,1.75)(4.5,1)(4.2,1.75)
%Le tour de la surface
\pscurve(0,11)(.5,8)(2,6.8)(3,4.25)(3,3.5)(3,1.75)(3.2,1)(4.5,0)(5.8,1)(6,1.75)(6,3.5)(6,4.25)
(6.1,5)(7.2,6.8)(7.8,11)
%Le chemin en zig-zag
%\psecurve[border=1pt](4.8,11.1)(3.9,10.7)(3,10.3)(2.7,11)(3.9,11.6)(7.5,11.6)(6.9,10.7)(6,10.3)
%\psecurve{->}(1.8,11.1)(.9,10.7)(0,10.3)(-.3,11)(1.2,11.9)(4.5,11.6)
%\psecurve[border=1pt](-.3,11)(1.2,11.9)(4.5,11.6)(3.9,10.7)(3,10.3)
%Gamma
%\psecurve[linewidth=2pt]{*-*}(.9,10.7)(.9,10.7)(2.4,8)(3.9,10.7)
%\psecurve[linewidth=2pt]{*-}(3.9,10.7)(3.9,10.7)(2.4,8)(.9,10.7)
%\pscurve[linewidth=2pt]{-*}(2.4,8)(5.4,6.8)(6,7.1)(6.9,10.7)
%\psline[linewidth=2pt]{*-*}(5.4,6.8)(4.5,5)
%\psline[linewidth=2pt]{*-*}(4.5,3.5)(4.5,2.5)
%\psline[linewidth=2pt]{*-*}(4.5,1)(4.5,0)
%\psccurve[linewidth=2pt](4.5,5)(5.1,4.25)(4.5,3.5)(3.9,4.25)
%\psccurve[linewidth=2pt](4.5,2.5)(5.1,1.75)(4.5,1)(3.9,1.75)
%\rput(4.5,6.1){$\Gamma$}
%Le huit du haut
%\psecurve[linewidth=.5pt](3.4,5.2)(2.8,5)(3.4,4.8)(4.5,5)(4.5,5)
%\psecurve[linestyle=dashed,dash=3pt 2pt](3.4,4.8)(2.8,5)(3.4,5.2)(4.5,5)(4.5,5)
%\psecurve(5.5,5.2)(6.1,5)(5.5,4.8)(4.5,5)(4.5,5)
%\psecurve[linestyle=dashed,dash=3pt 2pt](5.5,4.8)(6.1,5)(5.5,5.2)(4.5,5)(4.5,5)
%\gamma(h_2)
%\rput[b](6,12.4){$\gamma(h_2)$}
\psecurve[linestyle=dotted, border=1pt]{->}(4.8,11.1)(3.9,10.7)(3,10.3)(2.7,11)(3.2,11.9)(6,12.2)(7.5,11.6)
\psecurve[linestyle=dotted, border=1pt](3.2,11.9)(6,12.2)(7.5,11.6)(6.9,10.7)(6,10.3)
\pscurve[linestyle=dotted](2.4,8)(5.4,6.8)(6,7.1)(6.9,10.7)
%%\psecurve[linestyle=dotted](.9,10.7)(2.4,8)(3.9,10.7)(3.9,10.7)
%$\gamma(h_1)
%\rput[br](1.7,12){$\gamma(h_1)$}
\rput[b](6,12.4){$\gamma_r(h_2)$}
\psecurve[linestyle=dotted, border=1pt](4.8,11.1)(3.9,10.7)(3,10.3)(2.7,11)(3.2,11.9)(6,12.2)(7.5,11.6)
\psecurve[linestyle=dotted, border=1pt]{->}(6,10.3)(6.9,10.7)(7.5,11.6)(6,12.2)(3.2,11.9)
\pscurve[linestyle=dotted](2.4,8)(5.4,6.8)(6,7.1)(6.9,10.7)
\psecurve[linestyle=dotted](.9,10.7)(.9,10.7)(2.4,8)(3.9,10.7)

%$\gamma_r(h_1)
\rput[br](1.7,12){$\gamma_r(h_1)$}
\psecurve[linestyle=dashed,dash=2pt 3pt](1.8,11.1)(.9,10.7)(0,10.3)(-.3,11)(1.2,11.9)(4.5,11.6)
\psecurve[linestyle=dashed,dash=2pt 3pt, border=1pt]{->}(3,10.3)(3.9,10.7)(4.5,11.6)(1.2,11.9)(-.3,11)
\psecurve[linestyle=dotted](1.8,11.1)(.9,10.7)(0,10.3)(-.3,11)(1.2,11.9)(4.5,11.6)(3.9,10.7)(3,10.3)
\pscurve[linestyle=dashed,dash=2pt 3pt](.9,10.7)(2.4,8)(3.9,10.7)
%Pour h_1 (2.4,8)
\rput[t](3.8,7.7){$c_r(h_1)$}
\psecurve{->}(3.8,8.2)(4.8,8)(3.8,7.8)(2.4,8)
\psecurve(4.8,8)(3.8,7.8)(2.4,8)(2.4,8)
\psecurve[linestyle=dashed,dash=3pt 2pt](3.8,7.8)(4.8,8)(3.8,8.2)(2.4,8)(2.4,8)
%Pour h_2 (5.4,6.8)
\rput[tr](6.4,6.55){$c_r(h_2)$}
\psecurve{->}(6.4,6.95)(7.2,6.8)(6.4,6.65)(5.4,6.8)
\psecurve(7.2,6.8)(6.4,6.65)(5.4,6.8)(5.4,6.8)
\psecurve[linestyle=dashed,dash=3pt 2pt](6.4,6.65)(7.2,6.8)(6.4,6.95)(5.4,6.8)(5.4,6.8)
%Pour h_3
\psframe*[linecolor=white](2,3.95)(3.8,4.75)
\rput[r](3.7,4.35){$\gamma_r(h_3)$}
\psecurve{->}(4.5,5)(3.9,4.25)(4.5,3.5)(5.1,4.25)(4.5,5)(3.9,4.25)(4.5,3.5)
\rput[t](5.2,3.25){$c_r(h_3)$}
\psecurve{->}(5.2,3.65)(6,3.5)(5.2,3.35)(4.5,3.5)
\psecurve(6,3.5)(5.2,3.35)(4.5,3.5)(4.5,3.5)
\psecurve[linestyle=dashed,dash=3pt 2pt](5.2,3.35)(6,3.5)(5.2,3.65)(4.5,3.5)(4.5,3.5)
%Le $\gamma(h_4)$
\psframe*[linecolor=white](2,1.45)(3.8,2.25)
\rput[r](3.7,1.85){$\gamma_r(h_4)$}
\psecurve{->}(4.5,2.5)(3.9,1.75)(4.5,1)(5.1,1.75)(4.5,2.5)(3.9,1.75)(4.5,1)
%Le c du bas
\rput[l](5.9,1){$c_r(h_4)$}
\psecurve{->}(5.2,1.15)(5.8,1)(5.2,.85)(4.5,1)
\psecurve(5.8,1)(5.2,.85)(4.5,1)(4.5,1)
\psecurve[linestyle=dashed,dash=3pt 2pt](5.2,.85)(5.8,1)(5.2,1.15)(4.5,1)(4.5,1)
% Grande meridienne
%\psecurve(4.5,3.3)(3,3)(4.5,2.7)(6,3)(4.5,3.3)
%\psecurve[linestyle=dashed,dash=3pt 2pt](4.5,2.7)(3,3)(4.5,3.3)(6,3)(4.5,2.7)
% Les autres courbes meridiennes%
%\psecurve(3.6,4.4)(3,4.25)(3.6,4.1)(4.2,4.25)(3.6,4.4)
%\psecurve[linestyle=dashed,dash=3pt 2pt](3.6,4.1)(3,4.25)(3.6,4.4)(4.2,4.25)(3.6,4.1)
%\psecurve(5.4,1.9)(6,1.75)(5.4,1.6)(4.8,1.75)(5.4,1.9)
%\psecurve[linestyle=dashed,dash=3pt 2pt](5.4,1.6)(6,1.75)(5.4,1.9)(4.8,1.75)(5.4,1.6)
\psecurve(3.9,11.3)(4.8,11)(3.9,10.7)(3,11)(3.9,11.3)
\end{pspicture}

Define $C^+_{\ell}(h_i)$ from $C_{\ell}(h_i)$ by replacing $c_{8}(h_i) \times [p(c_{8}(h_i)),p_r(h_i)] \times \{\frac1{2k-1}\}$ by the union
$$\left(c_{1}(h_i) \times [p^+(c_1(h_i)),p_r(h_i)] \times \{\frac1{2k-1}\}\right) \cup$$ $$\cup \left(c_{2}(h_i) \times [p(c_{8}(h_i)),p_r(h_i)] \times \{\frac1{2k-1}\} \right)$$

Set $h_0=0$, then for any $i=1, \dots, g+k-1$, $C^+_{\ell}(h_i)$ is homologous to $C_{\ell}(h_{i-1})$.
%$$\cup_{t \in ]h_i,h_{i-1}]} C_{\ell}(t) \cup C^+_{\ell}(h_i)$$
%is a cobordism that identifies the homology classes of 
%$C_{\ell}(h_{i-1})$ and $C^+_{\ell}(h_i)$. 
Furthermore,
$$[C^+_{\ell}(h_i)]=[C_{\ell}(h_i)]- [c_1(h_i) \times 
\left(\gamma(h_i) \times \{\frac1{2k-1}\}\right)].$$

This shows that

$$[C_{\ell}(0)] = -\sum_{i=1}^{g(S)+k-1} [c_1(h_i) \times \left(\gamma(h_i) \times \{\frac1{2k-1}\}\right)].$$

Similarly, define $c_r(h_i)$ as the part of $c_{8}(h_i)$ where $[p(c),p_{\ell}(t)]$ changes. This part $c_r(h_i)$ is the right-hand side of $c_{8}(h_i)$ in the figure, it coincides with $c_1(h_i)$
when $i \geq k$. 
Define $\gamma_r(h_i)$ as the loop based at $p_{\ell}(h_i)$ that is the limit of $[p(c),p_{\ell}(t)]$ when $c$ approaches $c_r$ from above. 
%and let $\gamma_r(h_i)$ be the loop based at $p_{\ell}(h_i)$ 
% such that $[p^+(c_r(h_i)),p_{\ell}(h_i)]$ is homotopic to the 
%path composition of $\gamma_r(h_i)$ and $[p(c_{8}(h_i)),p_{\ell}(h_i)]$.
Again, $\gamma_r(h_i)=\gamma(h_i)$ when $i \geq k$.
Then, as above, we get that 
$$[C_r(0)] = \sum_{i=1}^{g+k-1} [\gamma_r(h_i) \times c_r(h_i) \times \{\frac1{2k-1}\}].$$
Thus, in $C_2(M)$, we have
$$[C_{\ell}(0)] + [C_r(0)]=$$
$$\sum_{i=1}^{g+k-1} \left( \ell(\gamma_r(h_i), c_r(h_i) \times \{\frac1{2k-1}\}) - \ell(c_1(h_i), \gamma(h_i) \times \{\frac1{2k-1}\}) \right)[S^2_{\mbox{\tiny diag}}].$$

When $i \geq k$, $$\ell(\gamma_r(h_i), c_r(h_i) \times \{\frac1{2k-1}\}) - \ell(c_1(h_i), \gamma(h_i) \times \{\frac1{2k-1}\})$$
$$=\ell(\gamma_r(h_i), c_r(h_i) \times \{\frac1{2k-1}\}) - \ell(\gamma_r(h_i), c_r(h_i) \times \{-\frac1{2k}\})$$
$$=-\langle \gamma_r(h_i),c_r(h_i)\rangle_S=-1.$$
When $i < k$, 
$$\ell(\gamma_r(h_i), c_r(h_i) \times \{\frac1{2k-1}\})$$
$$=\ell(c_r(h_i), \gamma_r(h_i) \times \{-\frac1{2k}\})$$
where $c_r(h_i)$ is homologous to $a \times \{\frac{2i}{2k-1}\}$ inside  a subsurface of $S$ that does not intersect $(\gamma_r(h_i) \times \{-\frac1{2k}\})$. Therefore, this expression equals
$$\ell(1 \times a, \gamma_r(h_i) \times \{-\frac1{2k}\})$$
$$=\frac1k \ell (1 \times \partial S,\gamma_r(h_i) \times \{-\frac1{2k}\})$$
$$=\frac1k\langle S \cup ([0,1] \times \partial S),\gamma_r(h_i) \times \{-\frac1{2k}\}\rangle_M$$
$$=\frac1k\langle \partial S,\gamma_r(h_i) \times \{-\frac1{2k}\}\rangle_{\partial A} =-\frac{i}k$$
since $\gamma_r(h_i) \times \{-\frac1{2k}\}$ meets $S$ only along its boundary.

When $i < k$, $c_1(h_i)$ is homologous to $\sqcup_{j=0}^{i-1} a \times \{\frac{2j}{2k-1}\} $ inside a subsurface of $S$
that does not intersect $(\gamma(h_i) \times \{\frac1{2k-1}\})$.
Therefore,
$$\ell(c_1(h_i) , \gamma(h_i) \times \{\frac1{2k-1}\})$$
$$=i\ell(1 \times a, \gamma(h_i) \times \{\frac1{2k-1}\})$$
$$=\frac{i}k\langle S \cup ([0,1] \times \partial S),\gamma(h_i) \times \{\frac1{2k-1}\}\rangle_M$$
$$=\frac{i}k\langle \partial S,\gamma(h_i) \times \{\frac1{2k-1}\}\rangle_{\partial A} =\frac{i}k .$$
Therefore,
$$\sum_{i=1}^{k-1} \left( \ell(\gamma_r(h_i), c_r(h_i) \times \{\frac1{2k-1}\}) - \ell(c_1(h_i) , \gamma(h_i) \times \{\frac1{2k-1}\}) \right)$$
$$=-2\sum_{i=1}^{k-1} \frac{i}k =1-k.$$

$$[C_{\ell}(0)] + [C_r(0)]=-(g(S)+k-1) [S^2_{\mbox{\tiny diag}}].$$
\eop
\newpage

\section{Comparison with the Walker invariant}
\setcounter{equation}{0}
\label{seccompwal}

This section is devoted to proving Theorem~\ref{thcompwal} as an application of the splitting formula for $Z_1$.
Recall that Kevin Walker has extended the Casson invariant to rational homology spheres in \cite{wal}. Let 
$\lambda_W$ \index{NN}{lambdaW@$\lambda_W$} denote his invariant normalized as in \cite{wal}. Here, we use the normalisation $\lambda=\frac{\lambda_W}{2}$ so that $\lambda$ \index{NN}{lambda@$\lambda$} coincides with the Casson invariant normalised as in \cite{akmc,gm2} for integral homology spheres. Thus, this section is devoted to proving
that for any rational homology sphere $M$,
$$Z_{1}(M)=\frac{\lambda(M)}{2}[\tata].$$

\subsection{Sketch of the proof of Theorem~\ref{thcompwal}}

Since $[\theta=\tata]$ freely generates $\CA_1(\emptyset)$, there exists a rational
invariant $\nu$ of rational homology spheres such that $Z_{1}=\frac{\nu}{2}[\theta].$

Recall from \inconstproporrev ~and Theorem~\ref{mainth} that $\nu$ satisfies:
\begin{itemize}
\item For any rational homology sphere, $\nu(-M)=-\nu(M)$.
\item For any rational generalised 2-clover, $\nu(D)=\frac{\ell(D;\theta)}{6}$.
\end{itemize}

Then Theorem~\ref{thcompwal} is the direct consequence of the two following 
propositions.

\begin{proposition}
\label{propcompwun}
\begin{itemize}
\item For any rational homology sphere, $\lambda(-M)=-\lambda(M)$.
\item For any rational generalised 2-clover, $\lambda(D)=\frac{\ell(D;\theta)}{6}$.
\end{itemize}
\end{proposition}

\begin{proposition}
\label{propcompwde}
Let $\kappa$ be a rational-valued invariant of rational homology spheres that satisfies the two following properties:
\begin{itemize}
\item For any rational homology sphere, $\kappa(-M)=-\kappa(M)$.
\item For any rational generalised 2-clover, $\kappa(D)=0$.
\end{itemize}
Then $\kappa$ vanishes identically.
\end{proposition}

\begin{remark}
For integral homology spheres, Proposition~\ref{propcompwde} is a consequence of the computation of the Goussarov-Habiro filtration of the space of integral homology spheres, that implies that the only rational invariants of integral homology spheres that vanish on integral generalised degree 2 clovers are constant, and 
therefore vanish provided that they vanish at the sphere $S^3$.
For rational homology spheres, the invariant $\mbox{Log}(|H_1|)$ is a non-constant invariant that vanishes at $S^3$ and on rational generalised degree two clovers, and we need the assumption on
the behaviour under orientation reversing, and a proof that is given in Subsection~\ref{sublemfilrat}.
\end{remark}

\noindent{\sc Proof of Proposition~\ref{propcompwun}}.
The first property of $\lambda$ is well-known \cite[Lemma 3.1]{wal}.
Let $D=(M;2;(A,B_0);(A^{\prime}, B^{\prime}_0))$ be a generalised clover. 
Let $B=M \setminus \mbox{Int}(A)$.
Let $B^{\prime}$ be obtained from $B$ by replacing $B_0$ by $B^{\prime}_0$.
Let $([0,1] \times \partial A)$ denote a collar of $\partial A$, disjoint from $B_0$, in $B$. 
Let $B_1= B \setminus \left( \partial A \times[0,1[ \right)$ and  
$B^{\prime}_1= B^{\prime} \setminus \left( \partial A \times[0,1[ \right)$.
Then $D_1=(M;2;(A,B_1);(A^{\prime}, B^{\prime}_1))$ is a rational generalised clover
such that 
$$\lambda(D)= \lambda(A \cup B) - \lambda(A^{\prime} \cup B)-
\lambda(A \cup B^{\prime})+\lambda(A^{\prime} \cup B^{\prime})=\lambda(D_1)$$ is computed in \cite[Theorem 1.3]{les}.

Similarly, $\nu(D_1)=\nu(D)$, and therefore
$\ell(D_1;\theta)=\ell(D;\theta)$.
Thus, it is enough to prove that $\lambda(D)=\ell(D_1;\theta)/6$. 
Let us compute $\ell(D_1;\theta)$.
Let $(\alpha_1,\dots, \alpha_g)$ and $(\beta_1,\dots, \beta_g)$ be two bases for
$\CL_A$ and $\CL_B$, respectively, such that $< \alpha_i,\beta_j>_{\partial A}$ is equal to
the Kronecker symbol $\delta_{ij}$ for any $i,j$ in $\{1,\dots,g\}$.
Then $\CI(A,A^{\prime})$ reads
$$\sum_{f:\{1,2,3\}\rightarrow \{1,2,\dots,g\}}
\CI(A,A^{\prime})(\alpha_{f(1)} \wedge \alpha_{f(2)} \wedge \alpha_{f(3)})
\beta_{f(1)}^{\ast} \otimes \beta_{f(2)}^{\ast} \otimes
\beta_{f(3)}^{\ast}$$
where $\beta_{f(i)}^{\ast}=<.,\beta_{f(i)}>_{\partial A}$, and $\CI(B_1,B_1^{\prime})=\CI(B,B^{\prime})$ reads
$$\sum_{h:\{1,2,3\}\rightarrow \{1,2,\dots,g\}}
\CI(B,B^{\prime})(\beta_{h(1)} \wedge \beta_{h(2)} \wedge \beta_{h(3)})
\alpha_{h(1)}^{\ast} \otimes \alpha_{h(2)}^{\ast} \otimes
\alpha_{h(3)}^{\ast}$$
where $\alpha_{h(i)}^{\ast}=<.,\alpha_{h(i)}>_{\partial B_1}=<\alpha_{h(i)},.>_{\partial A}$.
Let $\sigma: V(\theta) \rightarrow \{A,B_1\}$ be a bijection.
$$ \begin{pspicture}[.2](-1,0)(4,2)
\psline{*-*}(0,1)(3,1)
\pscurve{-}(0,1)(1,1.8)(2,1.8)(3,1)
\pscurve{-}(0,1)(1,.2)(2,.2)(3,1)
\rput[b](1,.3){\small 1}
\rput[b](1,1.1){\small 2}
\rput[b](1,1.9){\small 3}
\rput[b](2,.3){\small 3}
\rput[b](2,1.1){\small 2}
\rput[b](2,1.9){\small 1}
\rput[r](-.1,1){$\sigma^{-1}(A)$}
\rput[l](3.1,1){$\sigma^{-1}(B_1)$}\end{pspicture}$$

Then $\ell(D_1;\theta;\sigma)$ is the contraction of the tensor
$\CI(A,A^{\prime}) \otimes \CI(B,B^{\prime})$ with respect to the linking number and the edge data. This contraction maps $$<.,\beta_{f(i)}>_{\partial A} \otimes
<.,\alpha_{h(4-i)}>_{\partial B_1}\;\; \mbox{to}\;\; \ell(\beta_{f(i)},\alpha_{h(4-i)} \times \{1\})= \delta_{f(i)h(4-i)}.$$
Therefore, $\ell=\ell(D_1;\theta;\sigma)$ equals
$$\sum_{f:\{1,2,3\}\rightarrow \{1,2,\dots,g\}}
\CI(A,A^{\prime})(\alpha_{f(1)} \wedge \alpha_{f(2)} \wedge \alpha_{f(3)})
\CI(B,B^{\prime})(\beta_{f(3)} \wedge \beta_{f(2)} \wedge \beta_{f(1)})$$
$$=-6\sum_{\{i,j,k\} \subset \{1,2,\dots,g\}}
\CI(A,A^{\prime})(\alpha_{i} \wedge \alpha_{j} \wedge \alpha_{k})
\CI(B,B^{\prime})(\beta_{i} \wedge \beta_{j} \wedge \beta_{k}).$$
Since $\ell(D_1;\theta)=2\ell(D_1;\theta;\sigma)$,
$\ell(D_1;\theta)/6$ is equal to the right-hand side of the formula for 
$\lambda(D)$ in \cite[Theorem 1.3]{les}.
\eop

\subsection{Partial review of clover theory}
In order to prove Proposition~\ref{propcompwde}, we need to recall a few
facts about {\em clovers.\/} \index{TT}{clover}
 
\begin{figure}[!h]
\begin{center}
\includegraphics[scale=0.2, height=3cm]{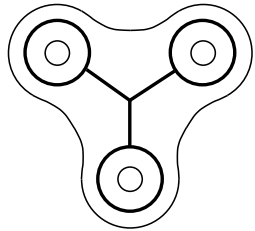}
\caption{The planar unoriented surface $N(Y)$} \label{fign(y)}
\end{center}
\end{figure}

Consider the planar unoriented surface $N(Y)$ of Figure~\ref{fign(y)} that is a two-dimensional neighborhood of the bold {\em \indexT{Y-graph}\/} whose {\em leaves\/} \index{TT}{Y-graph!leaf of} are the three bold circles. 
%$$\begin{pspicture}[.2](-2,-2)(2,1)
%\psline[linewidth=2pt](0,-.05)(0,-.6)
%\psline[linewidth=2pt](0,-.05)(-.5,.2)
%\psline[linewidth=2pt](0,-.05)(.5,.2)
%\pscircle[linewidth=2pt](1,.2){.5}
%\pscircle[linewidth=2pt](-1,.2){.5}
%\pscircle[linewidth=2pt](0,-1.1){.5}
%\pscircle(1,.2){.3}
%\pscircle(-1,.2){.3}
%\pscircle(0,-1.1){.3}
%\psccurve(-.7,-.7)(-.5,-1.7)(.5,-1.7)(.7,-.7)(1.4,-.3)(1.5,.7)(0,.55)(-1.5,.7)(-1.4,-.3)
%\rput[b](1,.3){\small 1}
%\end{pspicture}$$

\begin{figure}[!h]
\begin{center}
\includegraphics[scale=0.2, height=3cm]{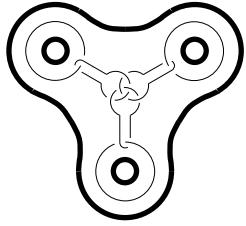}
\caption{Surgery link associated to a framed Y-graph} \label{figbor}
\end{center}
\end{figure}

A {\em degree $k$ \indexT{clover}\/} is an embedding $\phi$ of a disjoint union of $k$ copies $N(Y)^{(i)}$ of $N(Y)$, $i \in \{1,\dots,k\}$ in a rational homology sphere $M$.
With such a data we associate the integral generalised clover $$D=(M;k;(A^i)_{i \in \{1,\dots,k\}},(B^i)_{i \in \{1,\dots,k\}})$$ where $A^i$ is a regular neighborhood of $\phi(N(Y)^{(i)})$ and $B^i$ is obtained from $A^i$ by surgery on the six-component framed link in Figure~\ref{figbor}
%%%Figure~\ref{}
that inherits its framing from the embedding of $N(Y)$. The surgery on $D$ 
consists in replacing every $A^i$ by the corresponding $B^i$. This surgery does not depend on the orientations of the $N(Y)^{(i)}$. 

\bigskip

An orientation of the surface $N(Y)$ 
%an orientation of its outer boundary that in turn 
induces an orientation of the leaves and a cyclic order on them.
When the pictured $N(Y)$ is given the orientation induced by the standard orientation of the plane, the induced orientation and order are the counterclockwise
ones.
Changing the orientation of $N(Y)$ changes both the orientation of the leaves and their cyclic order. A $Y$-graph equipped with such an orientation is said to
be {\em oriented.\/}

\begin{figure}[!h]
\begin{center}
\includegraphics[height=2cm]{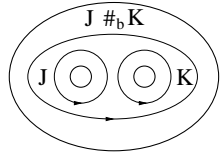}
\caption{Band sum of two framed knots} \label{2hole}
\end{center}
\end{figure}

A framed knot $J\sharp _b K$ is a {\em band sum\/} of two framed oriented knots
$J$ and $K$ if there exists an embedding of a $2$-hole disk
%\begin{itemize}
%\item 
that factors the three knot embeddings $J$, $K$ and $J\sharp _b K$ by the embeddings of the three
curves pictured in \mbox{Figure \ref{2hole}} into the $2$-hole disk, and
%\item 
that induces the three framings.
%\end{itemize} 

%Note that
%$$\ell(K^1\sharp _b K^2,K^1\sharp _b K^2)=
%\ell(K^1,K^1)+\ell(K^2,K^2)+2 \ell(K^1,K^2).$$

We shall use the following lemma that allows one to cut a leaf of a clover.

\begin{lemma}  {\bf (\cite[Theorem 3.1]{ggp},\cite[Lemma 4.15]{al})}
\label{lemcutleaf}
Let $K_1$, $K_2$, $K_3$ be three framed knots in an oriented 3-manifold $M$
that are the leaves of an oriented framed $Y$-graph $G$ in $M$. 
Assume that $K_1$ is a band sum of two framed knots $K^{2}_1$ and $K^{3}_1$. 
For $j=2$ and $3$, let $K_j^2$ and $K_j^3$ be two parallels of $K_j$ equipped with the framing $\ell(K_j,K_j)$ of $K_j$, and such that $\ell(K_j^2,K_j^3)=\ell(K_j,K_j)$.
Then there exist
two oriented disjoint framed $Y$-graphs $G^{2}$ and $G^{3}$ in $M$ whose framed leaves are
$K^{2}_1$, $K^{2}_2$, $K^{2}_3$ and $K^{3}_1$, $K^{3}_2$, $K^{3}_3$, respectively such that the surgery along $G$ is equivalent to the surgery along $G^{2} \cup G^{3}$.
\end{lemma}

Note that under the above assumptions, we have:
\begin{itemize}
\item for $i=2$ or $3$, and for $j,k \in \{2,3\}$, $\ell(K^i_j,K^i_k)=\ell(K_j,K_k)$,
\item for $j \in \{2,3\}$, $\ell(K_1,K_j)=\ell(K^2_1,K^2_j) + \ell(K^3_1,K^3_j)$,
\item $\ell(K_1,K_1)=\ell(K^2_1,K^2_1) + \ell(K^3_1,K^3_1) +2 \ell(K^2_1,K^3_1)$.
\end{itemize}

\subsection{Variation of $\kappa$ under surgery along a $Y$-graph}
\label{subvarrho}
%\begin{lemma}
%\label{lemclastriv}
%If a parallel with respect to the framing of some leaf of a $Y$-graph $G$ %bounds a disk in the complement of $G$, then the surgery along $G$ is trivial.
%\end{lemma}

In this subsection, we prove the following proposition.

\begin{proposition}
\label{propvary}
A real-valued invariant $\kappa$ of rational homology spheres that vanishes at rational generalised degree 2 clovers does not vary under a surgery along a $Y$-graph.
\end{proposition}

Since $\kappa(D)=0$ for any rational generalised 2-clover, the variation $\kappa^{\prime}$ of
$\kappa$ under surgery on an oriented framed $Y$-graph with leaves $K_1$, $K_2$, $K_3$ only depends on the lagrangian of the exterior of the $Y$-graph 
%$M \setminus \phi(A_1)$, 
that is determined by the linking matrix 
$\left[ \ell(K_i,K_j) \right]_{i,j \in \{1,2,3\}}$ of the framed leaves $K_1$, $K_2$ and $K_3$ of the framed $Y$-graph.

We are going to prove the following lemma:

\begin{lemma}
\label{lemreallink}
Any symmetric matrix $[a_{ij}]_{i,j \in \{1,\dots,n\}}$ with rational coefficients is the linking matrix of some framed link in some rational homology sphere.
\end{lemma}

In particular, any rational symmetric matrix $[a_{ij}]_{i,j \in \{1,2,3\}}$ is the linking matrix of some oriented $Y$-graph, and the variation of $\kappa$ by surgery on that framed graph is denoted by $\kappa^{\prime}([a_{ij}])$.
We are going to study the properties of $\kappa^{\prime}$ to prove that it identically vanishes.

%If a parallel with respect to the framing of some leaf bounds a disk 
%in the complement of a $Y$-graph, then the surgery 
%along this $Y$-graph is trivial.
%Therefore, if $a_{11}=a_{12}=a_{13}=0$, then $\kappa^{\prime}([a_{ij}])=0$  

First note that by symmetry, for any cyclic permutation $\sigma$ of $\{1,2,3\}$, $\kappa^{\prime}([a_{ij}])=\kappa^{\prime}([a_{\sigma(i)\sigma(j)}])$.

As a consequence of Lemmas~\ref{lemcutleaf} and \ref{lemreallink}, for any rational numbers $a_{\ast \ast}$ and $\ell^{23}$ involved below, we have
$$\kappa^{\prime}\left[ \begin{array}{ccc}
a^2_{11} + a^3_{11} + 2 \ell^{23}& a^2_{12} + a^3_{12}& a^2_{13} + a^3_{13}\\
a^2_{12} + a^3_{12}& a_{22} & a_{23}\\
a^2_{13} + a^3_{13}& a_{23} & a_{33}\\
\end{array}\right]=$$ $$\kappa^{\prime}\left[ \begin{array}{ccc}a^2_{11}  & a^2_{12} & a^2_{13}\\
a^2_{12} & a_{22} & a_{23}\\
a^2_{13} & a_{23} & a_{33}\\
\end{array}\right] 
+\kappa^{\prime}\left[ \begin{array}{ccc}a^3_{11}  & a^3_{12} & a^3_{13}\\
a^3_{12} & a_{22} & a_{23}\\
a^3_{13} & a_{23} & a_{33}\\
\end{array}\right].$$

Applying this when $a^3_{11}= a^3_{12} = a^3_{13}=\ell^{23}=0$
shows that, if $a_{11}=a_{12}=a_{13}=0$, then $\kappa^{\prime}([a_{ij}])=0$. 
Applying this property again when $a^3_{11}= a^3_{12} = a^3_{13}=a^2_{11} + 2 \ell^{23}=0$ shows that $\kappa^{\prime}$ does not depend on $a_{11}$. By symmetry, it does not depend on the other diagonal terms either.
Therefore,
$$\kappa^{\prime}([a_{ij}])=\kappa^{\prime}(a=a_{12}, b= a_{13}, c=a_{23})$$
where 
$$\kappa^{\prime}(a,b,c)=\kappa^{\prime}(c, a, b),$$
$$\kappa^{\prime}(0,0,c)=0,$$
and 
$$\kappa^{\prime}(a+a^{\prime},b+b^{\prime},c)=\kappa^{\prime}(a,b, c)+\kappa^{\prime}(a^{\prime},b^{\prime},c).$$
In particular,
$$\begin{array}{ll}\kappa^{\prime}(a,b,c)&=\kappa^{\prime}(a,0,c)+\kappa^{\prime}(0,b,c)=\kappa^{\prime}(c,a,0)+\kappa^{\prime}(b,c,0)\\
&=\kappa^{\prime}(c,0,0)+\kappa^{\prime}(0,a,0)+\kappa^{\prime}(b,0,0)+\kappa^{\prime}(0,c,0)=0.
%%&=a\kappa^{\prime}(1,0,c)+b\kappa^{\prime}(0,1,c)\\
%%&=(a+b)\kappa^{\prime}(1,0,c)
\end{array}$$
%where the last equality comes from the symmetry.
%This implies that $\kappa^{\prime}(a,b,0)=0$, and by symmetry that %$\kappa^{\prime}(1,0,c)=0$, and therefore that $\kappa^{\prime}$ identically %vanishes. 
This ends the proof of Proposition~\ref{propvary} assuming Lemma~\ref{lemreallink}.
\eop

\noindent{\sc Proof of Lemma ~\ref{lemreallink}:}
Let $p>0$ and $q$ denote coprime integers.
 Recall that 
the $H_1$ of the lens space $L(p,-q)$ is isomorphic to $\ZZ/p\ZZ$,
and that one of its generators $[\gamma_{p,q}]$ satisfies 
$$\ell([\gamma_{p,q}], [\gamma_{p,q}]) = \frac{q}{p} \;\;\mbox{mod}\; \ZZ,$$
and therefore when $\gamma_{p,q}$ is a knot equipped with a suitable framing
that represents $[\gamma_{p,q}]$,
$$\ell(\gamma_{p,q}, \gamma_{p,q}) = \frac{q}{p}.$$

Now, for $d>0$, consider the following connected sum
$$M(d)= L(d,-1) \sharp L(d,1) \sharp L(d,-1).$$
and three curves $\alpha$, $\beta$, $\gamma$, one in each factor, such that
$$\ell(\alpha,\alpha)=-\ell(\beta,\beta)=\ell(\gamma,\gamma)=\frac{1}{d}  \;\;\mbox{mod} \;\ZZ.$$
Then, for any integer $k$ we may choose curves $\delta(d,k), \varepsilon(d) \in M(d)$ that
are homologous to $(k\alpha+k\beta)$ and to $(\gamma-\beta)$, respectively,
such that
$$\ell(\delta(d,k),\delta(d,k))=\ell(\varepsilon(d),\varepsilon(d))=0$$
and
$$\ell(\delta(d,k),\varepsilon(d))=\frac{k}{d}.$$

Now, write the $a_{ij}$ as irreducible fractions with positive denominators $a_{ii}={q_i}/{p_i}$ and $a_{ij}=k_{ij}/d_{ij}$.
Set 
$$M= \sharp_{i \{1,2,\dots,n\}} L(p_i,-q_i) \sharp
\sharp_{i,j \in \{1,2,\dots,n\}, i<j}M(d_{ij}).$$
Choose knots $K_i$ that are homologous to
$$ \gamma_{p_i,q_i} + \sum_{j \in \{1,2,\dots,n\}, i<j}\delta(d_{ij},k_{ij})
+ \sum_{j \in \{1,2,\dots,n\}, j<i}\varepsilon(d_{ji}).$$
Then the linking matrix of $(K_1, \dots, K_n)$ is the wanted one modulo $\ZZ$.
Next crossing changes and framing changes adjust it to the wanted one.
\eop

\subsection{Proof of Proposition~\ref{propcompwde}}
\label{sublemfilrat}

Recall the following lemma. 

\begin{lemma}
\label{lemnondeg}
Let $N$ be a compact oriented $3$-manifold with boundary. 
Then $H_2(N,\partial N;\ZZ)$ is isomorphic to $\mbox{Hom}(H_1(N;\ZZ);\ZZ)$
by the isomorphism that maps the homology class $[F]$ of a surface $F$ with boundary in $\partial N$ to the algebraic intersection with $F$.
\end{lemma}
\bp Use the Poincar\'e duality to identify $H_{2}(N,\partial N;\ZZ)$ to $H^{1}(N;\ZZ)$ and the universal coefficients theorem to identify $H^{1}(N;\ZZ)$ to $\mbox{Hom}(H_1(N;\ZZ);\ZZ)$.\eop

Lemma~\ref{lemnondeg} yields the following characterisation of integral homology handlebodies.

\begin{lemma}
\label{lemcharZhan}
Let $H$ be a connected compact oriented $3$-manifold with connected boundary. 
If the map from $H_1(\partial H;\ZZ)$ to $H_1(H;\ZZ)$ induced by the inclusion
is a surjection, then $H$ is an integral homology handlebody.
\end{lemma}
\bp The long exact sequence associated to $(H, \partial H)$ yields the exact
sequence with integral coefficients:
$$H_2(H) \longrightarrow H_2(H,\partial H)\longrightarrow  H_1(\partial H)\longrightarrow H_1(H) \rightarrow 0$$
and shows that $H_1(H, \partial H;\ZZ)=0$. 
Then its dual $H^1(H, \partial H;\ZZ) =H_2(H)$ is also trivial,
and $H_2(H,\partial H)=\mbox{Hom}(H_1(H;\ZZ);\ZZ)$ is a free $\ZZ$-module whose rank is the genus of $\partial H$. 

We are now left with the proof that $H_1(H;\ZZ)$ has no torsion. 
If $H_1(H;\ZZ)$ had torsion, there would exist a primitive element $x$ of $H_1(\partial H;\ZZ)$, and a $k>0$ minimal such that $kx=0$ in $H_1(H)$, and $k>1$.
Under these assumptions,
$kx$ would be in the image of the boundary map 
$$H_2(H,\partial H)\longrightarrow  H_1(\partial H)$$
and its preimage $[\Sigma(kx)]$ would be primitive. 
Therefore Lemma~\ref{lemnondeg} implies that there would exist $y$ in $H_1(H)$ such that $\langle[\Sigma(kx)],y \rangle_H=1$. Now, this $y$ could be thought of as an element of $H_1(\partial H;\ZZ)$, and $\langle[\Sigma(kx)],y \rangle_H=\langle kx,y \rangle_{\partial H}=k\langle x,y \rangle_{\partial H}=1$. 
Then $k$ could not be larger than $1$.
\eop

Recall that the {\em linking form\/} $\ell_{\QQ/\ZZ}(M)$ \index{NN}{lQoverZ@$\ell_{\QQ/\ZZ}(M)$} of a rational homology sphere $M$ is the bilinear form
$$\ell_{\QQ/\ZZ}(M):H_1(M)^2 \rightarrow \frac{\QQ}{\ZZ}$$
that maps a pair of homology classes to the linking number of two of their representatives.

The form $\ell_{\QQ/\ZZ}(M)$ is {\em non-degenerate\/} in the following sense.
For any element $x$ of $H_1(M)$ of order $k>1$, there exists $y$ such that
$\ell_{\QQ/\ZZ}(x,y)=\frac1k.$ 
This can be seen by applying Lemma~\ref{lemnondeg}
to the exterior of a knot which represents $x$.

%Recall the following lemma.

%\begin{lemma}[\cite{al}]
%\label{lemclasptrans}
%For any two integral homology handlebodies $A$ and $B$ whose boundaries are %identified so that their lagrangians coincide, there exists a degree $k$ %clover
%in $A$ such that $B$ is obtained from $A$ by surgery on $D$.
%\end{lemma}

\begin{proposition}[Massuyeau \cite{mas}]
A real-valued invariant $\kappa$ of rational homology spheres that does not vary under surgery along $Y$-graphs 
%vanishes at rational generalised degree 2 clovers 
only depends on $(H_1(.);\ell_{\QQ/\ZZ}(.))$.
\end{proposition}
\bp
Let $M$ and $M^{\prime}$ be two rational homology spheres. Assume 
that there exists
an isomorphism $\phi: H_1(M^{\prime}) \longrightarrow H_1(M)$ such that
$$\ell_{\QQ/\ZZ}(M^{\prime})=\ell_{\QQ/\ZZ}(M) \circ \phi^2.$$

Then $L^{\prime}=(K^{\prime}_1, \dots, K^{\prime}_r)$ be a family of framed knots in $M^{\prime}$ whose homology classes
generate $H_1(M^{\prime})$.
Then there exists a framed link $L=(K_1, \dots, K_r)$ in $M$ such that the $K_i$
are homologous to the $\phi([K^{\prime}_i])$ with the same linking matrix as 
$L^{\prime}$. (The hypothesis implies that it is true in $\QQ/\ZZ$, perform crossing changes and framing changes so that it is true in $\QQ$.)
Glue a tree $T$ to $L$ along its endpoints so that each component of $L$
has exactly one point identified to an endpoint of $T$.
Then $L$ is part of a connected graph $T \cup L$ whose $H_1$ is freely generated by
the classes of the $K_i$. Then the regular neighborhood $N(T \cup L)$ 
of $L$ in $M$
is a handlebody whose complement is a homology handlebody $H$, thanks to Lemma~\ref{lemcharZhan}, because the map induced by the inclusion maps $H_1(\partial H)$ onto $H_1(H)$.

There is a similar diffeomorphic handlebody $N(T^{\prime} \cup L^{\prime})$ whose 
complement is again a homology handlebody $H^{\prime}$.
Using the diffeomorphism mentioned above, we may write 
$M=N(T \cup L) \cup H$ and  $M^{\prime}=N(T \cup L) \cup H^{\prime}$
where $H$ and $H^{\prime}$ are two handlebodies with the same lagrangian
determined by the linking matrices of $L$ and $L^{\prime}$. 
Recall the following lemma from \cite[Lemma 4.11]{al}.
\begin{lemma}
\label{lemclasptrans}
For any two integral homology handlebodies $A$ and $B$ whose boundaries are identified so that their lagrangians coincide, there exists a degree $k$ clover
in $A$ such that $B$ is obtained from $A$ by surgery on $D$.
\end{lemma}
\eop

According to this lemma, $H$ and $H^{\prime}$ can be obtained from one another by surgeries on
$Y$-graphs. Therefore $M^{\prime}$ can be obtained from $M$ 
by surgeries on $Y$-graphs. Then $\kappa(M)=\kappa(M^{\prime})$.
\eop

Now, we furthermore assume that $\kappa(-M)=-\kappa(M)$.
This implies that $\kappa(S^3)=0$.
Then since $\kappa$ vanishes at rational generalised 2-clover of genus $0$, $\kappa$ is additive under connected sum. 

In order to conclude the proof of Proposition~\ref{propcompwde} and hence of Theorem~\ref{thcompwal}, we are left with the proof of the following proposition.

\begin{proposition}
\label{proplinpair}
A real-valued invariant $\kappa$ of rational homology spheres 
\begin{itemize}
\item that only depends on $(H_1(.);\ell_{\QQ/\ZZ}(.))$
\item that is additive under connected sum, and
\item whose sign changes under orientation reversing,
\end{itemize} 
vanishes identically.
\end{proposition}

Since we now consider $\kappa$ as an invariant of $(H_1(M);\ell_{\QQ/\ZZ}(M))$, 
$\kappa$ is additive under orthogonal sum.
Since $\ell_{\QQ/\ZZ}$ maps a pair of elements with coprime order to zero, $(H_1(M);\ell_{\QQ/\ZZ}(M))$ is the orthogonal sum of its $p$-components
for prime integers $p$ (made of elements $\gamma$ such that there exists $k \in \NN \setminus \{0\}$ such that $p^k{\gamma}=0$).

The classification of linking forms of rational homology spheres has been started by Wall \cite[Theorem 4]{wall} and completed
by Kawauchi and Kojima~\cite{kk}. Part of the results in \cite{wall} recalled below, together
with the behaviour of $\kappa$ under orientation reversing and connected sum,
imply that $\kappa$ vanishes. Details are given below.

\begin{lemma} {\bf (\cite[Theorem 4]{wall}) }
Let $p$ be a prime integer greater than $2$. Assume that the order of any element of $H_1(M)$ is a power of $p$.
Let $n(p)$ be the smallest integer in $\{2, \dots, p-1\}$ that is not a square mod $p$. Then $(H_1(M);\ell_{\QQ/\ZZ}(M))$ is an orthogonal 
sum of modules of the forms
$$[p^k,s=n(p) \;\;\mbox{or} \;\;1]=(\ZZ/p^k\ZZ[x];\ell_{\QQ/\ZZ}(x,x)= \frac{s}{p^k}).$$
Furthermore the orthogonal sum of two copies of $[p^k,n(p)]$
is isomorphic to the orthogonal sum of two copies of $[p^k,1]$.
\end{lemma}
\bp
Let $k$ be the largest integer such that $p^k$ is the order of an element of 
$H_1(M)$. Then there exists $x$ such that $p^k\ell_{\QQ/\ZZ}(M)(x,x) \in \ZZ/p^k\ZZ$ is coprime with $p$. 
(Otherwise, if the order of $x$ is $p^k$, let $y$ be such that $p^k\ell_{\QQ/\ZZ}(M)(x,y)=1\;\mbox{mod} \;p^k$, then  $p^k\ell_{\QQ/\ZZ}(M)(x+y,x+y)=2\;\mbox{mod}\; p$.)
Now, the square from $\left(\ZZ/p^k\ZZ\right)^{\ast}$ to itself is a group morphism. Its image $\left(\left(\ZZ/p^k\ZZ\right)^{\ast}\right)^2$ has index two and is made of the squares. It does not contain $n(p)$.
Thus 
$$\left(\frac{\ZZ}{p^k\ZZ}\right)^{\ast}=\left(\left(\frac{\ZZ}{p^k\ZZ}\right)^{\ast}\right)^2
\cup \left(\left(\frac{\ZZ}{p^k\ZZ}\right)^{\ast}\right)^2 n(p).$$
Therefore, there exists $x \in H_1(M)$ such that $p^k\ell_{\QQ/\ZZ}(M)(x,x)$ equals $1$ or $n(p)$. This allows for a proof that $H_1(M)$ 
decomposes as in the statement by induction on the order of $H_1(M)$, because
the orthogonal of the subspace generated by $x$ is equipped with a non-degenerate linking form.

Now, since $(n(p)-1)=a^2$ in $\ZZ/p^k\ZZ$, in the orthogonal sum $(\ZZ/p^k\ZZ)e \oplus (\ZZ/p^k\ZZ)f$ of two copies of $[p^k,1]$, $p^k\ell_{\QQ/\ZZ}(ae+f,ae+f)=n(p)$. This orthogonal sum is isomorphic 
to the orthogonal sum of the module $[p^k,n(p)]$ generated by $(ae+f)$ and 
of its orthogonal that must be isomorphic to $[p^k,s=1 \;\mbox{or}\; n(p)]$.
Now, because of the effect of a change of generating system, $sn(p)$ must be a square, and therefore $s=n(p)$. 
\eop

Note that conversely, all the modules above are realised as the $H_1$ of connected sums of lens spaces.
In particular $2 \kappa([p^k,1])= 2 \kappa([p^k,n(p)])$. Then $\kappa$ is the same for all rational homology spheres $M$ with $H_1(M)=\ZZ/p^k\ZZ$. In particular,
for such a manifold, $\kappa(M)=\kappa(-M)=-\kappa(M)=0$.
Then, using the behaviour of $\kappa$ under orthogonal sum, $\kappa$ vanishes at all
the rational homology spheres without $2$-torsion in their $H_1$.

\begin{lemma} {\bf \rm (\cite[Theorem 4]{wall})}
For $n \in \{-3,-1,1,3\}$, and for $k \geq 1$, let $A^k(n)$ denote the cyclic group $\ZZ/2^k\ZZ[x]$ equipped with the pairing $\ell$ where $\ell(x,x)=n/2^k$
Assume that the order of any element of $H_1(M)$ is a power of $2$.
\begin{itemize}
\item
Then the orthogonal sum of $H_1(M)$ and of a finite number of copies of modules of the form $A^k(1)$, $k \geq 1$, is an orthogonal sum of modules of the form $A^k(n)$.
\item
For any $k \geq 1$, the orthogonal sum of $2$ copies of $A^k(1)$ is isomorphic 
to the orthogonal sum of $2$ copies of $A^k(-3)$.
\item
For any $k \geq 1$, the orthogonal sum of $4$ copies of $A^k(1)$ is isomorphic 
to the orthogonal sum of $A^k(3)$ and three other modules of the form $A^k(n)$.
\end{itemize}
\end{lemma}
\bp 
When $k \geq 3$, the square from $\left(\ZZ/2^k\ZZ\right)^{\ast}$ to itself is a group morphism whose kernel has the four elements, $1$, $-1$, $1+2^{k-1}$ and  $-1+2^{k-1}$. Therefore its image $\left(\left(\ZZ/2^k\ZZ\right)^{\ast}\right)^2$ has index $4$, since all its elements must be congruent to $1$ mod $8$, it is easy to see that and $\left(\left(\ZZ/2^k\ZZ\right)^{\ast}\right)^2$ is exactly made of the numbers that are congruent to $1$ mod $8$. In particular, any element of $\left(\ZZ/2^k\ZZ\right)^{\ast}$ belongs to $\left(\left(\ZZ/2^k\ZZ\right)^{\ast}\right)^2 \{-3,-1,1,3\}$.

Again, we try to diagonalise the linking form by induction on the order of $H_1(M)$. Let $k$ be the largest integer such that $2^k$ is the order of an element of 
$H_1(M)$. If there exists $x$ such that $2^k\ell_{\QQ/\ZZ}(M)(x,x)$ is odd,
then  $\ZZ/2^k\ZZ[x]$ is of the form $A^k(n)$,
and its orthogonal of is a module of lower order that satisfies the hypotheses.
Otherwise, for any $x$ of order $2^k$, $2^k\ell_{\QQ/\ZZ}(M)(x,x) \in \ZZ/2\ZZ^k$ is even.
However, for such an $x$, there exists $y$ such that $2^k\ell_{\QQ/\ZZ}(M)(x,y)=1$. Consider the orthogonal sum of $H_1(M)$ with 
$A^k(1)=\ZZ/2^k\ZZ[z]$ where $2^k\ell_{\QQ/\ZZ}(M)(z,z)=1$.
Then $H_1(M) \oplus A^k(1)$ splits as the orthogonal sum of the module generated by the two orthogonal elements $(x+z)$ and $(y-z)$, that is an orthogonal sum of modules of the form $A^k(n)$, and of its orthogonal whose order is lower than the order of $H_1(M)$.  The first assertion follows by induction on the order of $H_1(M)$.

The orthogonal sum 
$$\frac{\ZZ}{2^k\ZZ} a \oplus \frac{\ZZ}{2^k\ZZ} b$$
where the mentioned generators have self-linking $\frac{1}{2^k}$,
is isomorphic to the orthogonal sum $\frac{\ZZ}{2^k\ZZ} (a+2b) \oplus \frac{\ZZ}{2^k\ZZ} (b-2a)$ of two modules isomorphic to $A^k(-3)$.

The orthogonal sum 
$$\frac{\ZZ}{2^k\ZZ} a \oplus \frac{\ZZ}{2^k\ZZ} b \oplus \frac{\ZZ}{2^k\ZZ} c$$
where the mentioned generators have self-linking $\frac{1}{2^k}$
splits as the direct orthogonal sum of $A^k(3)$ that is generated by $(a+b+c)$ and its orthogonal that is  $\frac{\ZZ}{2^k\ZZ} x \oplus \frac{\ZZ}{2^k\ZZ} y$
equipped with a non-degenerate linking form where $2^k\ell_{\QQ/\ZZ}(M)(x,x)$ is even. As above, the orthogonal sum of $A^k(1)=\ZZ/2^k\ZZ[z]$ and $\left(\frac{\ZZ}{2^k\ZZ} (a+b+c)\right)^{\perp}$ is the orthogonal sum of the
orthogonal sum of two modules of the form $A^k(n)$ and its orthogonal that is 
necessarily cyclic and equipped with a non degenerate linking form, and that  is therefore of the form $A^k(n)$.
\eop

By the behaviour of $\kappa$ under orientation reversing we deduce that $\kappa(A^k(1))=-\kappa(A^k(-1))$ and that $\kappa(A^k(3))=-\kappa(A^k(-3))$.
By the second point of the above lemma, $2\kappa(A^k(1))=2\kappa(A^k(-3))$. 
$$\kappa(A^k(1))=\kappa(A^k(-3))=-\kappa(A^k(-1))=-\kappa(A^k(3)).$$
Then by the third point of the above lemma, $4\kappa(A^k(1))=r \kappa(A^k(1))$, where $r$ is an integral number less or equal than $2$.
Therefore $\kappa(A^k(1))=0$ and $\kappa(A^k(n))=0$, for all $n$ in $\{-3,-1,1,3\}$.
The first point of the lemma together with the additivity of $\kappa$ under orthogonal sum allows us to conclude that $\kappa$ vanishes at the manifolds $M$ such that the order of $H_1(M)$ is a power of two, and next that $\kappa(M)=0$ for any rational homology sphere $M$.

This ends the proof of Proposition~\ref{proplinpair} and hence the proof of Theorem~\ref{thcompwal}.

\subsection{Alternative definitions of the Casson-Walker invariant}

I thank Misha Polyak and Oleg Viro for encouraging me to write this last subsection.

As a corollary to Theorem~\ref{thcompwal}, we get
\begin{corollary}
For any rational homology sphere $M$, 
and for any trivialisation $\tau_M$ of 
$T(M \setminus \infty)$ that is standard near $\infty$,
$$\lambda(M)=\frac{\int_{C_2(M)}\omega(\tau_M)^3}{6} - \frac{p_1(\tau_M)}{24}.$$
%%%SIGNE
\end{corollary}
\bp 
According to Theorem~\ref{thcompwal}, $Z_1=\lambda /2 [\tata]$.
Then by definition of $Z_1$ \rconstthkktfra, and by the proof of \inconstpropxiun, $$Z_1(M)= \frac{\int_{C_2(M)}\omega(\tau_M)^3}{12}[\tata]+\frac{p_1(\tau_M)}{4}\xi_1.$$ 
%%%SIGNE
According to \inconstpropxiun, $12\xi_1=-[\tata]$. 
\eop

The integral $\int_{C_2(M)}\omega(\tau_M)^3$ can be rewritten in various different ways as follows.
Recall that a trivialisation $\tau_M$ of 
$T(M \setminus \infty)$ that is standard near $\infty$ induces a map 
$p_M(\tau_M)$ from $\partial C_2(M)$ to $S^2$. See \rconstsubfundform.

\begin{lemma}
Let $M$ be a rational homology sphere. Let $\tau_M$ be a trivialisation of 
$T(M \setminus \infty)$ that is standard near $\infty$.
Let $\omega_a(S^2)$, $\omega_b(S^2)$ and $\omega_c(S^2)$ be $3$ two-forms on $S^2$ whose integrals over $S^2$ equal $1$.
Let $\omega_a$, $\omega_b$ and $\omega_c$ be $3$ closed two-forms on $C_2(M)$ that coincide with $p_M(\tau_M)^{\ast}(\omega_{a}(S^2))$, $p_M(\tau_M)^{\ast}(\omega_{b}(S^2))$ and $p_M(\tau_M)^{\ast}(\omega_{c}(S^2))$ on $\partial C_2(M)$, respectively. 
Then $$\int_{C_2(M)}\omega(\tau_M)^3=\int_{C_2(M)}\omega_a \wedge \omega_b \wedge \omega_c.$$
\end{lemma}
\bp
The arguments are already in \const. However, since only antisymmetric forms were considered in \const, and since the proof is far quicker in this case, we give it.
It is enough to prove that the right-hand side does not depend on the choices of $\omega_a$, $\omega_b$ and $\omega_c$. By symmetry, it is enough to prove  
that it is independent of $\omega_a$. By \inconstlemnolosseta, it is enough
to prove that $\int_{C_2(M)}\omega_a \wedge \omega_b \wedge \omega_c$ does not change when ${\rm d} \eta$ is added to $\omega_a$, for some one-form $\eta$ on $C_2(M)$ that reads $p_M(\tau_M)^{\ast}(\eta_{S^2})$ on $\partial C_2(M)$ for some one-form $\eta_{S^2}$ on $S^2$.
By the Stokes theorem, the variation of $\int_{C_2(M)}\omega_a \wedge \omega_b \wedge \omega_c$ under the addition of ${\rm d} \eta$
reads
$$\int_{\partial C_2(M)}\eta \wedge \omega_b \wedge \omega_c=\int_{\partial C_2(M)}p_M(\tau_M)^{\ast}\left(\eta_{S^2} \wedge \omega_b(S^2)
\wedge \omega_c(S^2)\right)$$
and vanishes since $5$-forms vanish on $S^2$.
\eop

Let $d$ be a point in $S^2$. Consider the codimension $2$ submanifold
$p_M(\tau_M)^{-1}(d)$ of $\partial C_2(M)$. Since $H_3(C_2(M);\QQ)=0$, 
$p_M(\tau_M)^{-1}(d)$ is the boundary of a rational $4$-chain $\Sigma_d$.
If $M$ is a $\ZZ$-sphere, we may even assume that $\Sigma_d$ is a codimension $2$ submanifold of $C_2(M)$.

\begin{lemma}
Let $M$ be a rational homology sphere. Let $\tau_M$ be a trivialisation of 
$T(M \setminus \infty)$ that is standard near $\infty$.
Let $a$, $b$ and $c$ be three distinct points in $S^2$.
Let $\Sigma_a$, $\Sigma_b$ and $\Sigma_c$ be three rational $4$-chains in $C_2(M)$ with respective boundaries $p_M(\tau_M)^{-1}(a)$, $p_M(\tau_M)^{-1}(b)$ and $p_M(\tau_M)^{-1}(c)$.
Then $\int_{C_2(M)}\omega(\tau_M)^3$ is the algebraic intersection of $\Sigma_a$, $\Sigma_b$ and $\Sigma_c$ in $C_2(M)$.
\end{lemma}
\bp 
First note that the algebraic intersection of $\Sigma_a$, $\Sigma_b$ and $\Sigma_c$ is well-defined and independent of the involved choices because the boundaries are three disjoint fixed submanifolds of the boundary of $C_2(M)$ and because $H_4(C_2(M);\QQ)=0$. Thus, we must now compare two topological invariants of $(M;\tau_M)$. 

Without loss assume that $\Sigma_a$, $\Sigma_b$ and $\Sigma_c$ are in general position so that they intersect in a finite number of points that have
a neighborhood of the form $\RR^2_a \times \RR^2_b \times \RR^2_c$
(for three copies $\RR^2_a$, $ \RR^2_b$, $\RR^2_c$ of $\RR^2$) 
that is intersected by the support of $\Sigma_a$ along $0 \times \RR^2_b \times \RR^2_c$, by the support of $\Sigma_b$ along $\RR^2_a \times 0 \times \RR^2_c$,
and by the support of $\Sigma_c$ along $\RR^2_a \times \RR^2_b \times 0$.

Now, use the preceding lemma with three forms $\omega_a(S^2)$, $\omega_b(S^2)$ and $\omega_c(S^2)$ with disjoint supports concentrated near
$a$, $b$ and $c$ respectively, and with three forms $\omega_d$, for $d \in \{a,b,c\}$, such that:
\begin{itemize}
\item $\omega_d$ is supported in a small neighborhood $N(\Sigma_d)$ of $\Sigma_d$,
\item for any oriented dimension $2$-manifold $F$ whose boundary is disjoint
from $N(\Sigma_d)$, $\int_F \omega_d$ is the algebraic intersection of $F$ and $\Sigma_d$, and
\item in an affine neighborhood of an intersection point as above, $\omega_d$ vanishes along any two-dimensional plane that is not transverse to the affine 4-dimensional support of $\Sigma_d$.
\end{itemize}
%$\omega_d$ integrate as the  on any dimension $2$-manifold 
%whose boundary is in $\partial C_2(M)$ . 
\eop

Let $M$ be a rational homology sphere. Let $p_1(M)$ denote the image of the map
$p_1$ from the set of trivialisations of 
$T(M \setminus \infty)$ that are standard near $\infty$ to $\ZZ$.
If $j \in p_1(M)$, let $\tau_j(M)$ be a trivialisation of 
$T(M \setminus \infty)$ that is standard near $\infty$ such that $p_1(\tau_j(M))=j$.

If $0 \in p_1(M)$, (according to \inconstproppont, this is always the case when $M$ is a $\ZZ$-sphere), then $$\lambda(M)=\frac{\int_{C_2(M)}\omega(\tau_0(M))^3}{6}$$
and the previous lemma allows us to express $6\lambda(M)$ as the algebraic intersection of three rational $4$-chains in $C_2(M)$, (or even of three
codimension 2 submanifolds if $M$ is a $\ZZ$-sphere).

Otherwise, \inconstproppont ~ensures that $p_1(M)$ contains a subset of the form
$\{i-4, i\}$.
Then $$\lambda(M)=\frac{i}{24}\int_{C_2(M)}\omega(\tau_{i-4}(M))^3
+\frac{4-i}{24}\int_{C_2(M)}\omega(\tau_{i}(M))^3.$$
%at least one of the following subsets $\{0\}$, $\{-2,2\}$, $\{-3,1\}$ or %$\{-1,3\}$ 
%. 
In this case, the previous lemma allows us to express $24\lambda(M)$  as the sum of two  rational
algebraic intersection numbers in $C_2(M)$.

\newpage

\begin{center}{\bf \large Terminology}\end{center}
\addcontentsline{toc}{section}{Terminology}
\begin{theindex}

  \item clover, 48, 49

  \indexspace

  \item Euler number, 14

  \indexspace

  \item finite type invariant, 4
  \item form
    \subitem admissible, 8
    \subitem fundamental, 8

  \indexspace

  \item half-edge, 3
  \item homology sphere, 4

  \indexspace

  \item integral generalised clover, 5
  \item integral homology handlebody, 5

  \indexspace

  \item Jacobi diagram, 3
    \subitem automorphism of, 3
    \subitem orientation of, 4

  \indexspace

  \item Lagrangian, 5

  \indexspace

  \item $\QQ$-handlebody, 4
  \item $\QQ$-sphere, 4

  \indexspace

  \item rational generalised clover, 5
    \subitem degree of, 5
  \item rational homology handlebody, 4
  \item rational homology sphere, 4

  \indexspace

  \item special admissible forms, 17

  \indexspace

  \item Torelli homeomorphism, 9
  \item twist of a trivialisation, 14

  \indexspace

  \item weight system, 4

  \indexspace

  \item Y-graph, 48
    \subitem leaf of, 48

  \indexspace

  \item $\ZZ$-handlebody, 5
  \item $\ZZ$-sphere, 4

\end{theindex}

\begin{center}{\bf \large Notation}\end{center}
\addcontentsline{toc}{section}{Notation}

\begin{theindex}
\begin{multicols}{3}
  \item $A^i_t$, 29
  \item $\CA_n(\emptyset)$, 4
  \item AS, 4
  \item $A(t,u;s)$, 35
  \item $\sharp \mbox{Aut}(\Gamma)$, 4

  \indexspace

  \item $B^D_M$, 20

  \indexspace

  \item $C(a)$, 34, 35
  \item $c_A$, 20
  \item $c_B$, 20
  \item $\chi(S)$, 14
  \item $\chi(TS; .)$, 14, 19
  \item $C_I^i$, 9
  \item $C_2(X)$, 10

  \indexspace

  \item $D(\omega_0(M_I))$, 30

  \indexspace

  \item $E(\Gamma)$, 3, 6
  \item $e(S_0(a))$, 35
  \item $\eta(a^i_j)$, 10, 29
  \item $\eta(b^i_j)$, 30
  \item $\eta_{[-1,1]}$, 9

  \indexspace

  \item $F$, 18
  \item $F(a)$, 34, 35
  \item ${F}(c,\tau_b)$, 18
  \item $F_U$, 20

  \indexspace

  \item $G(a)$, 38

  \indexspace

  \item $H_g$, 5
  \item $H(\Gamma)$, 3, 6

  \indexspace

  \item $\CI(A^i,B^i)$, 6
  \item IHX, 4
  \item $\langle,\rangle_{\Sigma}$, 5

  \indexspace

  \item $\ell$, 10
  \item ${\cal L}_A$, 5
  \item $\lambda$, 47
  \item $\lambda_W$, 7, 47
  \item $\CL_A^{\ZZ/2\ZZ}$, 8, 15, 16
  \item $\ell(D;\Gamma)$, 6, 7, 10
  \item $\ell(D;\Gamma;\sigma)$, 6, 11
  \item $\ell_{\QQ/\ZZ}(M)$, 52

  \indexspace

  \item $M_I$, 8
  \item $M_J(D)$, 5

  \indexspace

  \item $\omega(c,\tau_b)$, 18
  \item $\omega(c;\tau,\tau_b)$, 17--19
  \item $\omega_j$, 8
  \item $\omega_M$, 30
  \item $\omega(M_I)$, 9, 10, 29, 31
  \item $\omega(M_i)$, 31
  \item $\omega(p^i)$, 29
  \item $\omega_{{\cal Q}}(c_A,\tau_M )$, 21
  \item $\omega_{{\cal Q}}(c_B,\tau_M )$, 21
  \item $\omega_{S^2}$, 8, 13
  \item $\omega(\tau_M)$, 8

  \indexspace

  \item $p(i)$, 9
  \item $p(a^i_j)$, 34
  \item $[p^i,\infty(v)]$, 29
  \item $p_1(c;\tau,\tau_b)$, 17
  \item $p_1(\tau_I^{\CC})$, 9

  \indexspace

  \item ${\cal Q}$, 21

  \indexspace

  \item $R_{\theta}$, 13

  \indexspace

  \item $S(b^i_j)$, 30
  \item $\Sigma$, 20
  \item $S(n)$, 14
  \item $S_2$, 35
  \item $S_{2,-4}$, 35
  \item $S_0(a)$, 35

  \indexspace

  \item $\tau$, 20
  \item $\tau_b$, 17, 20
  \item $\tau_{\CC}(c;\tau,\tau_b)$, 17, 20
  \item $\tau_j$, 8
  \item $\tau_j^{\CC}$, 8
  \item $\tau_M$, 8, 20
  \item $\CT_k$, 13
  \item $\CT_c$, 14
  \item $\theta$, 13
  \item $\theta(c)$, 13
  \item $T(x,y)$, 35

  \indexspace

  \item $V(\Gamma)$, 3, 6
  \item $v_i$, 13

  \indexspace

  \item $Z$, 3
  \item $Z_{KKT}$, 3
  \item $Z_{LMO}$, 3
  \item $Z(M)$, 17
  \item $Z(M_I)$, 9
  \item $z_n(c;\tau,\tau_b)$, 17, 18
  \item $Z_{n}(D)$, 7
  \item $Z(\omega(M_I))$, 9
\end{multicols}
\end{theindex}

%%\printindex{TT}{Terminology}
%%\printindex{NN}{Notation}
%\begin{multicols}{3}\end{multicols}

\end{document}